%% file: 0paper.tex
\documentclass[12pt]{article}
\usepackage{epsfig}
\title{Unbounded Orbits for Outer Billiards}
\author{Richard Evan Schwartz \thanks{\hskip 5 pt 
This research is supported by 
N.S.F. Grant DMS-0604426}}

\newtheorem{theorem}{Theorem}[section]

\newtheorem{lemma}[theorem]{Lemma}

\newtheorem{corollary}[theorem]{Corollary}

\def\startproof{{\bf {\medskip}{\noindent}Proof: }}

\def\endproof{$\spadesuit$  \newline}

\def\N{\mbox{\boldmath{$N$}}}%
\def\P{\mbox{\boldmath{$P$}}}%
\def\Q{\mbox{\boldmath{$Q$}}}%
\def\R{\mbox{\boldmath{$R$}}}%
\def\T{\mbox{\boldmath{$T$}}}%
\def\Z{\mbox{\boldmath{$Z$}}}%

\begin{document}

\maketitle

\begin{abstract}
The question of B.H. Neumann, which dates back to the $1950$s, asks
if there exists an outer billiards system with an unbounded orbit.
We prove that
outer billiards for the Penrose kite, the convex quadrilateral
from the Penrose tiling, has an unbounded orbit.  We
also analyze some finer properties of the orbit structure,
and in particular produce an uncountable family of unbounded orbits.
Our methods
relate  outer billiards on the Penrose kite to polygon
exchange maps, arithmetic dynamics, and self-similar
tilings.
\end{abstract}

\input{1intro}

\input{2broad}

\input{3shadow}

\input{4verify}

\input{5backwards}

\input{6return}

\input{7routines}

\input{8appendix}
\input{refs}

\end{document}

%% file: 1intro.tex
\section{Introduction}

\subsection{History of the Problem}

{\it Outer billiards\/} is a basic dymamical system which
serves as a toy model for
celestial mechanics. See Sergei Tabachnikov's book
[{\bf T\/}], and also the survey  $[{\bf DT\/}]$, for an
exposition of the subject and many references.

To define an outer billiards system,
one starts with a bounded convex set $S \subset \R^2$
and considers a point $x_0 \in \R^2-S$.  
One defines $x_1$ to be the point such that
the segment $\overline{x_0x_1}$ is tangent to $S$ at its
midpoint and $S$ lies to the right of the ray
$\overrightarrow{x_0x_1}$.  (See Figure 1.1 below.)
The point $x_1$ is not well-defined if
$\overline{x_0x_1}$ is tangent to $S$ along a segment.
This will inevitably happen sometimes when $S$ is a polygon.
Nonetheless, the outer billiards
construction is almost everywhere well defined.
The iteration $x_0 \to x_1 \to x_2...$
is called the {\it forwards outer billiards orbit\/} of $x_0$.
The backwards orbit is defined similarly.

\begin{center}
\psfig{file=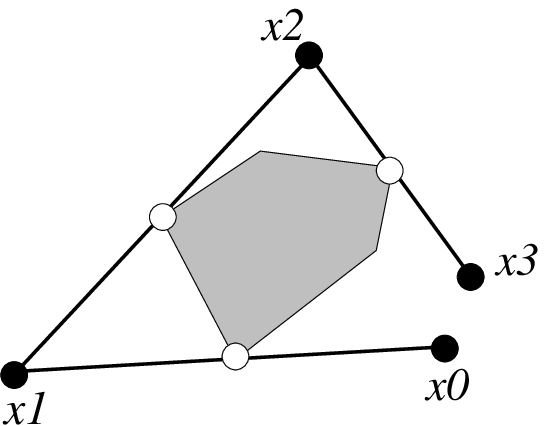}
\newline
{\bf Figure 1.1:\/} Outer Billiards
\end{center}

B.H. Neumann 
\footnote{My information on this comes from 1999 email
correspondence between Bernhard Neumann and Keith Burns,
and also from Berhnard's son Walter.}
introduced outer billiards during some lectures
for popular audiences given in the 1950s.
See, e.g. [{\bf N\/}].  J. Moser popularized the construction
in the 1970s. Moser [{\bf M\/}, p. 11]
attributes the
following question to Neumann {\it circa\/} 1960, though it is sometimes
called {\it Moser's Question\/}.
\newline
\newline
{\bf Question:\/} {\it Is there an outer billiards system with an unbounded
orbit? \/}
\newline

There have been several results related to this question.
\begin{itemize}
\item Moser [{\bf M\/}] sketches a proof, inspired by
K.A.M. theory, that outer
billiards on $S$ has all bounded orbits provided
that $\partial S$ is at least $C^6$ smooth and positively curved.
R. Douady gives a complete proof in his thesis, [{\bf D\/}].
See [{\bf B\/}] for related work.
\item In [{\bf VS\/}], [{\bf Ko\/}], and (later, but with different methods)
[{\bf GS\/}], it is proved
that outer billiards on a {\it quasirational polygon\/} has
all orbits bounded.  This class of polygons
includes polygons with rational vertices and
also regular polygons. In the rational case,
all defined orbits are periodic.
\item  Tabachnikov analyzes
the outer billiards system for the regular pentagon
and shows that there are some non-periodic (but bounded)
orbits.  See [{\bf T\/}, p 158] and the references there. 
\item Genin [{\bf G\/}] shows that all orbits are bounded
for the outer billiards systems associated to trapezoids.
He also makes a brief numerical study of
a particular irrational kite based on the
square root of $2$, observes possibly
unbounded orbits, and indeed conjectures that this is
the case.
\end{itemize}

\subsection{The Main Result}

The main goal of this paper is to prove:

\begin{theorem}
\label{main}
Outer billiards on the Penrose kite has an unbounded orbit.
In fact, both the forwards and backwards orbits of the point $p$ shown in Figure 1.2
are entirely defined and unbounded.
\end{theorem}

\begin{center}
\psfig{file=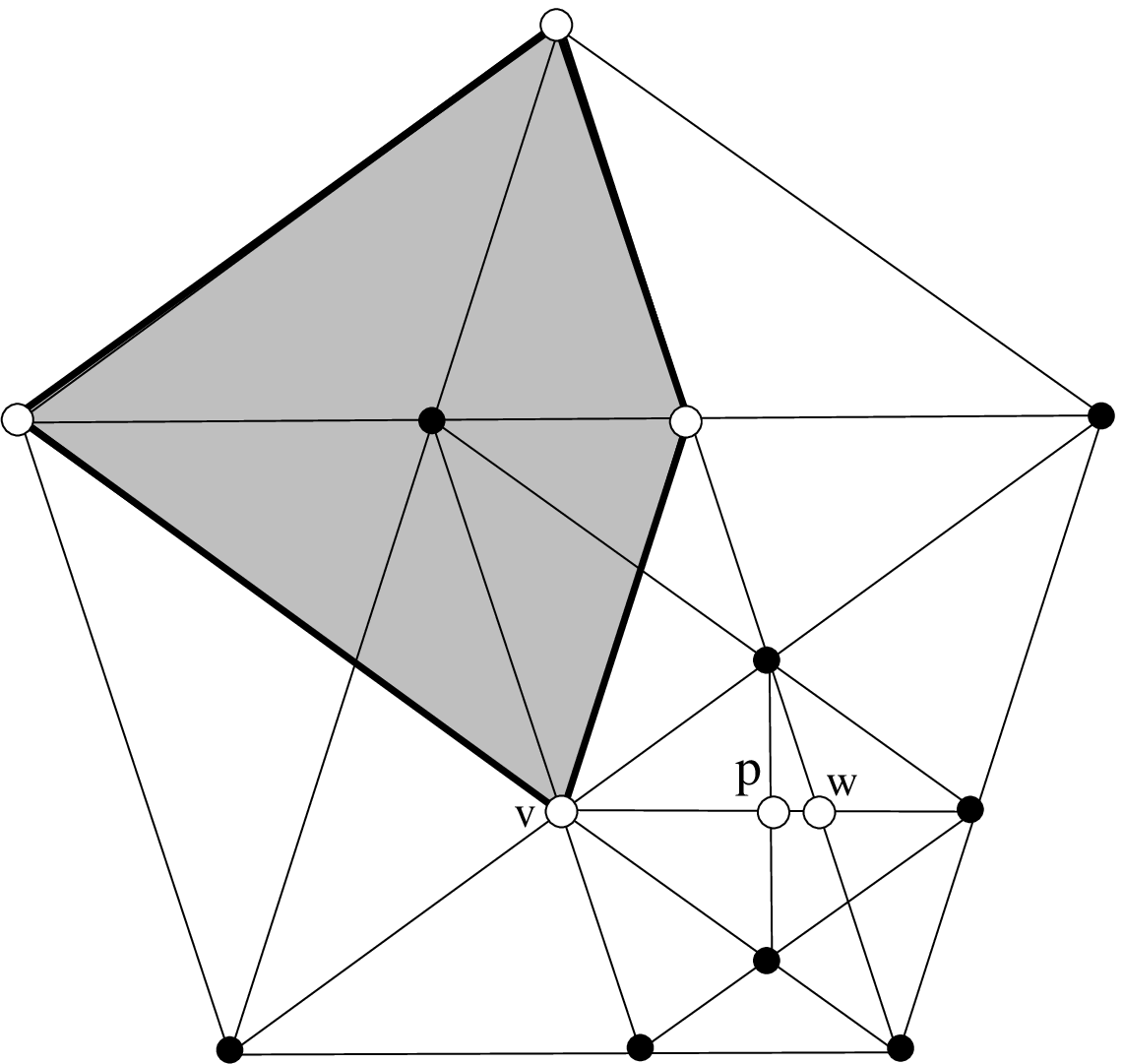}
\newline
{\bf Figure 1.2:\/} The Penrose Kite
\end{center}

The {\it Penrose kite\/} is the convex quadrilateral that
appears in the Penrose tiling. (In general, a {\it kite\/} is a
convex quadrilateral with bilateral symmetry.)  Figure 1.2 shows a classic
construction of the Penrose kite--the shaded figure--based
on a regular pentagon.  The additional lines show how
the point $p$ is constructed.
The significance of the points $v$ and $w$ will be
explained in Theorem \ref{aux} below.  In
\S 2.1 we will give more traditional coordinates for the
objects of interest to us.

\subsection{Outline of the Proof}

We will give a rigorous computer-assisted proof of
Theorem \ref{main}.  We think that our proof gets close
to the conceptual core of what is going on, but
we need a lot of computation to make it go.
Here are the main ideas. Let $\Upsilon$ denote the
square of the outer billiards map.

\begin{enumerate}
\item
We replace the Penrose kite by an affinely
equivalent kite $S$ whose vertices lie in the
ring $\Z[\phi]$, where $\phi$ is the golden ratio.
We show that the outer billiards map (and in particular
$\Upsilon$) is entirely defined on a
set ${\cal C\/}(\pm)$, consisting of certain
points whose first coordinates are positive and lie in
$\frac{1}{2}\Z[\phi]$ and whose second coordinates are $\pm 1$.
See Figure 2.1.  Let
$\Upsilon_R: {\cal C\/}(\pm) \to
{\cal C\/}(\pm)$ be the first return map of $\Upsilon$.
\item   Forgetting about the $y$-coordinate (which is always $\pm 1$)
we identify  ${\cal C\/}(\pm)$ with
the set of vertices in $\Z^2 \cap H$, where $H \subset \R^2$
is a certain halfplane.  We then encode the dynamics
of $\Upsilon_R$ by a graph $\Gamma \subset H $.
The vertices of $\Gamma$ lie in $\Z^2 \cap H$ and the edges
are short lattice vectors.
We call $\Gamma$ the {\it arithmetic graph\/}.  See
Figures 2.2 and 2.3 for pictures.
\item The structure of $\Gamma$ is
controlled by a certain partition $\cal P$ of the
square torus $\T^2$ into $26$ convex polygons, and an associated
dynamical system akin to a {\it polygon exchange map\/}.
See Figure 2.5.
This structural result, which we call the
{\it Arithmetic Graph Lemma\/}, is the central technical result
in the paper.
\item
The dynamical process described in Item 3 is compatible with
a kind of multi-valued contraction, defined on
various simply connected subsets of $\T^2$.
The graph of this correspondence is an irrationally
embedded complex line sitting inside $\T^2 \times \T^2$.
Rick Kenyon tells me that this picture is essentially
the same as what one sees in
the {\it cut and project\/} method from
the theory of quasi-crystals.
\item The compatibility discussed in Item 4
forces $\Gamma$ to behave like an aperiodic tiling
with an inflation rule.  Compare [{\bf Ke\/}]. That is,
when we dilate $\Gamma$ about the origin by $\phi^3$
we find that the dilated image is very closely
shadowed by the original image.  This structure
forces the connected
component $\Gamma_0$ of $\Gamma$ containing the
origin to travel unboundedly far away from $\partial H$, and this
translates into the unboundedness of the orbit of our point $p$
shown in Figure 1.2.
\end{enumerate}

\subsection{Orbit Structure}

Our proof of Theorem \ref{main} gives us information
about the fine orbit structure of outer billiards on
the Penrose kite. 
We will see (Lemma \ref{proper}) that the Euclidean norm of
the $n$th point of $O_+(p)$ is a proper function of $n$.
Here $O_+(p)$ is the forward orbit of $p$.
(A function $f: \N \to \N$ is {\it proper\/} if
$f^{-1}[1,n]$ is bounded for all $n \in \N$.)
In contrast,
the backwards orbit $O_-(p)$  returns densely
to a certain Cantor set.

Given a horizontal interval $I \subset \R^2$ and some $\lambda \in (0,1)$ let
$C(I,\lambda) \subset I$ denote the (Cantor set)
limit set  of the semigroup
generated by the two contractions of strength $\lambda$ that
map $I$ into itself and fix an endpoint of $I$.  We say that a
{\it gap\/} point of $C=C(\lambda,I)$ is an endpoint of any component of
$I-C$.  We let $C^*$ denote $C$ minus its countably many
gap points.

Let $\gamma_0$ (respectively $\gamma_1$) be the similarity of
strength $-\lambda$ (reversing orientation on $I$)
 that maps $C(\lambda,I)$ to its own right (respectively left) half.
Let $\Z_2$ denote the $2$-adic integers.
There is a unique homeomorphism $\theta_2: \Z_2 \to C$ that
conjugates
$\gamma_j$ to the map $x \to 2x+j$.  Let $p,v,w$ be as
in Figure 1.2. 
We have set things up so that $\theta_2(0)=p$
when $I=[v,w]$ and $\lambda=\phi^{-3}$.

\begin{theorem}
\label{aux}
Let $C=C([v,w],\phi^{-3})$, where $v$ and $w$ are as in
Figure 1.2 and $\phi$ is the golden ratio. 
The forwards and backwards orbits of every point of
$C^*$ are defined and unbounded.
\begin{itemize}
\item If $x \in C^*-\theta_2(0)$ then the $\Upsilon$-forwards
orbit of $x$ first returns to $C^*$ as $x_-$, such that
$\theta_2^{-1}(x_-)=\theta_2^{-1}(x)-1$.
The number of iterates between $x$ and $x_-$, as well as the maximum
Euclidean norm of such an iterate, is a proper function of
$1/\|x-\theta_2(0)\|$.
\item If $x \in C^*-\theta_2(-1)$ then the $\Upsilon$-backwards
orbit of $x$ first returns to $C^*$ as $x_+$, such that
$\theta_2^{-1}(x_+)=\theta_2^{-1}(x)+1$.
The number of iterates between $x$ and $x_+$, as well as the maximum
Euclidean norm of such an iterate, is a proper function of
$1/\|x-\theta_2(-1)\|$.
\end{itemize}
\end{theorem}

In the two items, we mean to consider just the even iterates in
the orbits.  Again, $\Upsilon$ is the square of the outer billiards map.
Theorem \ref{aux} immediately implies that there are uncountably
many unbounded orbits of the outer billiards map on the Penrose kite
with the following wild behavior:  If $V$ is any neighborhood of the
vertex $v$ and $U$ is any neighborhood of $\infty$, then both the
forwards and backwards orbits oscillate between $U$ and $V$ infinitely
often.

\subsection{Discussion}

Here are some remarks on the origins of Theorem \ref{main}.
Let $S'$ be the Penrose kite. It appears that
$\R^2-S'=A \cup B \cup C$ where
\begin{itemize}
\item $A$ is a dynamically invariant union of finite sided polygons.
All the orbits in a polygon are periodic with the same combinatorial type.
\item $B$ is a countable union of line segments consisting
of the points on which the outer billiard map is undefined.
$B$ is the so-called {\it discontinuity set\/}.
\item $C$ is a fractal set consisting entirely of
unbounded orbits.
\end{itemize}

Tabachnikov develops a similar picture for the regular
pentagon, except that his region $C$ consists of non-periodic but bounded
orbits.  The general regular $N$-gon seems to have a similar kind of
structure, though nobody has yet made a detailed study.
Tabachnikov has a beautiful picture of the case $N=7$
on the cover of his book [{\bf T\/}].
Inspired by Tabachnikov's picture, I designed a computer program,
{\it Billiard King\/}, \footnote{One can download
this Java program from my website: www.math.brown.edu/$\sim$res}
which draws these special sets for kites and explores their structure.

When the
kite has rational vertices, $C$ is empty and
$A$ is a locally finite tiling.  
I searched \footnote{Had I known earlier about Genin's numerical
study, mentioned above, I perhaps would have saved time.}
through the parameter space of rational
kites for examples 
where $A$ featured
small tiles having widely ranging orbits.  Eventually
I considered kites having vertices with Fibonacci
number coordinates, and this led to the Penrose kite.

I hope to find purely conceptual proofs for the results
in this paper, but my first goal is just to establish
them as true statements.
Many people do not like computer-aided proofs, so I would
like to say several things in favor of the proof I have
given.  First of all, it is the
best I could do and there are no other proofs
for results like these.
Second, the proof is not {\it just\/} a 
calculation$-$there are plenty of concepts in it.
Finally,  I was able to check essentially all 
the details in the proof using Billiard King, so
this proof has a lot of computational safeguards that
one might not see in a traditional proof.

I tried as hard as possible to write a proof that is
independent from Billiard King, but still the reader
of this paper would get a much greater appreciation
for what is going on by learning to use
Billiard King.  A few minutes playing with 
Billiard King is worth hours of reading the paper.

\subsection{Organization of the Paper}

\noindent
{\bf Part 1:\/} (\S 2-5)
In this part we prove all our main results modulo the
Arithmetic Graph Lemma. (See Item 3 of \S 1.3.)
In \S 2 we reduce Theorem \ref{main} 
to the statement that
$\Gamma_0$, the component of the arithmetic
graph of interest to us, admits what we call an
{\it inflation structure\/}.  Essentially
this means that certain smallish subsets
(which we call {\it genes\/})
of $\Gamma_0$ are, when dilated by
$\phi^3$, closely
shadowed by other subsets of $\Gamma_0$.
In \S 3 we tie the
existence of an inflation structure to the dynamical
behavior of the torus partition guaranteed by the
Arithmetic Graph Lemma, and thereby reduce the problem to
a finite calculation with integer arithmetic. 
In \S 4 we explain how we perform the calculation.
In \S 5, we prove Theorem \ref{aux} by taking a close look at
the specifics of our inflation structure.
\newline
\newline
{\bf Part 2:\/} (\S 6) 
In this part, we prove the Arithmetic Graph Lemma.
To this end, we factor the return map
$\Upsilon_R$ (from Item 2 of \S 1.3) as
the product of $8$ simpler maps, which we call
{\it strip maps\/}.  (Actually, the decomposition
involves a $9$th map as well.
See Equation \ref{mu} for a precise statement.)
We then show how each strip map actually
is compatible with a certain embedding of
$\R^2-S$ as a dense subspace of the $4$-torus.
Each strip map extends to act as a piecewise
affine transformation of the $4$-torus.
We then embed the dynamics of
$\Upsilon_R$ as a
$2$-dimensional geodesic slice of this $4$-dimensional
dynamical system.
The Arithmetic Graph Lemma
is a consequence of this structure
and some additional integer arithmetic calculations.
\newline
\newline
{\bf Part 3:\/} (\S 7-8)
In \S 7 we describe our calculations using pesudo-code,
in enough detail that the interested reader should be
able to reproduce them.
In \S 8 we include the data used in our calculations.
All the calculations are implemented in Billiard King,
a Java program.
Billiard King is a massive program, but we have
placed the smallish number of routines actually relevant
to the proof in separate files so that they are 
easier to survey. 

\subsection{Acknowledgements}

I thank Jeff Brock, Keith Burns, Peter Doyle, David Dumas,
Eugene Gutkin, Pat Hooper, Richard Kent, Rick Kenyon,
Curt McMullen, and Steve Miller for
helpful and interesting discussions related to this work.
I especially thank Sergei Tabachnikov for all
his help and encouragement.

\newpage

%% file: 2broad.tex
\section{The Proof in Broad Strokes}

\subsection{An Affinely Equivalent Shape}

The Penrose kite is affinely equivalent to the
quadrilateral with coordinates
\begin{equation}
\label{basicshape}
(0,1); \hskip 10 pt 
(-1,0); \hskip 10 pt
(0,-1); \hskip 10 pt
(A,0); \hskip 30 pt A=\phi^{-3}
\end{equation}
where $\phi$ is the golden ratio.
The affine equivalence maps the point $p$ in Figure 1.2 to the point
\begin{equation}
\label{basicpoint}
(\phi^{-2},-1).
\end{equation}

To clarify the logic of our argument, it
is useful to consider a general value of
$A \in (0,1)$ for the moment.
We define 
$T: \Z^2 \to \R$ via the formula
\begin{equation}
\label{map}
T(x,y)=2Ax+2y+\frac{1-A}{2}.
\end{equation}
The point in Equation \ref{basicpoint} is $(T(0,0),-1)$ when
$A=\phi^{-3}$.

\begin{lemma}
\label{basic existence}
\label{exist}
The outer billiards orbit of any point in the set
$$T(\Z^2) \times \{\pm 1\}$$
is entirely defined, both forwards and backwards.
In particular, the entire orbit of 
the point in Equation \ref{basicpoint} is
defined.
\end{lemma}

\startproof
Let $L$ denote the family of horizontal lines in $\R^2$ whose
$y$-coordinates are odd integers. 
The outer billiards map preserves $L$.  The only points
where the first iterate of
the outer billiards map is undefined are points of
the form $l \cap e$, where $l$ is a line of $L$ and $e$ is
a line extending an edge of our kite.  One can
check easily, given Equation \ref{basicshape}, that all such
points have first coordinates of the form $m+An$, where $m,n \in \Z$.
On the other hand, no point of $T(\Z^2) \times \{\pm 1\}$ has this
form, and no iterate of such a point has this form.
\endproof

Let $\Upsilon$ denote the square of the outer billiards map.
Let ${\cal C\/}=T(\Z^2) \times \Z_{\rm odd\/}$,
where $\Z_{\rm odd\/}$ is the set of odd integers.
The vector $\Upsilon(x)-x$ is always twice the difference
between two of the vertices of our shape $S$.  Given
Equation \ref{main}, we see that
the first coordinate of $\Upsilon(x)-x$ always has the
form $2m+2nA$ where $m$ and $n$ are integers.  The
second coordinate is always an even integer.  Hence
$\Upsilon$ is entirely defined on $\cal C$ and preserves
this set.

\subsection{The Arithmetic Graph}

Let 
\begin{equation}
{\cal C\/}(\pm)=(T(\Z^2) \cap \R_+) \times \{\pm 1\},
\end{equation}
where
$T$ is as in Equation \ref{map}, and $A=\phi^{-3}$.
(Henceforth we take $A=\phi^{-3}$.)
The set ${\cal C\/}(\pm)$ is dense in the union of two rays starting at
$(0,\pm 1)$ and parallel to $(1,0)$.    
The point in Equation \ref{basicpoint} belongs to
${\cal C\/}(-)$.

\begin{center}
\psfig{file=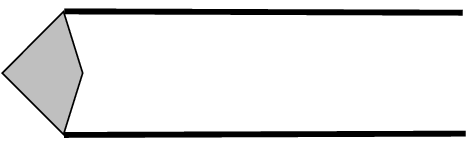}
\newline
{\bf Figure 2.1:\/} The objects $S$ and 
${\cal C\/}(\pm)$.
\end{center}

$\Upsilon$ does not preserve $\cal C(\pm)$, but we have \footnote{A ``runway'' like the one in
Figure 2.1, as well as
a map somewhat like $\Upsilon_R$, is also considered in [{\bf GS\/}].}
the first return map:

\begin{equation}
\Upsilon_R: {\cal C\/}(\pm) \to {\cal C\/}(\pm).
\end{equation}

Our idea is to encode the structure of 
$\Upsilon_R$ graphically.  We join
$x_1,x_2 \in \Z^2$ by a segment iff $T(x_j)>0$ for $j=1,2$, and
there are choices
$\epsilon_1,\epsilon_2 \in \{-1,1\}$ such that
\begin{equation}
\Upsilon_R(T(x_1),\epsilon_1)=(T(x_2),\epsilon_2).
\end{equation}
There is a subtlety in our definition, caused by the 
asymmetric roles played by $x_1$ and $x_2$ in the construction.
One might worry that this asymmetry could lead to an
inconsistency, whereby we are told to join $x_1,x_2$ by an
edge, but then told not to join $x_2,x_1$ by an edge!
However, we observe that reflection in the $X$-axis interchanges
${\cal C\/}(+)$ and ${\cal C\/}(-)$ and conjugates the outer
billiards map to its inverse.  Therefore
$$\Upsilon_R(T(x_1),\epsilon_1)=(T(x_2),\epsilon_2) \hskip 15 pt
\Longleftrightarrow \hskip 15 pt
\Upsilon_R(T(x_2),-\epsilon_1)=(T(x_1),-\epsilon_2).$$
Hence, our definition is consistent.

We let $\Gamma \subset H$ denote the graph resulting from our construction.
Here $H=T^{-1}(\R_+) \subset \R^2$ is an open halfplane containing
$\Gamma$. We call $\Gamma$
the {\it arithmetic graph\/}.

Figure 2.2 shows a small portion of $\Gamma$.
The dot is $(0,0)$.  The component $\Gamma_0$ containing
$(0,0)$ is
drawn in black.  The grid indicates $\Z^2$.  The
black line running along the bottom is $\partial H$.

\begin{center}
\psfig{file=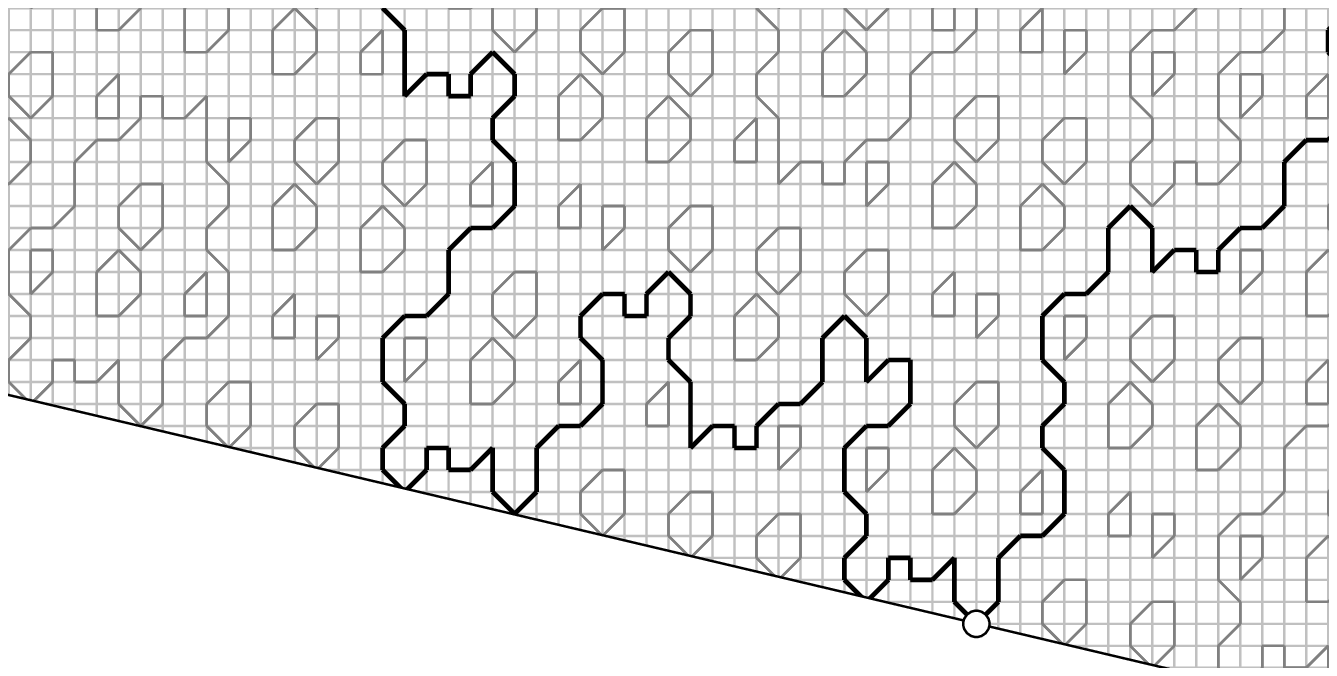}
\newline
{\bf Figure 2.2:\/} Part of the arithmetic graph
\end{center}

\begin{center}
\psfig{file=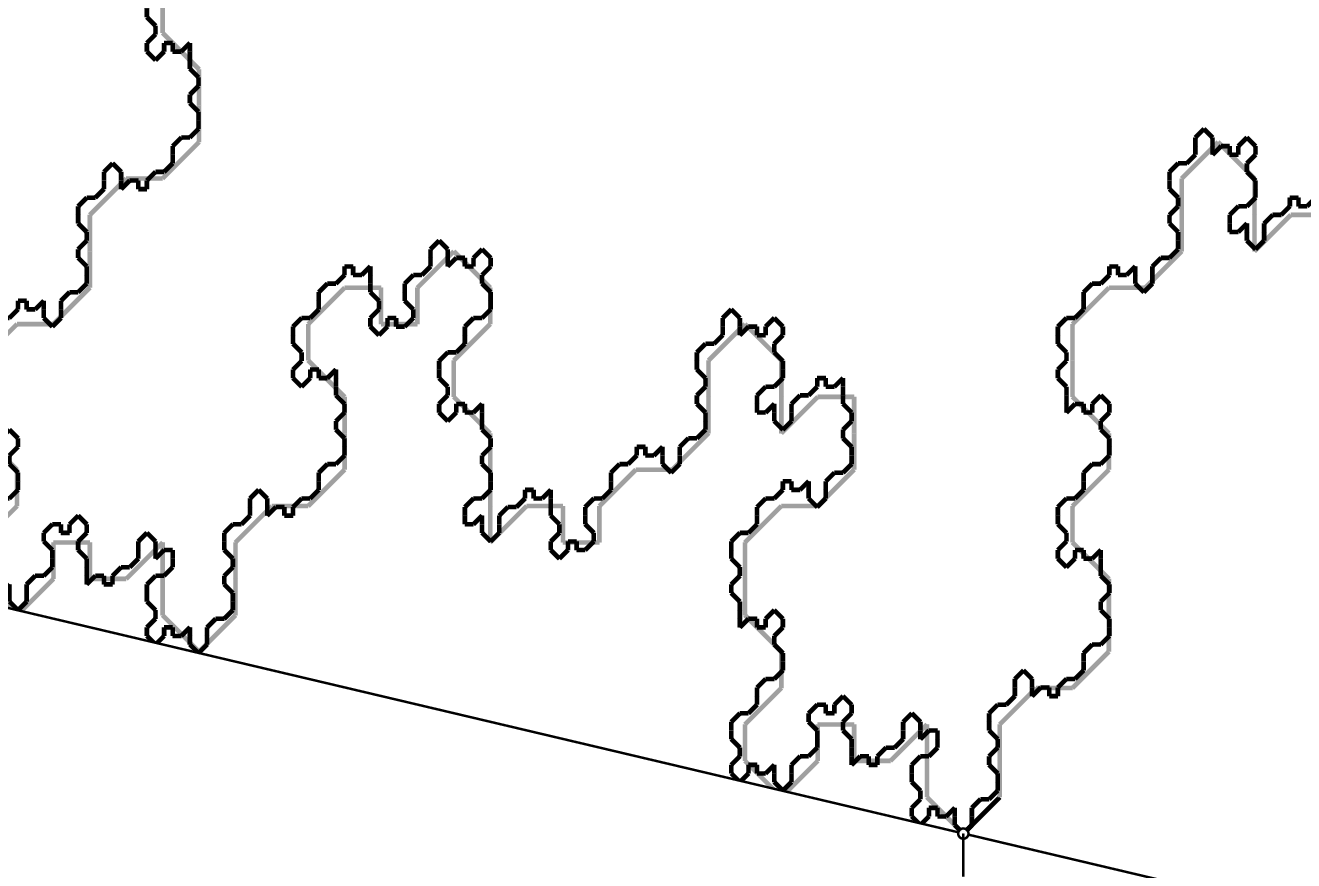}
\newline
{\bf Figure 2.3:\/} More of the component $\Gamma_0$.
\end{center}

Figure 2.3 shows more of $\Gamma_0$, again drawn in black.
The origin is denoted by a little vertical line segment
touching $\partial H$ at its top endpoint.
We have erased the other components and deleted the grid,
to get a clearer picture.  (The reader should use
Billiard King to see really great pictures of $\Gamma$.)
The curve $\Gamma_0$, and indeed all of $\Gamma$, behaves like the
self-similar inflationary tilings studied in [{\bf Ke\/}].  
One might call $\Gamma$ a ``large-scale fractal''.

Here is what we mean by
{\it inflationary\/}:
In Figure 2.3,
there is also a grey curve running alongside $\Gamma_0$.
This grey curve is the dilated image
$\phi^3\Gamma_0$.  In other words, we dilate $\Gamma_0$
by $\phi^3$ (about the origin), color it grey,
and superimpose it on the original picture.   From
what we can see of the picture, it looks like 
$\Gamma_0$ and $\phi^3\Gamma_0$ closely follow
each other.  One might say that
$\Gamma_0$ is {\it quasi-invariant\/} under a dilation.
If this is true$-$and we will prove it below$-$then
both ends of $\Gamma_0$ must exit every tubular neighborhood
of $\partial H$.  (Below we formulate this principle precisely.)
 Theorem \ref{main}
follows immediately from the fact that both ends of
$\Gamma_0$ exit every tubular neighborhood of
$\partial H$.

\subsection{The Arithmetic Graph Lemma}

Let $\T^2=(\R/\Z)^2$ be the square torus. Given
$p \in \R^2$ let $[p]$ denote the projection to $\T^2$.  We define
$\Psi: \Z^2 \to \T^2$ with the equation

\begin{equation}
\label{map2}
\Psi(x,y)=\bigg[\Big(\frac{T(x,y)}{2\phi},\frac{T(x,y)}{2}\Big)\bigg]=
\bigg[\Big(\phi^{-4}x+\phi^{-1}y,\phi^{-3}x\Big)+\frac{1}{2}\Big(\phi^{-3},\phi^{-2}\Big)\bigg].
\end{equation}
Here $T$ is the map from Equation \ref{map}.  The second equation is a direct calculation,
which we omit.  $\Psi(\Z^2)$ is dense in $\T^2$.

Given  $v \in \Z^2 \cap H$ we define
the {\it local type\/} of $v$ to be the translation equivalence
class of the union of edges of
$\Gamma$ emanating from $v$. It turns
out that there are $16$ types (including the type
where no edges emanate from $v$) but we prefer to
label these $16$ types by the integers
$1,...,26$ according to Figure 2.4.  Most of these
types can be seen in Figures 2.2 and 2.3.

\begin{center}
\psfig{file=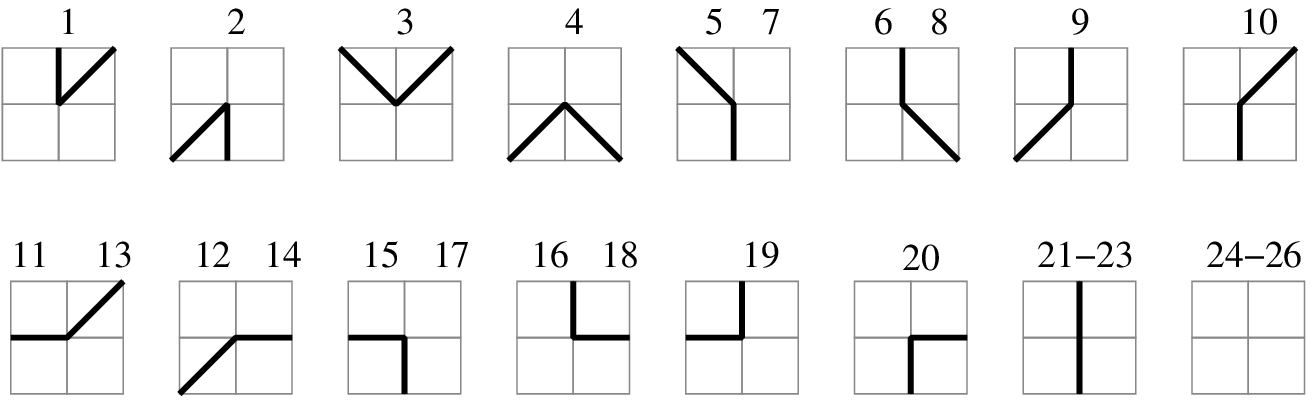}
\newline
{\bf Figure 2.4:\/} The local types
\end{center}

Figure 2.5 shows a partition of $\T^2$ into $26$ open
convex polygons
$\P_1,...,\P_{26}$.  

\begin{center}
\psfig{file=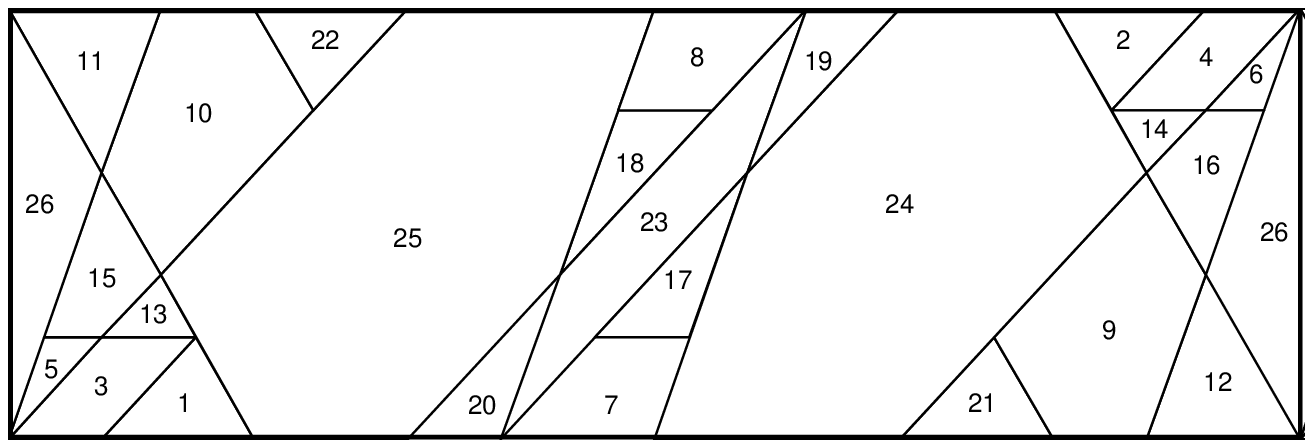}
\newline
{\bf Figure 2.5:\/} The Torus Partition
\end{center}

We have changed the aspect ratio
so as to get a nicer picture.  The true geometric
picture is contained in the unit square, with sides
identified as usual.   The origin $(0,0)$ is the
bottom left corner.   The only polygon that
``wraps'' around is $\P_{26}$.
In \S 6 we will prove

\begin{lemma}[Arithmetic Graph]
$\Psi$ maps each $x \in \Z^2 \cap H$ into some open polygon
of $\cal P$, and $x$ has local type $k$ if and only if
$\Psi(x) \in \P_k$.
\end{lemma}

The next lemma is not needed for our theorems, but it is
a nice fact to know and it illustrates the power of
the Arithmetic Graph Lemma.

\begin{lemma}
\label{embed}
The arithmetic graph is an embedded union of polygons
and polygonal arcs.
\end{lemma}

\startproof
(Sketch.)  First of all, we check the list of types that appear in the
Arithmetic Graph Lemma and we observe that every vertex in the graph has
valence 2.  Thus, the only kind of crossing that can occur is a transversal
crossing, where the cross point is not in $\Z^2$.  In particular,
there would have to be a point $(x,y)$ of type $a$ and a point
$(x+1,y)$ of type $b$, where $a \in \{4,6,8\}$ and
$b \in \{2,4,9,12,14\}.$   But $v=\Psi(x+1,y)-\Psi(x,y)$ is
the vector $(\phi^{-3},\phi^{-4})$.  This vector is
the one that points from the lower left corner of $\P_3$
to the upper right corner or $\P_3$ in Figure 2.5,
when it is anchored at the point $(0,0)$.   Looking
carefully at Figure 2.5 we see that no position of $v$ has
this property.
\endproof

\subsection{Inflation Structures}

As we mentioned above, the
basic idea in our proof that $\Gamma_0$ rises unboundedly
far away from $\partial H$ is to show that 
$\Gamma_0$ is quasi-invariant under
the dilation 
\begin{equation}
\Phi(x,y)=(\phi^3x,\phi^3y).
\end{equation}
In this section we formalize this idea.

Say that a {\it gene\/} is a (combinatorial) length $6$ connected
component of $\Gamma_0$.  
It turns out that there
are $75$ genes up to translation equivalence.  
Each gene $A$ has a {\it core\/}
$B$, consisting of the central path of length $2$.
The left hand side of Figure 2.6 shows a gene and its core.
We say that $p' \in \Z^2$ {\it shadows\/} $p \in \R^2$ if
$\|p-p'\|<4$.  (Here $4$ is a convenient cutoff.)
Let $B$ be a gene core, as above.  We say that a
finite polygonal path $A' \subset \Gamma$ {\it shadows\/}
$\Phi(B)$ if each endpoint of $A'$ shadows a corresponding
endpoint of $\Phi(B)$.  Figure 2.6 shows a
combinatorially accurate example.

\begin{center}
\psfig{file=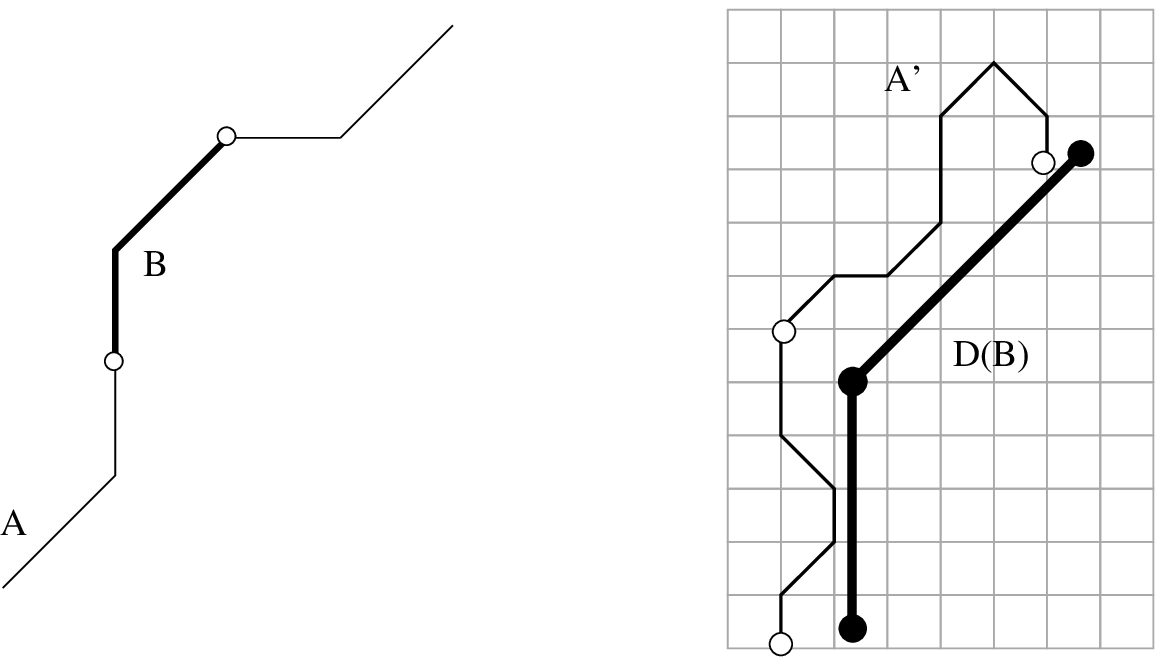}
\newline
{\bf Figure 2.6\/} genes, gene cores, and shadowing
\end{center}

Let $\cal G$ denote the set of all genes.   Let $\cal X$ denote the set of all
finite polygonal paths of $\Gamma$.
We say that an
{\it inflation structure\/} is an assignment $\chi$, to each
gene $A \in \cal G$ a path $A'=\chi(A)$ in $\cal X$ which shadows
$\Phi(B)$.  We emphasize that $\chi(A)$ only shadows the
dilated core $\Phi(B)$, but perhaps depends on all of
$A$ for its definition.  We also emphasize that the translation
type of $\chi(A)$ only depends on the gene type of $A$.
Billiard King computes part of an inflation structure
and allows the interested user to interact with it.

We say that a inflation structure is {\it coherent\/} if
\begin{itemize}
\item If $A$ is \underline{the} gene whose center vertex is $(0,0)$, then
$A \subset \chi(A)$. We call $A$ the {\it zeroth gene\/}.
\item If $A_1$ and $A_2$ are two consecutive genes, then
$\chi(A_1) \cup \chi(A_2)$ is a single polygonal path.
In other words, $\chi(A_1)$ and $\chi(A_2)$ splice
together seamlessly, just as $A_1$ and $A_2$ splice together.
\end{itemize}

\begin{lemma}
\label{inflrise}
If $\Gamma_0$ admits a coherent inflation structure then either
direction of
$\Gamma_0$ rises unboundedly far away from $\partial H$,
and hence Theorem \ref{main} is true.
\end{lemma}

\startproof
We see explicitly by direct computation that either
direction of $\Gamma_0$ contains
a polygonal arc which starts at $(0,0)$ and rises $100$ units
away from $\partial H$.  Let $\Gamma_1$ be one of these
two directions.
Say that we have already shown that $\Gamma_1$ rises up
$d \geq 100$ units away from $\partial H$.
Let $A_0,A_1,A_2,...$ be the consecutive genes of $\Gamma_1$,
rising up to this height,
chosen so that $A_0$ is the zeroth gene.
Then $A_0 \subset \chi(A_0)$ and the
consecutive paths $\chi(A_0)$, $\chi(A_1)$, $\chi(A_2)$,...
fit together.  These paths therefore trace out
$\Gamma_1$  and rise
up at least (say) $\phi^3d-100$, a quantity larger than $2d$.
Thus, we have now shown that $\Gamma_1$ rises up
$2d$ units from $\partial H$ in the same direction.
Iterating, we get our
result. \endproof

\begin{lemma}
\label{proper}
If $\Gamma_0$ admits a coherent inflation structure
then the distance from the $N$th point of
$O_+(p)$ to the origin is a proper function of $N$.
\end{lemma}

\startproof
Let $\Gamma_0^+$ be the component of
$\Gamma_0-(0,0)$ whose first (say) $100$ iterates are
contained in the positive quadrant $\R^2_+$.  (See Figure 2.3.)
An argument just like the one given for
Lemma \ref{inflrise} shows that $\Gamma_0^+ \subset \R^2_+$.
As we trace out $\Gamma_0^+$, a non-periodic curve
with integer vertices,
we can only return to any given compact neighborhood of
$(0,0)$ finitely many times.  Hence, the distance
from the $N$th vertex of $\Gamma_0^+$ to
$(0,0)$ is a proper function of $N$.  Since
$\Gamma_0^+ \subset \R^2_+$ and
$\partial H$ has negative slope, we now see that
the distance from the $N$th vertex of
$\Gamma_0^+$ to $\partial H$
is a proper function of $N$.  Hence the
Euclidean norm of the $N$th point of
$O_+(p) \cap {\cal C\/}(\pm)$ is a proper function of $N$.
But now observe that the norms of the points in
$O_+(p)$ between two consecutive 
points of $O_+(p) \cap {\cal C\/}(\pm)$ are, up to
a uniformly bounded factor, the same as the norms as
the two points of $O_+(p) \cap {\cal C\/}(\pm)$.
(See Figure 6.1. and the Pinwheel Lemma of \S 6.2.)
\endproof

\newpage

%% file: 3shadow.tex
\section{Manufacturing an Inflation Structure}

In this chapter we assume the Arithmetic Graph Lemma
and use it to delve more deeply into the
connection between the arithmetic graph $\Gamma$
and the torus partition $\cal P$.

\subsection{Pointed Strands and Dynamical Polygons}

Say that a
{\it pointed strand\/} of $\Gamma$ is a pair
$(X,x)$, where $X$ is a finite polygonal arc of $\Gamma$
having length at least $2$, 
and $x$ is a vertex of $X$.   Say that two pointed strands
$(X_1,x_1)$ and $(X_2,x_2)$ are {\it equivalent\/} if
there is a translation of $\Z^2$ carrying one
to the other, and
$\Psi(x_1)$, $\Psi(x_2)$ belong to the
same polygon $\P_j$ of $\cal P$.
Say that a {\it pointed strand type\/} of {\it length\/} $n$ is an
equivalence class of pointed strands, where the strands have $n$ segments.
Let $\Sigma_{\sigma}$ denote those $x \in \Z^2 \cap H$
such that $\sigma=[(X,x)]$.

\begin{lemma}[Dynamical Polygon]
There exists an open convex polygon 
$P_{\sigma} \subset \T^2$ such that
$a \in \Sigma_{\sigma}$ if and only if
$\Psi(a) \in P_{\sigma}$.
\end{lemma}

\startproof
We first consider the effect of keeping the strand but
changing the basepoint.  Suppose
$\sigma_1=[(X,x_1)]$ and
$\sigma_2=[(X,x_2)]$, where
$x_1$ and $x_2$ are adjacent vertices of $X$.   
We have
$x_2=x_1+(\epsilon_1,\epsilon_2)$ where
$\epsilon_j \in \{-1,0,1\}$.
Equation \ref{map2} gives us
$\Psi(x_2)=\Psi(x_1)+C$, where $C$ only depends on
$\epsilon_1$ and $\epsilon_2$.
If we have managed to prove this lemma for $\sigma_1$
then we define
$P_{\sigma_2}=P_{\sigma_1}+C$ and the lemma follows
for $\sigma_2$.  We call this the
{\it vertex slide argument\/}.

By the Arithmetic Graph Lemma, this lemma is true for the
local types.  The vertex slide argument now establishes this
lemma for all pointed strand types of length $2$.
Now suppose that
$\sigma=[(X,x)]$ has length at least $3$.
By the vertex slide argument, we can assume that
$x$ is an interior vertex of $X$.
Hence $X=X_1 \cup X_2$ where $x \in X_j$ and $X_j$ is a shorter
strand than $X$.   By induction there are open convex
polygons $P_{\sigma_1}$ and $P_{\sigma_2}$ such that
$a \in \Sigma_{\sigma_j}$ iff $\Psi(a) \in P_{\sigma_j}$.
Let $P_{\sigma}=P_{\sigma_1} \cap P_{\sigma_2}$.
Evidently $P_{\sigma}$ is an open convex polygon.
If $a \in \Sigma_{\sigma}$ then $a \in \Sigma_{\sigma_1} \cap
\Sigma_{\sigma_2}$ and by induction $\Psi(a) \in P_{\sigma_1} \cap P_{\sigma_2}=
P_{\sigma}$.  Conversely suppose that $\Psi(a) \in P_{\sigma}$. Then
$\Psi(a) \in P_{\sigma_j}$ and hence
there is a strand $A_j$ such that
$\sigma_j=[(A_j,a)]$.  Letting $A=A_1 \cup A_2$ we
have $\sigma=[(A,a)]$.  Hence
$a \in \Sigma_{\sigma}$.
This completes the induction step.
\endproof

We call $P_{\sigma}$ a {\it dynamical polygon\/} because,
as we explain in \S 4, it can be produced by a dynamical
process.

\subsection{Inflation Maps}

In this section and the next we
build some machinery for describing the interaction
between the dilation map
$$\Phi(x,y)=(\phi^3 x,\phi^3 y)$$
and the map
$\Psi: \Z^2 \to \T^2$ given in Equation \ref{map2}.  That is,
$$\Psi(x,y)=
\bigg[\Big(\phi^{-4}x+\phi^{-1}y,\phi^{-3}x\Big)+\frac{1}{2}\Big(\phi^{-3},\phi^{-2}\Big)\bigg].
$$
Our goal is to define, in a canonical way, a kind of
integral approximation to $\Phi$.   Such approximations are
vital to the creation of an inflation structure.

We say that a {\it small polygon\/} is an open convex
polygon of $\T^2$ that is contained inside a square of
side length $1/2$.  The polygons of interest to us will
all be small.
Let $P$ be a small polygon.
We say that a map $\gamma: P \to \T^2$
is a {\it special similarity\/} if, relative to local
Euclidean coordinates on $P$ and $\gamma(P)$ the
map $\gamma$ is just multiplication by
$-\phi^{-3}$.  Put another way, $\gamma$ is a holomorphic
map on $P$ whose complex derivative is constantly equal
to $-\phi^{-3}$.  

Define
\begin{equation}
\Sigma(P)=\Psi^{-1}(P) \cap \Z^2.
\end{equation}
If $P=P_{\sigma}$ then $\Sigma(P) \cap H=\Sigma_{\sigma}$, the
set defined in the last section.
We say that a map
$\beta: \Sigma(P) \to \Z^2$ is an {\it inflation map\/} if
\begin{equation}
\Psi(\beta(a))=\gamma(\Psi(a)); \hskip 30 pt
\forall a \in \Sigma(P),
\end{equation}
for some special similarity $\gamma: P \to \T^2$.
In the next section we prove the existence of such maps.
Here we investigate some of their abstract properties.

\begin{lemma}
\label{tweak1}
Suppose that $\beta$ is an inflation map defined relative to $P$.
Let $\widehat \beta$ be the new map defined by the equation
$\widehat \beta(a)=\beta(a)+b_0$ for some $b_0 \in \Z^2$.
Then $\widehat \beta$ is also an inflation map defined relative
to $P$.
\end{lemma}

\startproof
We compute that
$$\Psi(\widehat \beta(a))=\Psi(\beta(a))+C=\gamma(\Psi(a))+C=
\widehat \gamma(\Psi(a)).$$
Here $C$ is some constant and $\widehat \gamma$ is some
other special similarity which differs from $\gamma$ only
by a translation.  Our equation makes sense because
$\T^2=R^2/\Lambda$ is canonically an abelian group.
\endproof

\begin{lemma}
\label{tweak2}
Suppose that $\beta$ is an inflation map defined relative to $P$.
Given $a_0 \in \Z^2$ let $\widehat P=P+\Psi(a_0)$.  Define
$\widehat \beta: \Sigma(\widehat P) \to \Z$ as
$\widehat \beta(a)=\beta(a-a_0)$.  Then $\widehat \beta$ is
an inflation map defined relative to $\widehat P$.
\end{lemma}

\startproof
If $a \in \Sigma(\widehat P)$ then $\Psi(a) \in \widehat P$.  But then
$\Psi(a-a_0)=\Psi(a)-\Psi(a_0) \in P$.  Thus
$\beta(a-a_0)$ is always defined.
We compute
$$\Psi(\widehat \beta(a))=\Psi(\beta(a-a_0))=
\Psi(\beta(a))-C = \gamma(\Psi(a))-C=\widehat \gamma(\Psi(a)).$$
Here $C$ is some constant.
The rest of the proof is as in the previous lemma.
\endproof

\begin{lemma}
\label{agree}
Let $\beta_1$ and $\beta_2$ be two inflation maps defined relative
to the same polygon $P$.  If $\beta_1(a)=\beta_2(a)$ for some
$a \in \Sigma(P)$ then $\beta_1=\beta_2$.
\end{lemma}

\startproof
Let $\gamma_j$ be such that
$\gamma_j(\Psi(a))=\Psi(\beta_j(a))$.  Note that
$\gamma_1$ and $\gamma_2$ must differ by a translation.
Thus, these two maps are equal if they agree at any point.
But 
$\gamma_1(\Psi(a))=\Psi(\beta_1(a))=\Psi(\beta_2(a))=
\gamma_2(\Psi(a)).$
Hence $\gamma_1=\gamma_2$.   Now we know that
$\Psi(\beta_1(a))=\Psi(\beta_2(a))$ for all
$a \in \Sigma(P)$.  To finish our proof,
we note that $\Psi$ is pretty obviously injective on
$\Z^2$.
\endproof

Now we lay the groundwork for explaining how
$\beta$ is an integral approximation to $\Phi$.
Given an inflation map $\beta$ defined relative to $P$ and
$a \in \Sigma(P)$, let
\begin{equation}
v_{\beta}(a)=\Phi(a)-\beta(a).
\end{equation}
Here $v_{\beta}$ measures the discrepancy between $\beta$ and
the dilation $\Phi$.  Say that $\beta$ is
$K$-{\it pseudo-Lipschitz\/} if
\begin{equation}
\|v_{\beta}(a_1)-v_{\beta}(a_2)\| \leq K \|\Psi(a_1)-\Psi(a_2)\|.
\end{equation}

\begin{lemma}
\label{lip}
Suppose that the inflation map $\beta$, defined relative to $P$,
is $K$-pseudo-Lipschitz.  
Then any inflation map $\widehat \beta$ defined relative to $P$
is $K$-pseudo-Lipschitz.
\end{lemma}

\startproof
In light of Lemma \ref{agree}, we must have, for all
$a \in \Sigma(P)$, the identity.
$\widehat \beta(a)=\beta(a)+a_0$ for some $a_0 \in \Z^2$.
The result is obvious from this identity.
\endproof

\subsection{The Shadow Lemma}

Recall that convex polygon $P \subset \T^2$ is
{\it small\/} if $P$ is contained inside an open square of radius
$1/2$ in $\T^2$. 

\begin{lemma}[Shadow]
Let $P$ be a small polygon. There exists
an $4$-pseudo-Lipschitz inflation map $\beta$
defined relative to $P$.
\end{lemma}

It suffices to prove the Shadow Lemma for the map
$$\Psi_0(x,y)=
\bigg[\Big(\phi^{-4}x+\phi^{-1}y,\phi^{-3}x\Big)\bigg]$$
which differs from $\Psi$ by a translation of $\T^2$.
That is, we will produce a special similarity $\gamma$
such that $\Psi_0(\beta(a))=\gamma(\Psi_0(a))$ for all $a \in \Sigma(P)$.

Given that $P$ is small, there are intervals
$I_1,I_2 \subset \R$, each having length less than
$1/2$ with the following property:  For
$a=(x,y) \in A$, we have unique integers $m$ and $n$ such that
\begin{equation}
\label{lattice0}
\phi^{-4}x+\phi^{-1}y= m + \epsilon_1; \hskip 30 pt
\phi^{-3} x = n + \epsilon_2; \hskip 30 pt
\epsilon_j \in I_j.
\end{equation}
We write $\mu \sim \mu'$ if $\mu-\mu' \in \Z$.
Since $\phi^3 \sim \phi^{-3}$ and $x \in \Z$ we have
\begin{equation}
\label{lattice00}
\phi^3 x \sim \epsilon_2.
\end{equation}
Since $\phi^3 \sim 2\phi^{-1}$ and $-2\phi^{-4} \sim 3 \phi^{-3}$
and $x,y \in \Z$ we have
\begin{equation}
\phi^3 y \sim 2\phi^{-1} y =2(\phi^{-4}x+\phi^{-1}y)-2\phi^{-4}x \sim
2(m+\epsilon_1)+3 \phi^{-3}x \sim 2 \epsilon_1+3 \epsilon_2.
\end{equation}
Accordingly we define
\begin{equation}
\beta(x,y)=(\widehat x,\widehat y); \hskip 30 pt
\widehat x= \phi^3x-\epsilon_2; \hskip 30 pt
\widehat y= \phi^3y-2\epsilon_1-3\epsilon_2.
\end{equation}
By construction $(\widehat x,\widehat y) \in \Z^2$.  

\begin{lemma}
\label{needsmall}
For any $a_1,a_2 \in \Sigma(P)$ we have
$$\|v(a_1)-v(a_2)\|<4\|\Psi_0(a_1)-\Psi_0(a_2)\|.$$
\end{lemma}

\startproof
Let $(\epsilon_{1j},\epsilon_{2j})$ be constants associated to the
point $a_j=(x_j,y_j)$ as in Equation \ref{lattice0}. 
Let $\delta_j=\epsilon_{1j}-\epsilon_{2j}$ for $j=1,2$.
First, we have
$$\Psi_0(a_1)-\Psi_0(a_2)=(\delta_1,\delta_2).$$  Second, we have
$$v(a_j)=\Phi(a_j)-\beta(a_j)=(\phi^3x-\widehat x,\phi^3 y-\widehat y)=
(\epsilon_{2j},2\epsilon_{2j}+3\epsilon_{2j}).$$
Hence
$$v(a_1)-v(a_2)=(\delta_2,2\delta_1+3\delta_2).$$
The matrix associated to the linear map
$(\delta_1,\delta_2) \to (\delta_2,2\delta_1+3\delta_2)$
has $L_2$-norm equal to $\sqrt{14}<4$ and hence the map itself
is $4$-Lipschitz.  
\endproof

We compute
\begin{equation}
\label{lattice4}
\phi^{-3} \widehat x=\phi^{-3}(\phi^3x-\epsilon_2)=x-\phi^{-3}\epsilon_2 \sim -\phi^{-3}\epsilon_2.
\end{equation}

\begin{lemma}
\label{lattice5}
$\phi^{-4}\widehat x + \phi^{-1} \widehat y \sim -\phi^{-3} \epsilon_1$.
\end{lemma}

\startproof
We have
\begin{equation}
\label{lattice6}
\phi^{-4} \widehat x+\phi^{-1}\widehat y=
\phi^{-1} x+\phi^2 y   -\phi^{-4}\epsilon_2 -2 \phi^{-1}\epsilon_1 -3 \phi^{-1}\epsilon_2.
\end{equation}
But 
\begin{equation}
\label{lattice7}
\phi^{-1}x + \phi^2 y \sim \phi^{-1}x +\phi^{-1} y=
(\phi^{-1}-\phi^{-4})x+(\phi^{-4}x +\phi^{-1}y)=$$
$$(2\phi^{-3})x +(\phi^{-4}x +\phi^{-1}y) \sim
\epsilon_1+2\epsilon_2.
\end{equation}
Our last equation comes from Equations \ref{lattice0} and \ref{lattice00}.
Plugging Equation \ref{lattice7} into Equation
\ref{lattice6} and grouping terms, we have
$$\phi^{-4} \widehat x + \phi^{-1} \widehat y \sim
(-2 \phi^{-1}+1)\epsilon_1 + (-\phi^{-4}-3\phi^{_1}+2)\epsilon_2=-\phi^{-3}\epsilon_1.$$
This completes the proof
\endproof

Equation \ref{lattice4} and Lemma \ref{lattice5} together show that
$\Psi_0(\widehat x,\widehat y)=-\phi^{-3}\Psi_0(x,y)$ in local isometric coordinates.
This is what we wanted to prove.  This completes the proof of the Shadow Lemma.

\subsection{Inflation Maps and Dynamical Polygons}

Now we put together the material from the previous two
sections.  The result we prove here, an immediate
consequence of what we have already done, provides the
key step in constructing a coherent inflation
structure from a finite amount of computational
information.

Let $\sigma_1$ and $\sigma_2$ be two pointed strand types.
Let $P_{\sigma_1}$ and $P_{\sigma_2}$ be the associated
dynamical polygons.  We write
\begin{equation}
\sigma_1 \longrightarrow_{\beta} \sigma_2
\end{equation}
if there is a inflation map $\beta$, defined relative
to $P_{\sigma_1}$ and having the property that

\begin{equation}
\gamma(P_{\sigma_1}) \subset P_{\sigma_2}.
\end{equation}
Here $\gamma$ is the special similarity associated to $\beta$.

\begin{lemma}
\label{shadow3}
Suppose $\sigma_1 \longrightarrow_{\beta} \sigma_2$.
Suppose also that $a_1,a_2 \in \Z^2 \cap H$ are such that
$a_2=\beta(a_1)$ and
$\sigma_1=[(A_1,a_1)]$ for some polygonal arc $A_1$.
Then there exists a polygonal arc $A_2 \subset \Gamma$ such that
$\sigma_2=[(A_2,a_2)]$.
\end{lemma}

\startproof
We have $a_1 \in \Sigma_{\sigma_1} \subset \Sigma(P_{\sigma_1})$ 
by the Dynamical Polygon Lemma.
By definition
$\Psi(a_2)=\gamma(\Psi(a_1)) \subset P_{\sigma_2}$.  By the Dynamical Polygon Lemma
we have $\sigma_2 \in \Sigma_{\sigma_2}$.  This means that
there is an arc $A_2$ of $\Gamma$ such that
$\sigma_2=[(A_2,a_2)]$.
\endproof

\subsection{The Inflation Generator}

Now we turn our attention to the genes we discussed at the end of \S 2.
Recall that a {\it gene\/} is a polygonal arc of length $6$ contained
in $\Gamma_0$.  Here $\Gamma_0$ is the component of $\Gamma$ that
contains $(0,0)$.
Each gene $A$ gives rise to a pointed strand
$(A,a)$ where $a$ is the central vertex of $A$.
Each gene $A_j$ therefore gives rise to the
dynamical polygon $P_j$ associated to the
pointed strand type $[(A_j,a_j)]$.   It turns out that
there are $75$ such dynamical polygons, corresponding
to $75$ distinct combinatorial types of gene. The corresponding
polygons are the dark-shaded polygons in Figure 3.1.
(We will explain the significance of the $24$ light-shaded
polygons in the next chapter.)
One can see easily
from the picture, or from inspecting the list in
\S \ref{polylist1}, that any two points in the same
dynamical polygon $P_j$ are within $1/4$ of each other.
In particular, these polygons are all small.

\begin{center}
\psfig{file=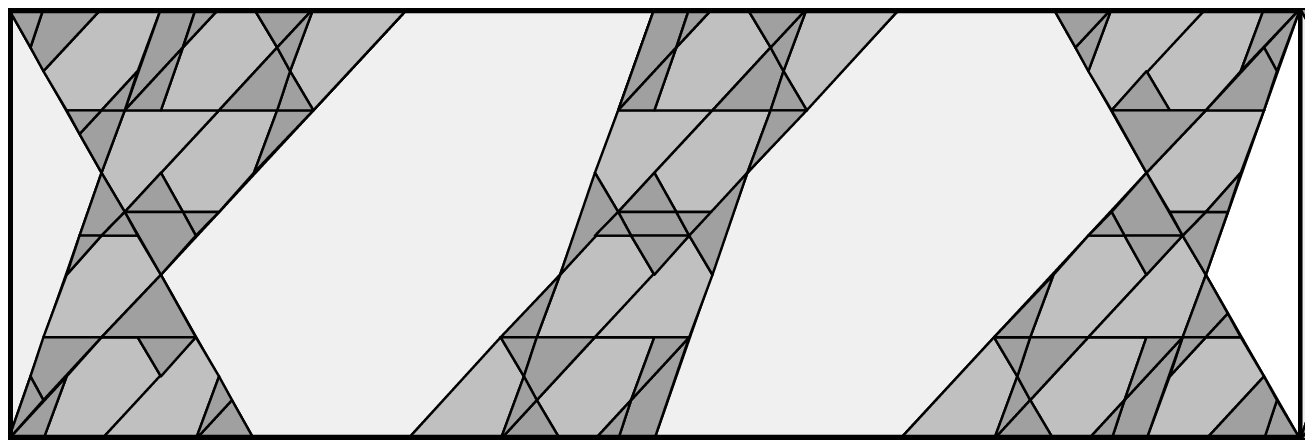}
\newline
{\bf Figure 3.1:\/} dynamical polygons corresponding to genes
\end{center}

We say that the gene $A$ (with central vertex $a$ and
gene core $B$) 
is {\it nicely shadowed\/} by
a polygonal arc $A'$ if there is a vertex $a'$ of $A'$
such that the following is true:
\begin{itemize}
\item There exists an inflation map $\beta$ with the properties that
$\beta(a)=a'$ and
$[(A,a)] \longrightarrow_{\beta}[(A',a')]$.
\item Each endpoint of $\Phi(B)$ is within $3$ units of an endpoint of $A'$.
\end{itemize}

We say that an {\it inflation generator\/} is a list $G$ of
$75$ pairs of the form $\{(A_i,A_i')\}$ such that
each gene of $\Gamma_0$ is translation equivalent to
one gene $A_i$ on the list, and
$A_i' \subset \Gamma$ nicely shadows $A_i$.

Let $A$ be some gene.  There is some gene $A_j$ on our list
such that $A$ and $A_j$ have the same type. 
Let $a$ be the central vertex of $A$ and let
$a_j$ be the central vertex of $A_j$.
Let $\beta_j$, $A'_j$, and $a'_j$ be the objects
associated to $A_j$, as above.
From Lemma \ref{shadow3} there
is a pointed strand $(A',a')$ such that
$[(A'_j,a'_j)]=[(A',a')]$ and
$a'=\beta_j(a)$.  We let
$\chi(A)=A'$.

\begin{lemma}
\label{inflate1}
$\chi$ is an inflation structure.
\end{lemma}

\startproof
For an arbitrary gene $A$, with gene core $B$, we need to show
that each endpoint of $\chi(A)$ is within $4$ units of the
corresponding endpoint of $\Phi(B)$.  We write $A'=\chi(A)$.
Let $b$ be one of the
endpoints of $B$.  For ease of labelling let's assume that
$A_1$ is the gene on our list of $75$ that has the same
type as $A$.  Let $B_1$ be the core of $A_1$ and let
$b_1$ be the vertex of $B_1$ corresponding to $B$.
Let $\beta$ be the inflation map associated to $A_1$.
Let $b_1'$ be the endpoint of $A_1'$ that is within
$3$ units of $B(b_1)$.  Let $P_1$ denote the
dynamical polygon associated to $[(A_1,a_1)]$ and
let $\widehat P_1$ be the dynamical polygon
associated to $[(A_1,b_1)]$.  Note that
$P_1$ and $\widehat P_1$ are translates of each other.

Consider a new map $\widehat \beta$ defined by the
rule
$$
\widehat \beta(x)=\beta(x+a_1-b_1)+b_1'-a_1'
$$
By Lemmas \ref{tweak1} and \ref{tweak2},
together with the fact that $\widehat P_1$ is the
relevant translate of $P_1$, we see that
$\widehat \beta$ is a inflation map defined
on $\Sigma(\widehat P_1)$.  By the Shadow Lemma
$\beta$ is $4$-pseudo-Lipschitz, and hence
so is $\widehat \beta$.
By translation equivalence, we have
$a-b=a_1-b_1$ and
$a'-b'=a_1'-b_1'.$
Therefore
\begin{equation}
\label{lucky}
\widehat \beta(b)=\beta(b+a_1-b_1)+b_1'-a_1'=
\beta(b+a-b)+b'-a'=b'.
\end{equation}
Note that
$$\|v_{\widehat \beta}(a_1)\|=\|b_1'-\Phi(b_1)\| <3$$
by the second property of our inflation generator.
Therefore
$$\|b'-\Phi(b)\|=^*\|\widehat \beta(b)-\Phi(p)\|=
\|v_{\widehat \beta}(b')\| \leq
\|v_{\widehat \beta}(b)\|+
\|v_{\widehat \beta}(b')-v_{\widehat \beta}(b_1')\|<$$
$$
3 + 4\ {\rm diam\/}(P_1)<3+4(1/4)=4.$$
The starred equality comes from Equation \ref{lucky}.
The $1/4$ in this last calculation comes from the fact
that any two points in the same dynamical polygon are
within $1/4$ of each other.
\endproof

\subsection{Coherence}
\label{coh}
Here we explain how to check that the inflation
structure $\chi$ is coherent.
Our inflation generator gives rise to
a list $\beta_1,...,\beta_{75}$ of inflation maps,
where $\beta_j$ is defined on 
$\Sigma(P_j)$.  Here $P_j$ is the dynamical
polygon associated to $[(A_j,a_j)]$.

Say that an {\it extended gene\/} is a polygonal arc of
$\Gamma$ having length $7$.  An extended gene is just the
union of two consecutive genes.   Say that an
{\it extended gene type\/} is an equivalence class,
up to translation, of extended genes.  It turns out that
there are $89$
extended gene types.  Given an extended gene $X$ we define
$\chi(X)=\chi(X_1) \cup \chi(X_2)$, where $X_1$ and $X_2$
are the two genes comprising $X$.

\begin{lemma}[Coherence]
Suppose there is a list $X_1,...,X_{89}$ of
extended genes, representing all the types,
such that $\chi(X_j)$ is a polygonal arc
for all $j=1,...,89$.  Then $\chi$ is coherent.
\end{lemma}

\startproof
Let $\widehat X$ be an arbitrary extended gene and let
$X$ be the translation equivalent extended gene on our list
of $89$.  It suffices to show that $\chi(X)$ and
$\chi(\widehat X)$ are translates of each other.
For each object $Y$ we associate to
$X$, we let $\widehat Y$ be the corresponding object for $\widehat X$.

Let $X_1$ and $X_2$ be the two genes comprising $X$. 
Let $x_1$ and $x_2$ be the two center points of these genes. 
Let $\beta_1$ and $\beta_2$ be the inflation maps
associated to the gene types of $X_1$ and $X_2$.
To prove that $\widehat X$ is translate equivalent to $X$ it
suffices to prove
\begin{equation}
\label{affgoal2}
\beta_1(\widehat x_1)-\beta_2(\widehat x_2)=
\beta_1(x_1)-\beta_2(x_2).
\end{equation}
This is what we will do.

Since $\Psi$ is affine and 
injective on $\Z^2$ it suffices to show
that
\begin{equation}
\label{affgoal}
\Psi(\beta_1(\widehat x_1))-\Psi(\beta_2(\widehat x_2))=
   \Psi(\beta_1(x_1))-\Psi(\beta_2(x_2)).
\end{equation}
Let $p_j=\Psi(x_j)$ and $\widehat p_j=\Psi(\widehat x_j)$.
Since $X$ and $\widehat X$ are translation equivalent, we have
$\widehat x_1-x_1=\widehat x_2-x_2$.  Since
$\Psi$ is affine, we have
\begin{equation}
\widehat p_1-p_1=\widehat p_2-p_2.
\end{equation}

Let $\gamma_j: P_j \to \T^2$ be the special
similarity associated to $\beta_j$. 
We work in local Euclidean coordinates, so that
$\gamma_j(x,y)=-\phi^{-3}(x,y)+C_j$ for
some constant $C_j$.  Then
$$\gamma_1(\widehat p_1)-\gamma_1(p_1) =
-\phi^{-3}(\widehat p_1-p_1)=
-\phi^{-3}(\widehat p_2-p_2)=
\gamma_2(\widehat p_2)-\gamma_2(p_2).$$
Hence
\begin{equation}
\label{affine}
\gamma_1(\widehat p_1)-\gamma_2(\widehat p_2)=
\gamma_1(p_1)-\gamma_2(p_2).
\end{equation}
But we also know that
\begin{equation}
\label{affine2}
\gamma_j(p_j)=\Psi(\beta_j(x_j)); \hskip 30 pt
\gamma_j(\widehat p_j)=\Psi(\beta_j(\widehat x_j)).
\end{equation}
Equations
\ref{affine} and \ref{affine2} combine to establish
Equation \ref{affgoal}.
\endproof

\newpage

%% file: 4verify.tex
\section{Computer Aided Verification}

\subsection{Overview}

Recall that a gene is a polygonal arc of length $6$ contained in $\Gamma_0$,
the component of the arithmetic graph containing $(0,0)$. 
We are trying to verify the existence of an inflation generator.
This inflation generator consists of a length $75$ list of the form
$\{(A_j,A'_j)\}$, where $A_j$ is a gene and $A'_j$ is a path that nicely
shadows $A_j$ in the sense of \S 3.5.
We associate some auxilliary objects
to our list, namely:
\begin{itemize}
\item Let $P_j$ be the dynamical polygon associated to $A_j$.
\item Let $\beta_j$ be the inflation map associated
to $(A_j,A'_j)$.  
\item Let $\gamma_j: P_j \to \T^2$ be the special similarity
associated to $\beta_j$.  Here
$\gamma_j \circ \Psi=\Psi \circ \beta_j$, whenever all
maps are defined.
\end{itemize}

Given these basic objects, here are the things we need to check:
\begin{enumerate}
\item Each endpoint of $A_j'$ is within $3$ units of the
corresponding endpoint of $\Phi(B_j)$, the dilation of the core $B_j$ of $A_j$.
We have a list of all the vertices involved and we
just check this directly.  See \S 8.
\item For each gene $A_j$, with center vertex $a_j$, we have
correctly computed
the dynamical polygons $P_j$ associated to
$[(A_j,a_j)]$.
\item
We have $\gamma_{j}(P_{j}) \subset P_{j}'$, where
$P_{j}'$ is the dynamical polygon associated to 
the pointed strand type $[(A_j',a_{j}')]$.  Our method
will not require us to compute $P'_{j}$ explicitly.
\item There is a complete list $X_1,...,X_{89}$ of representatives of extended
genes such that $\chi(X_j)$ is a polygonal arc for $j=1,...,89$.
Here {\it complete list\/} means that every extended gene type
is represented.  
\item Our list $A_1,...,A_{75}$ of gene types is exhaustive and
our list $X_1,...,X_{89}$ of extended gene types is exhaustive.  
There are no other gene types or extended gene types.
\end{enumerate}

In this chapter we will explain the main theoretical points
of our verifications.  In \S 7 we will write enough pseudo-code so that
the interested reader can see explicitly how we do all the
calculations.

\subsection{Computing the Polygons}

Let $A=A_j$ denote one of the genes on our list.
Let $a_{-2},...,a_2$ denote the $5$ interior
vertices of $A$.   We are interested
in computing the dynamical polygon $P$ associated
to the pointed strand type $[(A,a_0)]$.  Here we
explain how $P$ is computed.  We don't actually need to know
how to compute $P$ for our proof$-$we just need to verify the
answer is correct$-$but the method of computation suggests
how we verify that the answer is correct.

We define a {\it dynamical translation\/} of the torus
$\T^2$ to be a map of the form
\begin{equation}
\label{dt}
[(x,y)] \to [(x,y)+(\epsilon_1 \phi^{-4}+\epsilon_2 \phi^{-1},\epsilon \phi^{-3})]; \hskip 30 pt
\epsilon_1,\epsilon_2 \in \{-1,0,1\}.
\end{equation}
This operation makes sense because $\T^2=\R^2/\Z^2$  is canonically
an abelian group.  We remind the reader that
$[(x,y)]$ is the equivalence class of the point
$(x,y)$ in $\T^2$.
 We say that $(\epsilon_1,\epsilon_2)$ is the
{\it type\/} of the dynamical translation. 

Recall that $\cal P$ is partitioned into $26$ open polygons (and
their boundaries)
as in Figure 2.5.
We associate to $A$ a list $\Q_{-2},...,\Q_{2}$ of $5$ open polygons
of the partition $\cal P$. 
Here $\Q_j$ is the polygon having
the same local type as the vertex $a_j$.   Depending on $A$
there are $4$ dynamical translations
$T_{-2},T_{-1},T_1,T_2$ with the property that
$x \in P$ if and only if 
\begin{itemize}
\item $x \in \Q_0$.
\item $T_1(x) \in \Q_1$
\item $T_2(T_1(x)) \in \Q_2$.
\item $T_{-1}(x) \in \Q_{-1}$.
\item $T_{-2}(T_{-1}(x)) \in \Q_{-2}$.
\end{itemize}
We call these four conditions {\it Property X\/}.
Given this information we have
\begin{equation}
P=\Q_0 \cap T_1^{-1}(\Q_1 \cap T_2^{-1}(\Q_2)) \cap
T_{-1}^{-1}(\Q_{-1} \cap T_{-2}^{-1}(\Q_2)).
\end{equation}
Billiard King simply computes this intersection.

If we replace each set $\Q_j$ by its closure
$\overline{\Q_j}$, then we arrive at a criterion
for when $x \in \overline P$.  We call this
criterion {\it Property\/} $\overline X$.

Figure 3.1 (which we repeat as Figure 4.1 below, for convenience)
shows the $75$ dynamical polygons in dark grey.
We list the actual coordinates in \S \ref{polylist2}.
The $24$ light polygons comprise the complement.  Here
is their interpretation:  A vertex $p$ in the arithmetic
graph lies on a closed polygon of length $5$ or $7$
if and only if $\Psi(p)$ lies in the interior of
one of the light polygons.  We do not use these
light polygons in this paper, but we note that they ``round out''
Figure 4.1 by explaining what is left over after we plot
our $75$ dynamical polygons.  

\begin{center}
\psfig{file=Pix/partition2.ps}
\newline
{\bf Figure 4.1:\/} The partition
\end{center}  

For the first stage of our verification, we use the
Billiard King algorithm to compute $P$ to high precision.
Next, we use the fact that the vertices of
$P$ are elements of $\frac{1}{2}\Z[\phi]$
to guess the exact expression for
the vertices.  Let $P_G$ denote this guess.  Our
goal is simply to verify that $P_G=P$. 

First we check that each vertex of $P_G$ satisfies
Property $\overline X$. 
This tells us that $\overline{P_G} \subset \overline P$.
If $\overline{P_G}$ is a proper subset of $\overline P$
then there is some edge $e$ of $P_G$ such that every
point $x$ on the interior of $e$ satisfies Property $X$.
We rule this out by explicitly choosing one point
per each edge of $P_G$ and showing that it fails
to have Property X.  In \S 7 we explain the
calculation in detail.

When we run the computation for each of the $75$ cases it works.
Of course, the fact that Figure 4.1 is a partition of the torus
gives extremely strong visual evidence that we have guessed
correctly in the first place.

\subsection{Checking the Shadowing Property}

Our method here is quite similar to our method for verifying
that $P_G=P$.   
Let $A=A_j$ be one of our genes.  Let $\beta=\beta_j$, etc. be the objects
associated to $A$.  Our goal is to check that
$\gamma(P) \subset P'$. 

Let $\Q'_{-1},\Q'_0,\Q'_1,...$ be the sequence of polygons
in the partition $\cal P$ associated to the vertices of
$A'$.  We set up the indices so that the
$0$th vertex of $A'$ is the distinguished point $a'$
that shadows $\Phi(a)$. (Here $a$ is the center of the gene $A$.)

Just as we did for $A$, we
generate the sequence of
dynamical translations for the strand $A'$.  This is a sequence whose
length varies with the choice of gene.
The length typically varies from $10$ to $20$.
The condition that $x \in P'$ amounts to checking that $x$ satisfies
what we call {\it Property\/} $X'$:

\begin{itemize}
\item $x \in \Q_0'$.
\item $T_1(x) \in \Q'_1$ and $T_{-1}(x) \in \Q'_{-1}$
\item $T_2(T_1(x)) \in \Q'_2$ and $T_{-2}(T_{-1}(x)) \in \Q'_{-2}$
\item $T_3(T_2(T_1(x))) \in \Q'_3$ and $T_{-3}(T_{-2}(T_{-1}(x))) \in \Q'_{-3}$
\item etc.
\end{itemize}

We also have the corresponding Property $\overline {X'}$,
the closed version.
We simply check that each vertex of $\gamma(P)$
satisfies Property $\overline{X'}$.  This means
that the closure of $\gamma(P)$ is contained in the
closure of $P'$.  Hence $\gamma(P) \in P'$.
Billiard King computes all these quantities
and displays them visually so that the user can
see in each case that
$\gamma(P) \subset P'$. 

\subsection{Checking the Coherence}

We will use the forwards direction of $\Gamma_0$ to check
the coherence. 
In order to make the construction to follow we first
need to know something about a certain finite
portion of $\Gamma_0$.  We verify by direct
inspection that  $\Gamma_0$ contains a polygonal
arc of combinatorial length $2^{14}$, connecting $(0,0)$ to
a point in the positive quadrant. 

In this section we will use the notation and
terminology from \S \ref{coh}
and the Coherence Lemma.  In particular
$X_1,...,X_{89}$ is a complete list of representatives of the
extended gene types. 

For $N \leq 13$ let $\Gamma_0^N$ denote the first $2^N$ segments of
$\Gamma_0$, starting from $(0,0)$ and moving in the direction
of the positive quadrant.  We call
$N$ {\it sufficiently large\/} if each extended
gene type has a representative on $\Gamma_0^N$.
We show by a direct computation that
$N=10$ is sufficiently large.  Once we have
our list of $89$ gene types it is a completely
straightforward matter of checking that they
all occur on $\Gamma_0^{10}$
{\it A forteriori\/}, the
polygonal arc $\Gamma_0^{10}$ contains
each of our genes $A_1,...,A_{75}$. 

Our inflation generator gives rise to the inflation
structure $\chi$.  Our goal is to prove that
$\chi$ is coherent.   Here is the strategy.
Given any gene $A \subset \Gamma_0^{10}$, we
let $a$ be the center point of $A$.   We know
that $A$ has the same gene type as some gene
$A_j$ in our inflation generator.  The center
point $a_j$ of $A_j$ is such that $\Phi(a_j)$ is
close to the point $\beta_j(a_j)$.  We record
the gene type of the gene centered at $\beta_j(a_j)$
and then find a point $a' \in \Gamma_0^{13}$ which
is near $\Phi(a)$.  We then verify by direct
computation that $\beta_j(a)=a'$.

Now we know that $\chi(A) \subset \Gamma_0$.  We
check explicitly that the points
$a'_1,...,a'_{1024}$ occur in order on
$\Gamma_0^{13}$. This means that
the two strands $\chi(A_1)$ and $\chi(A_2)$ overlap
whenever $A_1$ and $A_2$ are consecutive genes
on $\Gamma_0^{10}$.    But then
$\chi(A_1) \cup \chi(A_2)$ is a polygonal arc.
This works for all consecutive genes on
$\Gamma_0^{10}$, including the $89$ we
need for the Coherence Lemma.   Applying
the Coherence Lemma, we see that $\chi$ is
coherent.  See \S 7 for more details.
\newline

\subsection{Checking the Completeness of the Lists}

We compute explicitly that $\Gamma_0^{10}$ contains $75$
gene types and $89$ extended gene types.  We also compute
explicitly that $\Gamma_0^{13}$ has $75$ gene types and
$89$ extended gene types.  That is, when we go out
$8$ times as far, we see no new genes or gene types.
We also recall that $\Gamma_0^{10}$ is shadowed by a
subset of $\Gamma_0^{13}$.   This means that the
process of replacing $A_j$ by $\chi(A_j)$, for each
of $j=1,...,75$ produces no new gene types and
no new extended gene types.  Iterating, and
applying induction we see that our list of
gene types and extended gene types is complete,
in the forwards direction.  We then make all the
same checks in the backwards direction.
\newline
\newline
{\bf Remark:\/}
We could take a different approach to the
completeness of our list of genes.
Figure 4.1 shows a partition of the grey polygons
of $\cal P$ into $75$ light grey polygons and $24$
dark grey polygons.  As we mentioned above, the
dark grey polygons correspond to points in
the arithmetic graph contained on closed
$5$-gons or closed $7$-gons.  Our list of
genes is complete because the remaining part of
$\cal P$ is completely partitioned by the
$75$ corresponding dynamical polygons.   Any additional
gene type would have a dynamical polygon that
overlaps with one of the ones we already have
plotted.   We do not insist on this approach
because it requires an analysis of the $24$
dark gray polygons in the picture, something
we have not attempted.  A similar picture would
reveal the completeness of the list of
$89$ extended gene types.

\newpage

%% file: 5backwards.tex
\section{Proof of Theorem \ref{aux}}

\subsection{Existence of the Orbits}

Under our affine map $S' \to S$, discussed in \S 2.1, the
Cantor set $C$ from Theorem \ref{aux} is the Cantor set
on the line $\{y=-1\}$ with
similarity constant $\phi^{-3}$ and endpoints $(0,-1)$ and $(2\phi^{-3},-1)$.
For ease of discussion we identify $C$ with a subset of $\R$ by
dropping the second coordinate.

\begin{lemma} The orbit of every point in $C^*$ is entirely defined.
\end{lemma}

\startproof
 Lemma \ref{exist}, in the case $A=\phi^{-3}=2\phi-3$, shows that
the outer billiards map is defined on and preserves
$(\R-2\Z[\phi]) \times \Z_{\rm odd\/}$.  Hence, 
it suffices to prove that
$C^* \cap 2\Z[\phi]=\emptyset$.

We argue by {\it descent\/}.
Suppose $x=a+b \phi \in C^*$, with $a,b \in 2\Z$. We must have $a \not = 0$.
We take $|a| \geq 2$ as small as possible. 
If $a<-2$ then 
$$x'=2\phi^{-3}-x = (-6-a)+(4-b)\phi \in C^*$$ has $|a'|<|a|$.
Hence $a \geq -2$.  If $a=-2$ then we must have $b>0$.  But
$-2+2\phi$ is already too large to lie in $C^*$. Hence $a \geq 2$.
\newline
\newline
{\bf Case 1:\/}
Suppose $x$ lies in the left half of $C^*$.
Then $\phi^3x \in C^*$. We compute
$$\phi^3x=
(a+2b)+(2a+3b)\phi.$$
By minimality, $|a+2b| \geq a$.  If
$a+2b \geq a$ then $b \geq 0$ and $a+b\phi>1$.  Hence $x \not \in C$.
If $a+2b \leq -a$ then
$b \leq -a$. Hence $a+b\phi<0$ and $x \not \in C$.
\newline
\newline
{\bf Case 2:\/}
Suppose $x$ lies in the right half of $C^*$.
The map $\alpha(x)= 2-\phi^3 x$ maps the right half
of $C^*$ back into $C^*$.   We compute
$$\alpha(x)=(2-a-2b)+(2-2a-3b)\phi.$$
By minimality $|2-a-2b| \geq a$.  If
$2-a-2b \geq a$ then $a+b \leq 1$.  Since $a$ and $b$ are
even, this forces $a+b \leq 0$.   But then $a+b\phi<0$ and
$x \not \in C^*$.
If $2-a-2b \leq -a$ then $b>0$ and $a+b\phi>2$.  Again
$x \not \in C$.
\newline

This takes care of all the cases.
\endproof

\subsection{A Gap Phenomenon}
\label{further}

Let $\Gamma_-$ denote the backwards half of $\Gamma_0$.
We want to see that $\Gamma_-$ rises away from $\partial H$
and then back to $\partial H$ infinitely often (in a
precise way).   Our main idea is to observe that some
initial portion of $\Gamma_-$ rises and falls, and then
to propagate this property using the self-similarity.
One worry is that the approximate nature of 
the self-similarity of $\Gamma_-$ causes the ``lowest points
of approach'' to drift away from $\partial H$.
A {\it gap phenomenon\/} comes to the rescue.

Let
$$\psi(x)=\bigg[\Big(\frac{x}{2\phi},\frac{x}{2}\Big)\bigg].$$
Equation \ref{map2} says that $\Psi=\psi \circ T$.   Recall that $\P_1,...,\P_{26}$ are the
open polygons comprising the partition $\cal P$ of $\T^2$.  The polygon
$\P_3$ is the parallelogram at the bottom left of Figure 2.5 and
$\psi([0,2 \phi^{-3}])$ is the long diagonal of
$\P_3$. 

\begin{lemma}
Suppose that $x \in (0,16)$ and $\psi(x) \in \P_3$.  Then
$x \in (0,2\phi^{-3})$.
\end{lemma}

\startproof
Looking at Figure 2.5 we see that $\P_3$ is the little parallelogram
in the bottom left corner.  The corner vertex of $\P_3$ is $(0,0)$
and the opposite vertex is $(\phi^{-4},\phi^{-3})$.  So, if
$\psi(x) \in \P_3$ then $\delta(x/2) \in (0,\phi^{-3})$.
Hence $x \not \in [2 \phi^{-3},2]$.   It is now an easy exercise to
check that the segment $\psi([2,16])$ is disjoint from $\P_3$.
One way to do this exercise is to plot the $\psi$-image of
the $2000$ maximally and
evenly spaced points in the
interval $[2,16]$ and observe that none of them is within $.0001$ of
$\P_3$.  Since $\psi$ is $1$-Lipschitz, our plot guarantees that  no point
of $[2,16]$ lies in $\P_3$.  We omit the details.
\endproof

Given a point $p \in H$, the halfplane containing the arithmetic graph,
let $v(p)$ denote the length of the vertical line segment connecting
$p$ to the line through the origin parallel to $\partial H$.  For
instance $v(0,m)=m$.  (The line $\partial H$ lies slightly below the
origin.)

\begin{corollary}[Gap Phenomenon]
Suppose $p$ is a point of the arithmetic graph having type $3$.
If $v(p)<7$ then $v(p)<\phi^{-3}$.
\end{corollary}

\startproof
Note that $T(0,7)<15$, and the fibers of $T$ are parallel to $\partial H$.
Hence $T(p) \in (0,15)$.  By definition, $\psi(T(p)) \in \P_3$.  Hence
$T(p) \in (0,2\phi^{-3})$.  But $T(0,\phi^{-3})>2 \phi^{-3}$. In other
words, if $v(p')=\phi^{-3}$ then $T(p')>2 \phi^{-3}$. 
Hence $v(p)<\phi^{-3}$.
\endproof

\subsection{Infinite Return}

In this section we prove that $O_-(p)$ returns to every neighborhood
of $p$ infinitely often.  We think of this as a warm-up to the
proof of Theorem \ref{aux}.
Let $A_0$ be the gene through the origin.
Let $A_0'$ be the strand that shadows it.  Figure
5.1 shows the picture.

\begin{center}
\psfig{file=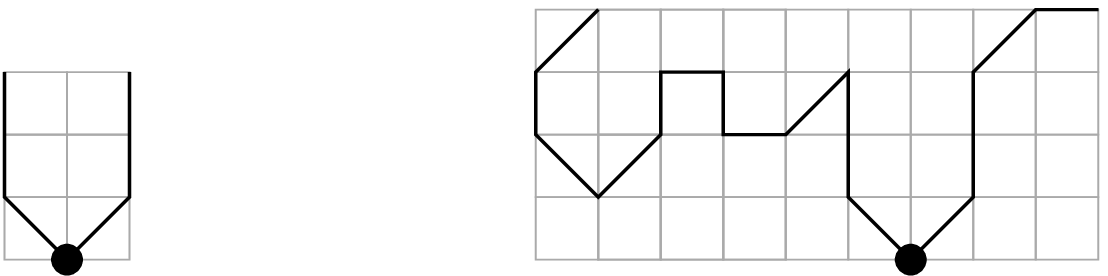}
\newline
{\bf Figure 5.1:\/} The gene $A_0$ and its shadow
$A'_0$.
\end{center}  

The dynamical polygon $P_{0}$ associated to $A_{0}$ is contained in
$\P_3$, because the center of $A_{0}$ has type $3$.  Thus $P_{0}$
is a small tile in the bottom left corner of Figure 2.5.  Here we
show a plot, and describe all the plotted objects
in order of decreasing size.  Each polygon is contained in the
previous one.  Referring to Figure 5.2:

\begin{itemize}
\item  The big $L$ is the corner of our fundamental domain
for $\T^2$. 
\item  The big parallelogram is $\P_3$. 
\item  The big triangle is
$P_{0}$. 
\item  The small parallelogram is
$P'_{0}$, the polygon
associated to $A'_{0}$. 
\item The small triangle is $\gamma(P_0)$, where
$\gamma$ is the special similarity fixing $\Psi(0,0)$.
\end{itemize}

\begin{center}
\psfig{file=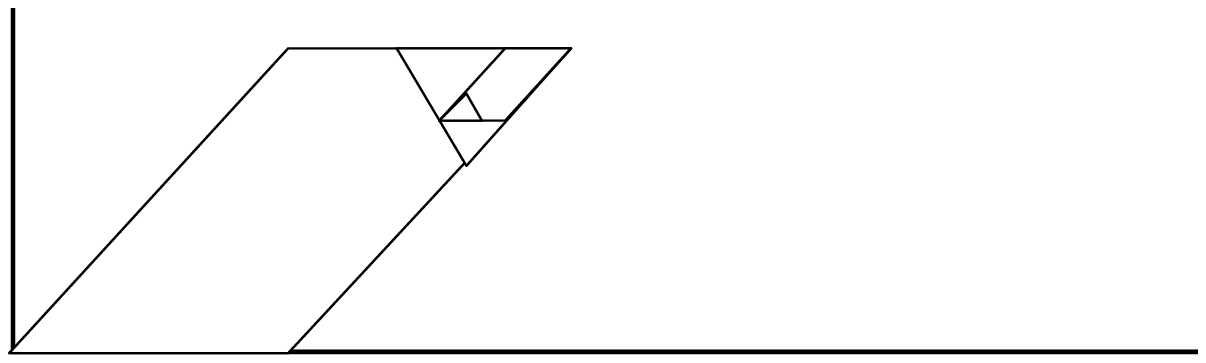}
\newline
{\bf Figure 5.2:\/} The associated dynamical polygons.
\end{center}  

Let $A$ be some occurrence of $A_{0}$ on $\Gamma_-$.  Let
$a$ be the center point of $A$.  The first site
is $(-21,5)$.   Let $A'=\chi(A)$ be the
strand that shadows $A$.  Let $a'$ be the center point of $A'$.  
Our analysis shows that $\Psi(a)$ and $\Psi(a')$ both lie in $P_{0}$.
The point is that $P'_{0} \subset P_{0}$.    Let's
assume by induction that $v(a)<\phi^{-3}$.  Then, by the
triangle inequality, and the definition of our inflation
structure,
$$v(a')<\phi^3 \times \phi^{-3} + 4=5.$$
Since $a'$ also has type $3$, our Corollary above now says that
$v(a')<\phi^{-3}$.
We can now start with the copy of $A_0$ based at $a'$.  This produces
a further point $a''$, centered on a copy of $A_0$, such that
$v(a'')<\phi^{-3}$.  Thus we produce distinct  points
$a,a',a'',... \in \Gamma_-$ all
within $\phi^{-3}$ of $\partial H$.
Let $a^{(n)}$ be the $n$th point produced by
this process and
let $x^{(n)}=T(a^{(n)})$.
The points $x^{(n)}$ all belong to $O_-(p)$.

By construction
$$\psi(x^{(n+1)})=\gamma(\psi(x^{(n)})); \hskip 30 pt
x^{(n)}, x^{(n+1)} \in (0,2 \phi^{-3}) \subset (0,1/2).$$
These two conditions imply that
$x^{(n+1)}=\gamma(x^{(n)})$.  Hence
$x^{(n)}=\gamma^n(x)$, where
$x=T(-21,5)$.  But $\gamma$ is a contraction
fixing $p$.  Hence $x^{(n)}$ converges to $p$.

\subsection{Cantor Set Structure}

Our argument above came from considering just the gene $A_0$. We get
the Cantor set structure by considering the action of two genes.
Notice that $A'_0$ contains two genes having a core of type $3$.
Figure 5.3 shows the other gene $A_1$ with this property.

\begin{center}
\psfig{file=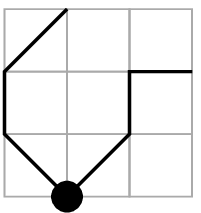}
\newline
{\bf Figure 5.3:\/} The gene $A_1$.
\end{center}  

Let $A_1'$ be the strand that shadows $A_1$.  It turns out that
$A_1'$ and $A_0'$ have the same type. Moreover, the center point
$a_1'$ of $A_1'$ is the center of another copy of $A_0$.  Figure
5.4 shows the corresponding dynamical polygons.  This time
$P_1$, the little triangle in the bottom left corner, is
disjoint from $P'_1$.  The map $\gamma$ is the same one as
for $P_0$, and the parallelogram is again $P'_0=P'_1$.
The tiny triangle inside $P'_1$ is $\gamma(P_1)$.

\begin{center}
\psfig{file=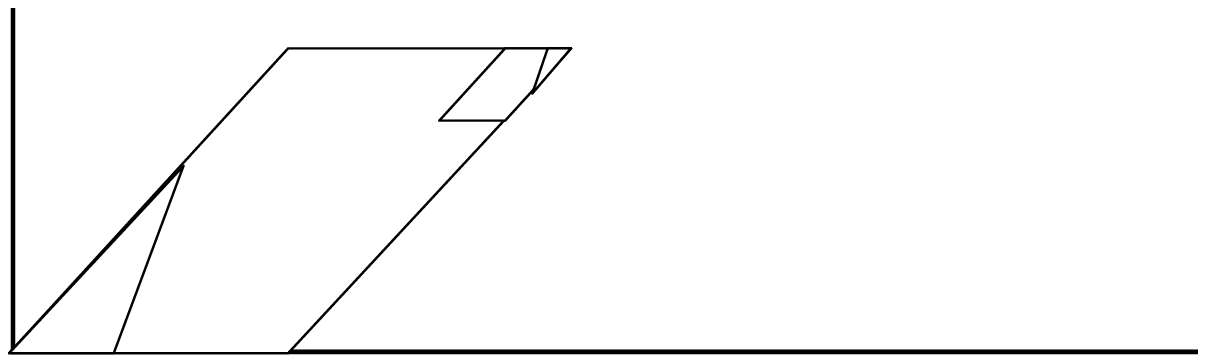}
\newline
{\bf Figure 5.4:\/} The associated dynamical polygons.
\end{center}  

We set $K=\P_3$ and $K_0=P'_0$ and $\gamma_0=\gamma$.
There are several geometric pieces of information we now record.
First, $\gamma_0(K)=K_0$.
Second, let $a_0$ and $a_1$ be the two type-$3$ points of
$A_0'$.  For instance, we might take $A'_0$ such that
$a_0=(0,0)$ and $a_1=(-5,1)$. 
From Equation \ref{map2} we compute that
\begin{equation}
V_0:=\Psi(a_0)-\Psi(a_1)=[(5\phi^{-4}-\phi^{-1},5\phi^{-3})]=(\phi^{-4}-\phi^{-7},\phi^{-3}-\phi^{-6}).
\end{equation}
Let $\tau(x)=x-V_0$.  
Geometrically, the translated parallelogram $K_1=\tau(K_0)$ has $(0,0)$ as a vertex
and fits exactly
into the lower left corner of $K$.  
Figure 5.5 shows a schematic picture.
Let $\gamma_1=\tau \circ \gamma_0$.  Then
$\gamma_j(K)=K_j$.  We define $K_{00}=\gamma_0(K_0)$ and
$K_{01}=\gamma_0(K_1)$, etc.  Figure 5.5 shows a schematic
picture of the first few of these sets.

\begin{center}
\psfig{file=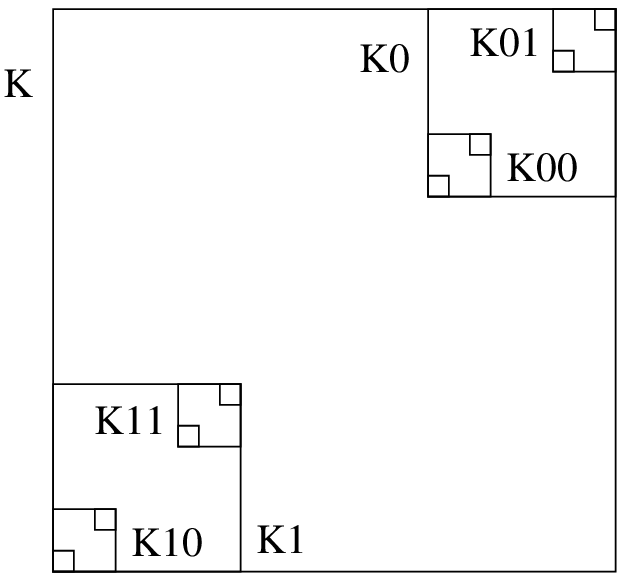}
\newline
{\bf Figure 5.5:\/} A schematic picture
\end{center}  

The reader might expect $K_{00}$ to be on the top and
$K_{11}$ on the bottom, rather than in the middle.  This
twisting of the labelling comes from the fact that the
multiplier of $\gamma$ is $-\phi^{-3}$, a negative number.
The fixed point of 
$\gamma$ is the nested intersection $K_0 \cap K_{00} \cap K_{000}$...
The limit set $C'$ of the semigroup $\langle \gamma_0,\gamma_1\rangle$
is a Cantor set whose endpoints are the far opposite vertices of $K_0$.

Now we can explain the Cantor set structure.
Suppose that $X$ is a type $3$ point $\Gamma_-$ such 
that $X$ is the center of a gene $A$ that is a copy of either $A_0$ or $A_1$.
Suppose also that $v(X)<\phi^{-3}$.
Let $X_0$ and $X_1$ be the two points of $A'$ having type $3$. Our
notation is such that $X_0$ is the center of a copy of the
gene $A_0$ and $X_1$ is the center of the copy of the gene
$A_1$.  The point $X_0$ shadows the dilated point
$\phi^3X$, and the point $X_1$ is several units to the right.
Our arguments above show that
\begin{itemize}
\item $v(X_0)$ and $v(X_1)$ are both less than $\phi^{-3}$.
\item $\Psi(X_1) \subset K_1$ and $\Psi(X_0) \subset K_0$.  
\item $\tau(\Psi(X_0))=\Psi(X_1)$.
\item $\Psi(X_j)=\gamma_j(\Psi(X))$.
\end{itemize}

Now we consider a tree-like induction process.
We start with the point $X=(0,0) \in \Gamma_-$.
Considering the above inflation process, we produce
the points
$$X_0=X; \hskip 30 pt
X_1=(-5,1).$$
Inflating again, we produce the points
$$X_{00}=X_0;
\hskip 15 pt X_{01}=X_1; \hskip 15 pt
X_{10}=(-21,5); \hskip 15 pt 
X_{11}=(-26,5).
$$
Inflating again, we produce the $8$ points,
corresponding to the binary strings of length $3$,
and so on.   
In general, we produce one point $X_{\beta} \in \Gamma_-$ for each
finite binary string $\beta$.  (These points are
not all distinct, e.g. $X_{00}=X_0$.) 
By construction
$$\Psi(X_{\beta}) \subset K_{\beta}.$$
Therefore, the closure of the vertex set of
$\Psi(\Gamma_-)$ contains the
Cantor set $C'$.

Let $I_0=[0,2 \phi^{-3}]$.   We have already remarked that
$\psi(I_0)$ is long diagonal of $\P_3$, the segment that
connects the two endpoints of  $C'$.
Recall that $T$ maps the vertices of $\Gamma_-$
into the backwards orbit $O_-(p)$.  Also
$\Psi=\psi \circ T$.  Therefore
$C'=\psi(C)$, where $C \subset I_0$ is an
affine image of $C'$.  Indeed,
$C$ is the Cantor set from
Theorem \ref{aux}.  (We are identifying
$[0,\infty) \times \{-1\}$ with $[0,\infty)$, as
discussed at the beginning of the chapter.)
We have realized $C$ as a set contained in the closure
of $O_-(p)$.  Hence $O_-(p)$ returns
densely to $C$.

\subsection{$2$-adic Structure}

Since $\psi(I_0)$ is the long diagonal of $\P_3$, the
two conditions
$x \in \Z^2 \cap H$ and $T(x)<2 \phi^{-3}$ imply that
$x$ has type $3$.   We say that such points are
{\it basepoints\/} of $\Gamma$.

We see inductively that our shadowing
construction accounts for all the basepoints of $\Gamma_-$.  
That is, every such point arises on a strand
of $\Gamma_-$ that shadows an inflated gene.
Therefore, we can index the basepoints
by binary strings $\beta$.   Each binary
string $\beta$ represents the integer $n(\beta)$ as usual.
For instance
$$n(11001)=(1,1,0,0,1) \cdot (16,8,4,2,1)=25.$$
Examining our construction we see that the two basepoints
$X_{\beta}$ and $X_{\beta'}$ are equal if and only if
$n(\beta)=n(\beta')$.  Moreover, as we travel along
$\Gamma_-$ away from $(0,0)$ we encounter the points
$X_{\beta}$ in order of the integers they represent!

Let $\theta_2: \Z_2 \to C$ be the homeomorphism from
Theorem \ref{aux}.  By construction, we have
$\theta_2^{-1}(0)=p$.  Let $T: \Z^2  \to [0,\infty)$ be the
map from Equation \ref{map}.   Let $X$ and $X_0$ and
$X_1$ be the points referred to in our shadowing
construction above.
We have already mentioned that
$$\Psi(X_0)=\gamma_j(X); \hskip 30 pt
j=0,1.$$
Given the defining property of $\theta_2$, we have
$$
\theta_2^{-1}(T(X_j))=2\theta_2^{-1}(T(X))+j; \hskip 30 pt
j=0,1.
$$
But then we have
$$
\theta_2(T(X_{\beta}))=n(\beta).
$$
But, by construction $T(X_{\beta})$ is the $n$th point of
$O_-(p) \cap I_0$, which is the same as the
$n$th point of $O_-(p) \cap C$.   Hence $\theta_2$ maps the $n$th point of
$O_-(p) \cap C$ to $n$.  This is a special case of the
first statement of Theorem \ref{aux}.

Now we will deal with the question of return times and excursion
distances for the points of $O_-(n)$.   Given a basepoint
$X$, let $\Gamma_+(X)$ be the forwards portion of
$\Gamma_0$ which starts at $X$. 
Let $X$ and $X_0$ and $X_1$ be as in the previous section, so
that $X_0$ is the basepoint that shadows $\phi^3X$.  

 By
the inflation property, $\Gamma_+(X_0)$ closely follows
$\phi^3\Gamma_+(X)$.  Also, one can see by looking at 
a single example that the initial portion of $\Gamma_+(X_0)$
rises up at least $10$ units from $\partial H$ before
coming back towards $\partial H$.  Moreover, the
 next basepoint (after $X_0$) encountered by
$\Gamma_+(X_0)$ is the one that shadows the dilation
of the next basepoint encountered by $\Gamma_+(X)$.
These properties show that
$\Gamma_+(X_0)$ rises up about $\phi^3$ times as high
as $\Gamma_+(X)$, and takes about $\phi^3$ times as long
to return to the next basepoint as does
$\Gamma_+(X)$.  The next basepoint encountered by
$\Gamma_+(X_0)$ is the one that shadows the dilation
of the next basepoint encountered by $\Gamma_+(X)$.
On the other hand, $\Gamma_+(X_1)$ just travels a
few units to the right before returning to
the basepoint $\Gamma_+(X)$.  

Given any binary string $\beta$, let $\nu(\beta)$ denote
the number of $0$s on the right of $\beta$.  For
instance $\nu(11000)=3$.  Equivalently,
$\nu(\beta)$ equals the highest power of
$2$ dividing $n(\beta)$.   Given $X=X_{\beta}$
let $\nu=\nu_{\beta}$. 
Applied inductively,
our arguments show that
$\Gamma_+(X)$ rises up roughly $\phi^{3\nu}$ units
before returning to the next basepoint after 
roughly $\phi^{3\nu}$ units of time. Here
$X=X_{\beta}$.  Hence the excursion
distances and return times for the forward orbit
of the $n$th point $x_n$ of $O_-(p) \cap C$ are proper
functions of the $2$-adic distance from
$\theta_2^{-1}(x_n)$ to $0$.   

By simply reversing the direction of $\Gamma_+(X)$ we
see that the excursion distances and return times for
the forwards orbit of $x_n$ are proper functions of
the $2$-adic distance from
$\theta_2(x_n)$ to $-1$.  The point here is that
$\theta_2$ maps $n$ and $n-1$ respectively to the
endpoint of $\Gamma_+(\beta)$ and the first
basepoint it encounters when travelling to the right.

There is one more observation we want to make.
If $X$ is a basepoint and $\nu(X)$ is very large,
then a very long initial portion of $\Gamma_+(X)$
looks exactly like $\Gamma_+(0,0)$, the forward
portion of $\Gamma_0$ that starts at $(0,0)$.
This is because both strands are ``produced''
by many iterates of the same inflation process.
An equivalent way to see this is that any finite
portion of $\Gamma_0$ is determined by some
small dynamical polygon in $\T^2$ containing
$\Psi(0,0)$.  If $X$ is some point with
$\nu(X)$ large, then $\Psi(X)$ belongs to this
same dynamical polygon. 

\subsection{Extension to the Cantor Set}

Now we put everything together and finish the proof of
Theorem \ref{aux}.

\begin{lemma}[Rising]
For each positive constant $K$ and each lattice point
$X \in \Gamma_-$, there exists a constant
$K'=K'(K,X)$ with the following property:
If we take $K'$ steps in either direction along
$\Gamma_-$, starting at $X$, then we rise at least
$K$ units above $\partial H$.
\end{lemma}

\startproof
Consider the forwards direction.  The backwards direction is similar.
It this lemma is false then 
we can encounter a long string of basepoints $\{X_i\}$ such that
$\nu(X_i)$ is always small.  But the map
$x \to x-1$, when iterated, brings any point very close
$2$-adically to $0$ within a uniform number of steps.
Hence, we don't have to walk very long from $X$ before we
hit a basepoint with high $\nu$-value, and then we rise
up steadily away from $\partial H$ for a long time.
This comes from the fact that a long initial portion of
$\Gamma_+(X)$ is a translate of a long initial portion
of $\Gamma_+(0,0)$, as we remarked above.
\endproof

\begin{lemma}
For any $y \in C^*$ and any $K>0$ there are
positive integers $n_+=n_+(y,K)$ and
$n_-=n_-(y,K)$ such that the $(n_+)$th point of
$O_+(k)$ and the $(n_-)th$ point of $O_-(y)$ are at both least $K$ units
from the origin.
\end{lemma}

\startproof
We can find a sequence $\{x_n\} \in O_-(p)$
converging to $y$.   For $n$ sufficiently large,
the first $K'$ iterates of $x_n$ and $y$ will have
the same combinatorial structure.  Here $K' \to \infty$
as $n \to \infty$.   But then the first $K'$ iterates
of $y$ remain (say) within $1$-unit of the corresponding
$K'$ iterates of $x_n$.  Letting $K'$ be as in the
Rising Lemma, we see then that some point of
$x_n$ rises up $K$ units and therefore some iterate
of $y$ rises up at least $K-1$ units.  Since
$K$ is arbitrary, the orbit of $y$ (in either
direction) is unbounded.
\endproof

We continue with the notation from the lemma.  Let
$x_{n,-}$ denote the forward return to $C^*$ of the point $x_n$.
As long as $y \not = p$, there is some uniform $K$
such that $x_n$ and $x_{n,-}$ are separated by at
most $K$ units.  Here $K$ depends on $y$ but not on $n$.
This we see that $y_-$ and $y$ are separated by at most
$K$ units.  For $n$ large enough the orbit of 
$x_n$ has the same combinatorial structure as the
orbit of $y$ for the first $K$ iterates.  Hence
$y_--y=x_{n,-}-x_n$.   Therefore, by continuity,
$\theta_2^{-1}(y_-)=\theta_2^{-1}(y)-1$.  
The same continuity argument takes care of the statement
about the return times and excursion distances for
the forward orbit.  
This proves the first statement of Theorem \ref{aux}, and
the second statement has
essentially the same proof.

\newpage

%% file: 6return.tex
\section{Proof of the Arithmetic Graph Lemma}

\subsection{The Computational Evidence}

Before we launch into the proof of the Arithmetic Graph Lemma
we discuss how \footnote{The directness of
this discussion is 
misleading.  We arrived at the ideas here only after having
exhausted every half-baked
and useless scheme we could think up.}
 we discovered that it was true. Our proof of
the Arithmetic Graph Lemma is logically independent from
the way we discovered it, but the computations we
made during the discovery process will serve
us below at a certain step in our proof, namely the
step where we show 
that Lemmas \ref{ag} and \ref{ag2} below are equivalent
to the Arithmetic Graph Lemma.

Recall that $\Gamma \subset H$, a certain halfplane in
$\R^2$.  We thought of $\Gamma$ as a union of
paths of physical particles, and we got the 
idea to partition $H$ into parallel bands, each of
which has width $\phi$.  We thought of these bands as
semi-permeable membranes, letting various particles
pass through and deflecting others.  Experimentally,
we noticed that the permeability of the $n$th band
depended somehow on the decimal part of $n/\phi$.  In particular,
the bands indexed by Fibonacci numbers are the least
permeable and behave the most like mirrors.
Thus, the quantity $\delta(n/\phi)$, an
$\R/\Z$ valued coordinate, seemed relevant
to the ``physics'' of a particle in the $n$th band.
Here $\delta(x)$ is the decimal part
of $x$. 

The other quantity that seemed to influence the ``physics''
of a particle
was the relative position of the particle within
the band that contained it.  Like our first coordinate,
this is another coordinate that takes values in
$\R/\Z$.  Thus, the local properties of $\Gamma$
seemed to be determined by two coordinates in
$\R/\Z$, which we considered as a point in
$\T^2$.   When we wrote down the formulas
for these coordinates and tuned them, we arrived at our map
$\Psi$ from Equation \ref{map2}.

Once we had $\Psi$ we performed an experiment.  We
took some huge sample of $\Gamma$, and fixed some
local type.  Sampling lattice points $p$,
we plotted $\Psi(p)$ if and only
if $p$ had that local type.  We noticed that the resulting
plot was very dense in either a single convex polygon,
or else the union of several convex polygons, depending
on the type.  (Later we indexed the types so that each type
corresponds to a single polygon.)  It
was easy to guess the vertices of the polygons from
the plots, and
this is how we determined the partition $\cal P$.
Indeed, Billiard King
draws $\Gamma$ in two ways, using the partition and using
the dynamics, and the picture is the same.
Now we turn to the proof.

\subsection{Factoring the Return Map}

Say that a {\it strip\/} is a region $A \subset \R^2$ bounded by $2$ parallel lines,
$\partial_0 A$ and $\partial_1 A$.  Let $V$ be a vector such that
$\partial_0A+V=\partial_1 A$. 
Given the pair $(A,V)$ we (generically) define a map
$E: \R^2 \to A$ by the formula
$E(x)=x-nV$ where $n$ is the unique integer such that $E(x) \in A$.
This map is well defined unless $x$ lies in a discrete infinite
family of parallel lines.  

There is a unique affine functional 
$f(x,y)=a_1x+a_2y+a_3$ such that
$f_L(V)=1$, and
$f(x,y) \in (0,1)$ iff $(x,y) \in A$.
Here $f_L(x,y)=a_1x+a_2y$ is the linear
part of $f$.  Given $f$ we have the following
explicit formula:
\begin{equation}
\label{floor}
E(p)=p-{\rm floor\/}(f(p))\ V.
\end{equation}
Equation \ref{floor} is defined unless
$f(p)$ is an integer.  We say that $\alpha=(a_1,a_2,a_3)$ is
the {\it triple associated to\/} $(A,V)$.

\begin{center}
\psfig{file=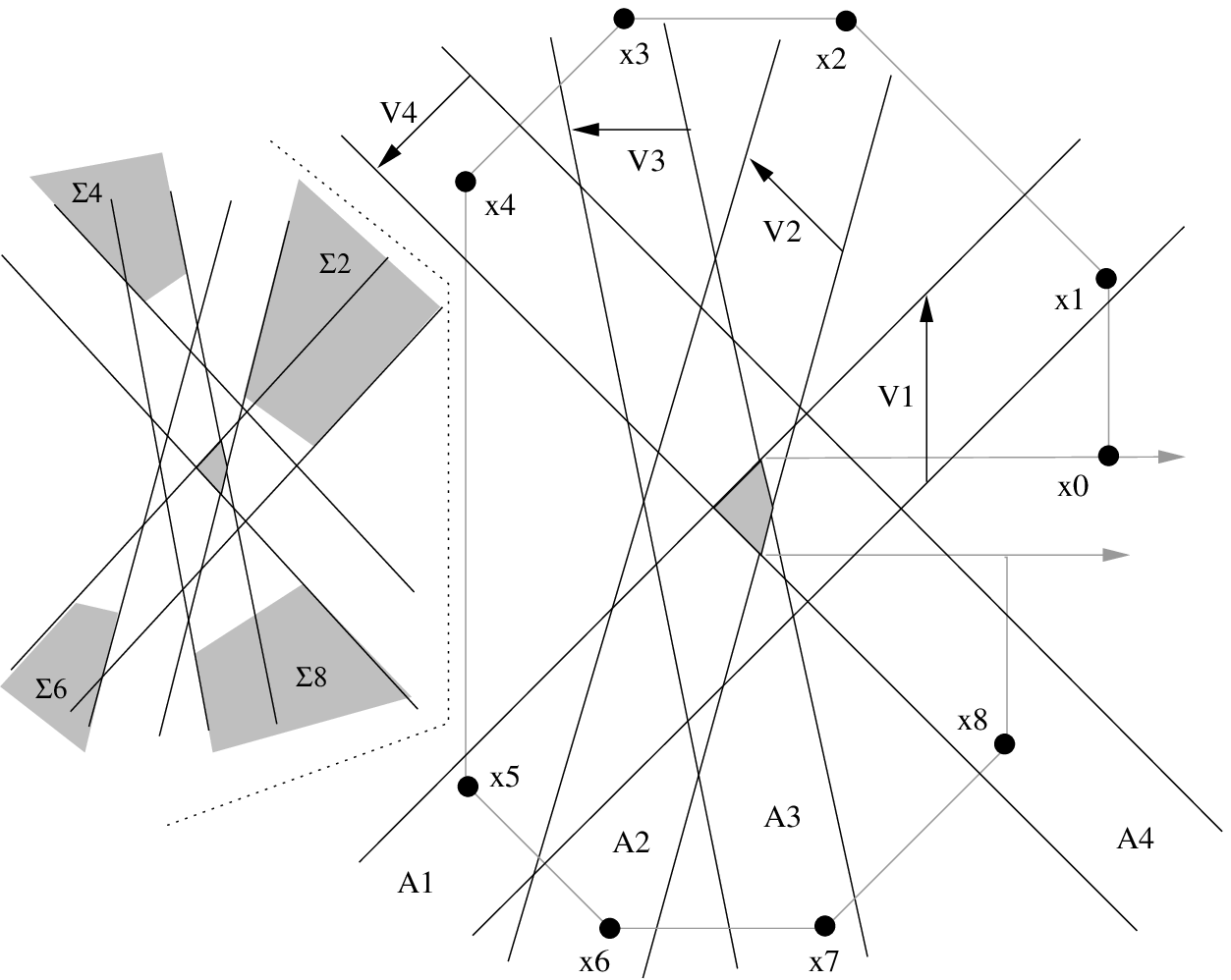}
\newline
{\bf Figure 6.1:\/} The pinwheel
\end{center}

Each strip $A_j$ in Figure 6.1 is obtained by extending an edge of
the side of $S$, and then rotating this extended side through one
of the vertices that does not contain the edge.  Each vector
$V_j$ is twice the difference between a pair of vertices on
the kite.  Let $E_1,E_2,E_3,E_4$ be the corresponding
strip maps.  Here is the data associated
to these maps:
\newline
\noindent
\begin{tabular}{ll}
\\
 $\alpha_1=(-1/4,+1/4,+3/4)$ & $V_1=(0,4)$ \\ 
 $\alpha_2=(-\phi/4,+1/2-\phi/4,+1/2-\phi/4)$ & $V_2=(-2,+2)$ \\
 $\alpha_3=(-\phi/4,-1/2-\phi/4,+1/2-\phi/4)$ & $V_3=(+4-4\phi,0)$. \\
 $\alpha_4=(-1/4,-1/4,+3/4)$ & $V_4=(-2,-2)$. \\
\\
\end{tabular}
\newline
Let $E_{j+4}=E_j$ for $j=1,2,3,4$.

\begin{lemma}[Pinwheel]
Let $x_0 \in {\cal C\/}(\pm)$. Let  $x_j=E_j(x_{j-1})$.
Then $x_8$ and $\Upsilon_R(x_0)$ lie
on the same vertical line. 
\end{lemma}

\startproof
Let $V_{j+4}=-V_j$ for $j=1,2,3,4$.
For any $x$ on which $\Upsilon$ is defined, we have
$\Upsilon(x)-x=V_x$, where $V_x$ is twice one of the vectors
pointing from one vertex of $S$ to another.
The vectors $V_1,...,V_8$ have this form.

The lines comprising our strips divide $\R^2$ into a finite number
of bounded regions and a finite number of unbounded regions.
Let $K$ denote the closure of the union of bounded regions.
One can check, with a little experimentation, that the
finitely many iterates $p,\Upsilon(p),...,\Upsilon_R(p)$
avoid the set $K$ provided that $p=(x_0,\pm 1)$ and
$x_0>4$.   We will consider such points first.

The complement of $K$ is divided into 
$8$ sectors $\Sigma_1,...,\Sigma_8$. Each sector
$\Sigma_j$ is bounded by one line of $A_{j-1}$
and one line of $A_j$, and contains a non-compact subset
of $A_{j-1}$.  The left hand side of Figure 6.1 
indicates half of these sectors.
It is easy to check that
$q \in \Sigma_j$ implies $V_q=V_j$.
This result immediately implies our lemma for points
$(x_0,\pm 1)$ where $x_0>4$ and all maps are defined.

When $x_0<4$ we can just check the few cases by hand.
Indeed, to be sure we sampled $1$ million evenly spaced
points in the two intervals $(0,10) \times \{\pm 1\}$.
\endproof

The Pinwheel Lemma lets us {\it factor\/} the return map.
The composition $E_8...E_1$ maps 
$x_0 \in {\cal C\/}(\pm)$ to $x_8 \in \R^2$,
and the $y$-coordinate of $x_8$ is $\pm 1+4k$.
In each case, $\Psi_R(x_0) \in {\cal C\/}(\pm)$. Hence
\begin{equation}
\label{mu}
\Upsilon_R=\zeta \circ (E_4 \circ ...\circ E_1)^2; \hskip 50 pt
\zeta(x,\pm 1+4k)=(x,\pm 1).
\end{equation}

\subsection{Four Dimensional Compactificaton}

We say that a strip map $E$, with associated triple
$(a_1,a_2,a_3)$ and vector $V=(v_1,v_2)$, is {\it special\/} if
\begin{itemize}
\item $a_1,a_2,a_3 \in \frac{1}{4}\Z[\phi]$ and $a_1 \pm a_2 \in \frac{1}{2}\Z[\phi]$.
\item $V_1, V_2 \in 2\Z[\phi]$ and $v_1 \pm v_2 \in 4\Z[\phi]$.
\end{itemize}
By inspection, the strip maps from Equation \ref{mu} are all special.

Let $T^4_8=(\R/8\Z)^4.$
Consider the following
embedding of $\R^2$ into $T^4_8$:
\begin{equation}
\label{wide}
\widetilde \psi(x,y)=
\big[(x+y,x-y,(x+y)/\phi,(x-y)/\phi)\big]_8.
\end{equation}  
Here $[v]_8$ denotes the image of $v \in \R^4$ in
$T^4_8$. 
We only care about the next result for our $4$ maps, but the argument is easy
to make in general.
We give explicit formulas for the extensions
of our $4$ maps in \S \ref{partitioncalc}.

\begin{lemma}[ Extension]
\label{lipschitz}
Let $E$ be a special strip map. 
There is a finite
union $Y$ of flat $3$-dimensional tori
in $T^4_8$ together with a
map $\widetilde E: T_8^4-Y \to T_8^4$ such
that $\widetilde \psi \circ E=\widetilde E \circ \widetilde \psi$
whenever both maps are defined.  
$\widetilde E$ is locally affine on each component
of $T_4^8-Y$ and the linear part of 
$\widetilde E$ is defined over $\Z[\phi]$ and independent of component.
\end{lemma}

\startproof
We fix  $p_1 \in \R^2$ on which $E$ is defined,
and consider points $p_2$ such that
$\widetilde \psi(p_1)$ and $\widetilde \psi(p_2)$ are
close together in $T_8^4$.  Setting $p=p_2-p_1$, this means
that the coordinates of $\widetilde \psi(p)$ are
all near $0$.  We make this assumption about the
coordinates of $\widetilde \psi(p)$ and treat these
coordinates as
small real numbers.

Let $p=(x,y)$. 
Given $\tau_j \in \Z$ for $j=1,2,3,4$ we write
\begin{equation}
\langle\tau_1,\tau_2,\tau_3,\tau_4\rangle=
\tau_1 \bigg(\frac{x+y}{8}\bigg)+
\tau_2 \bigg(\frac{x-y}{8}\bigg)+
\tau_3  \bigg(\frac{x+y}{8\phi}\bigg)+
\tau_4  \bigg(\frac{x-y}{8\phi}\bigg)
\end{equation}
The coefficients of $\tau_j$ are nearly integers.  Hence
$\delta(\langle \tau_1,...,\tau_4 \rangle)$ is near $0$, and also
is a $\frac{1}{8}\Z$-linear combination of the coordinates of $\widetilde \psi(p)$.
We call this the {\it reduction principle\/}.
Also,
$\phi\langle \tau_1,\tau_2,\tau_3,\tau_4\rangle=\langle \tau'_1,\tau_2',\tau_3',\tau_4'\rangle$
because $\phi=1+\phi^{-1}$.  Call this the {\it multiplication principle\/}.

Let $f$ be the functional associated to $E$.
Since $\Z[\phi]=\Z[\phi^{-1}]$ we have
$$a_1=\frac{1}{4}s_1+\frac{1}{4\phi}t_1; \hskip 30 pt
a_2=\frac{1}{4}s_2+\frac{1}{4\phi}t_2; \hskip 15 pt
s_1,t_1,s_2,t_2 \in \Z.$$
The special property of our map implies that
$u_{\pm}=\frac{1}{2}(u_1 \pm u_2) \in \Z$ for $u=s,t$. We have
\begin{equation}
\label{maincalc}
f(p_2)-f(p_1)=
\bigg(s_1\frac{x}{4} + s_2 \frac{y}{4}\bigg)+
\bigg(t_1 \frac{x}{4\phi}+ t_2 \frac{y}{4\phi}\bigg)=
2\langle s_+,s_-,t_+,t_-\rangle.
\end{equation}
From equation \ref{maincalc} and the reduction principle,
$f(p_1)-f(p_2)$ is extremely near an integer, and
much closer to an integer than $f(p_1)$.  Hence
$f(p_2)$ is not an integer and
$E$ is defined on $p_2$.

Let $N_j$ be the floor of $f(p_j)$. 
The properties of $f(p_1)-f(p_2)$ and $f(p_j)$ just mentioned imply that
$N_2-N_1$ is the integer closest to $f(p_2)-f(p_1)$.
Equation \ref{maincalc} gives us
$N_2-N_1=2K$ with 
\begin{equation}
\label{Kval}
K=\beta-\delta(\beta) \in \Z; \hskip 20 pt
\beta=\langle s_+,s_-,t_+,t_-\rangle.
\end{equation}

Let $V=(v_1,v_2)$ be the vector associated to $E$.
Equation \ref{floor} yields
\begin{equation}
E(p_2)=E(p_1)+p-2KV.
\end{equation}
To prove the existence of $\widetilde E$ it suffices to show that
each coordinate of 
$\widetilde \psi(2KV)$ is a $\Z[\phi]$ linear
combination of the coordinates of $\widetilde \psi(p)$.
Let $2KV=(w_1,w_2)$.
By hypotheses $v_1 \pm v_2 \in 4\Z[\phi]$.
Hence $\phi^a(w_1 \pm w_2) \in 8\Z[\phi]$ for any
exponent $a \in \Z$.  (We only care about $a=0,-1$.)
In light of this observation, it suffices to show that
$\delta(K\phi)$ is near $0$ and also a
$\frac{1}{8}\Z[\phi]$
linear combination of the coordinates of $\widetilde \psi(p)$.

Equation \ref{Kval} gives
\begin{equation}
\delta(K\phi)=
\delta(\phi \beta-\phi \delta(\beta))=
\delta(\phi \beta)-\delta(\phi \delta(\beta))=^*\delta(\phi \beta)-\phi\delta(\beta).
\end{equation}
The starred inequality comes from the fact that $\delta(\beta)$ is near $0$.
From the reduction principle and the real multiplication principle we see that
both $\delta(\beta)$ and $\delta(\phi \beta)$ are near $0$ and
$\frac{1}{8}\Z$ linear combinations of the coordinates
$\widetilde \psi(p)$.  This proves what we want about the coordinates
if $\widetilde \psi(2KV)$.

It only remains to find the domain of definition for $\widetilde E$ on
$T_8^4$.
The map $E$ is defined on the complement of an infinite
union of evenly spaced lines in $\R^2$.  These
lines consist of points $(x,y)$ satisfying the integral linear equations
$\langle s_3,s_4,t_3,t_4\rangle \in \Z$.   Applying the map
$\widetilde \psi$ we see that the image of these
points is contained in a finite union of flat embedded $3$-tori.
Our construction shows that $\widetilde E$ is defined and locally
affine on each component of $T_8^4-Y$.
\endproof

\subsection{Geometry of the Return Map}

We begin with a consequence of the Pinwheel Lemma.  Let
$\pi_1:\R^2 \to \R$ be projection onto the first coordinate.

\begin{lemma}
\label{finiteness}
The map $\pi_1(\Upsilon_R(x,1))-x$ only takes on finitely
many values.
\end{lemma}

\startproof
To see this note that the ``octagonal spiral''
$\Sigma$ in Figure 6.1 connecting the points
$x_0,...,x_8,\zeta(x_8)$ is within a bounded distance of being symmetric
about the origin.  In other words, the image of $\Sigma$ under
reflection through the origin lies within a uniformly bounded
tubular neighborhood of $\Sigma$.  In terms of the displacement
vectors $x_1-x_0,x_2-x_1,...,$ we see that the opposite sides of
$\Sigma$ cancel, up to a uniformly bounded error.  Hence
there are uniformly bounded integers $\alpha_j$, depending
on $(x,1)$, such that
$\Upsilon_R(x,1)-(x,1)=\sum \alpha_j V_j$.
\endproof

We now return to the theme taken up in the previous section.
The analog of the Extension Lemma is not quite
true for the map $\zeta$ in Equation \ref{mu} but
nonetheless a similar result holds.
The domain for 
$\zeta$ is the union $L$ of horizontal lines whose $y$-coordinates
are odd integers.  
Consider the map $\psi: {\cal C\/}(\pm) \to \T^2$ given by
\begin{equation}
\label{psimap}
\psi(x,\pm 1)=[(x/(2\phi),x/2)]
\end{equation}

\begin{lemma}
\label{zetastar}
There is a locally affine map
$\widetilde \zeta: T_8^4 \to \T^2$ such that
$\widetilde \zeta \circ \widetilde \psi=
\psi \circ \zeta.$
\end{lemma}

\startproof
We define $\widetilde \zeta: \R^4 \to \T^2$ by the formula
$$\widetilde \zeta(x_1,x_2,x_3,x_4)=\bigg[\frac{x_3+x_4}{4},\frac{x_1+x_2}{4}\bigg].$$
The image in $\T^2$ is unchanged if we add multiples of $8$ to the $x_i$ coordinates.
Hence $\widetilde \zeta$ factors through a locally affine
map (which we give the same name)
$\widetilde \zeta: T_8^4 \to \T^2$.  Letting $m$ be an odd integer, we compute
$$\widetilde \zeta \circ \widetilde \psi(x,m)=\widetilde \zeta\bigg(x+m,x-m,\frac{x+m}{\phi},
\frac{x-m}{\phi}\bigg)=$$
$$
[(x/2\phi),x/2)]=\psi(x,\pm 1)= \psi \circ \zeta(x,m).$$
The choice of sign for $\pm 1$ depends on the congruence of $m$ mod $4$.
\endproof

Let 
\begin{equation}
\widetilde T^2(\pm)= {\rm closure\/}(
\widetilde \psi({\cal C\/}(\pm))) \subset T_8^4.
\end{equation}
Notice that the first and third (or second and fourth)
coordinates of $\widetilde \psi$
completely determine the whole image of
${\cal C\/}(+)$ in $T^4_8$.
Thus $\widetilde T^2(+)$
is a flat $2$-dimensional torus.
The same goes for $\widetilde T^2(-)$.

Let $\widetilde E_j$ be the extension to
$T_8^4$ of $E_j$.  We give the explicit
formula in \S \ref{partitioncalc}.  Referring
to Lemma \ref{zetastar} we define
\begin{equation}
\label{extension}
\widetilde \Upsilon_R: \widetilde T^2(\pm) \to \T^2; \hskip 40 pt
\widetilde \Upsilon_R= \widetilde \zeta \circ (\widetilde E_4 \circ ... \circ \widetilde E_1)^2.
\end{equation}
Combining the  Extension Lemma with Lemma
\ref{zetastar} we have
\begin{equation}
\label{semiconjugate}
\widetilde \Upsilon_R \circ \widetilde \psi=\psi \circ \Upsilon_R
\end{equation}

\begin{lemma}[Constancy Lemma]
Suppose that $X \subset \widetilde T^2(\pm )$ is a path
connected open set
on which the map $\Upsilon_R$ is entirely defined.
Suppose that $p_1,p_2 \in {\cal C\/}(\pm)$ are points
such that $\widetilde \psi(p_j) \in X$ for $j=1,2$.
Then 
$$\pi_1(\Upsilon_R(p_1)-p_1)=
\pi_1(\Upsilon_R(p_2)-p_2).$$
\end{lemma}

\startproof 
Assume without loss of generality that the path is on
$\widetilde T^2(+)$.
Since $\Upsilon_R$ is entirely defined on $X$ it is also
continous on $X$.  The map
$$v(\alpha)=\widetilde \Upsilon_R(\alpha)-\widetilde \zeta(\alpha) \in \T^2$$
is a continuous function of $\alpha \in X$.

A dense subset of $X$ 
has the form $\widetilde \psi(y,1)$, where $y \in \R$.
For such points we have
$$v(\widetilde \psi(y,1))=^1
\widetilde \Upsilon_R \circ \widetilde \psi(y,1)-
\widetilde \zeta \circ \widetilde \psi(y,1)=^2
\psi \circ \Upsilon_R(y,1)-\psi \circ \zeta(y,1)=^3$$
\begin{equation}
\label{mc}
\psi \circ \Upsilon_R(y,1)-\psi(y,1)=^4
\psi \circ \pi_1(\Upsilon_R(y,1)-y,1)
\end{equation}
The first equality is by definition; the second
one is Equation \ref{semiconjugate} and Lemma \ref{zetastar}; the third one is
the fact that $\zeta$ is the identity on
${\cal C\/}(\pm)$; the fourth one comes
from linearity, and from the fact that
$\psi$ only depends on the first coordinate.
By Lemma \ref{finiteness}, the last expression in
Equation \ref{mc} takes on only finitely
many values.  Therefore,
$v(y,1)$ takes on only finitely many values.

Now we know that the continuous map $v$ only
takes finitely many values on a dense subset of
the path connected subset $X$. Therefore
$v$ is constant on $X$.
Since $v$ is constant on $X$, the calculation in
Equation \ref{mc} gives the same answer for all
$y \in {\cal C\/}(+)$ such that
$\widetilde \psi(y) \in X$.  The conclusion of the lemma
follows immediately from setting $y=x_1,x_2$.
\endproof

\subsection{The Proof Modulo the Question of Definedness}

Referring to Equations \ref{map2} and \ref{psimap}, we have
$\Psi=\psi \circ T$, assuming that we identify the range of $T$
with either ${\cal C\/}(+)$ or ${\cal C\/}(-)$ by padding the
second coordinate with a $\pm 1$.  Given
$x_+ \in {\cal C\/}(+)$, let
$x_-$ be the point in 
${\cal C\/}(-)$ with the same first coordinate.
Let
\begin{equation}
v(x_{\pm})=\pi_1(\Upsilon_R(x_{\pm})-x_{\pm}).
\end{equation}

\begin{lemma}
\label{ag}
The point $\psi(x_+)=\psi(x_-)$ lies in an open polygon of
the partition $\cal P$ of $\T^2$, and this polygon determines
the unordered pair $\{v(x_+),v(x_-)\}$.
\end{lemma}

Computing just a single point per polygon$-$we
computed millions$-$we 
check that the  quantities given in Lemma \ref{ag},
when hit with the map $T^{-1}$, match the types shown 
in Figure 2.5.  Therefore
Lemma \ref{ag} (and this small computation) implies
the Arithmetic Graph Lemma.

The polygons of $\cal P$ do not determine the
{\it ordered\/} pair of numbers in Lemma \ref{ag}.
For this we need to pass to a finite cover.
The map
$\widetilde \zeta: T_8^4 \to \T^2$ from
Lemma \ref{zetastar} gives a finite
covering map from $\widetilde T^2(+)$ to $\T^2$.
We lift the partition $\cal P$ to a finite partition
$\widetilde P$ of $\widetilde T_2(+)$.  We
do the same thing for
$\widetilde T_2(-)$.

\begin{lemma}
\label{ag2}
The point $\widetilde \psi(x_+)$ lies in
an open polygon of $\widetilde {\cal P\/}$ and
this polygon determines $v(x_+)$. Likewise
for $v(x_-)$.
\end{lemma}

A small amount of computation shows that
Lemma \ref{ag2} implies Lemma \ref{ag}.
Hence Lemma \ref{ag2} implies the Arithmetic Graph Lemma.

\begin{lemma}
\label{defined}
Suppose that
$\widetilde \Upsilon_R$ is defined on each open
polygon of the partition $\widetilde {\cal P\/}$
of $\widetilde T^2(\pm)$.  Then Lemma \ref{ag2} is true.
\end{lemma}

\startproof
It follows immediately from Lemma \ref{defined} and
from the Constancy Lemma that
{\it if\/} $\widetilde \psi(x_+)$ lies in an open polygon of
$\widetilde {\cal P\/}$, {\it then\/} $v(x_+)$
is determined by this polygon.
Likewise for $x_-$.  We just have to rule out the
possibility that $\widetilde \psi(x_+)$ lies in the
boundary of one of the polygons. The point
$\widetilde x_-$ has the same treatment.

Suppose that $\widetilde \psi(x_+)$ lies in the
boundary of some tile of $\widetilde {\cal P\/}$.
We know $\Upsilon_R$ is defined at $x_+$.  But then
$\widetilde \Upsilon_R$ would be defined and
continuous in a neighborhood of
$\widetilde \psi(x_+)$ on $\widetilde T^2(+)$. 
But then two adjacent tiles of
$\widetilde P$ on $\widetilde T^2(+)$
would determine the same $v$ values.
But then two adjacent tiles of
$P$ would determine the same
unordered pair of $v$ values.
We check computationally
that this does not happen.  Contradiction.
\endproof

\subsection{Verifying the Definedness}
\label{binary}

At this point we have reduced the Arithmetic
Graph Lemma to verifying the hupotheses of
Lemma \ref{defined}.  We will explain how
we make the verification.

We know that $\widetilde E_{j+4}=\widetilde E_j$ for
$j=1,2,3,4$, but for notational purposes, it is
convenient to just write $\widetilde E_1,...,\widetilde E_8$.  We let
$\widetilde F_k=\widetilde E_k...\widetilde E_1$.
We interpret $\widetilde F_0$ as the identity map.
Also, $\widetilde \Upsilon_R=\widetilde \zeta \circ \widetilde F_8$.
To show that $\widetilde \Upsilon_R$ is defined
at some point $x$ it suffices to show that
the maps $\widetilde F_1,...,\widetilde F_8$ are
all defined on $x$.

We will see from the explicit description of
our map $\widetilde E_k$, given in \S \ref{partitioncalc},
that there is a union $Y_k$ of two $3$-tori, such
that $T_8^4-Y_k$ consists of two connected components
$C_k(0)$ and $C_k(1)$ on which the map
$\widetilde E_k$ is entirely defined and locally affine.

If $\widetilde \Upsilon_R$ is defined on a point
$x \in T_8^4$ then we can define a
length $8$ binary sequence
$\epsilon_1,...,\epsilon_8$, by the property
that $\widetilde F_{k-1}(x) \in C_k(\epsilon_k)$
for $k=1,...,8$.
In short, we can associate a canonical
binary sequence of length $8$ to any
point $x$ on which $\widetilde \Upsilon_R$
is defined.  We call this sequence the
{\it itinerary\/} of $x$.

Here we define a slightly more general notion
of an itinerary. 
Given an itinerary $\epsilon_1,...,\epsilon_8$ we can
define $\widetilde E_k$ on $Y_k$ by demanding that
$\widetilde E_k$ extends continuously to the closure
of $C_k(\epsilon_k)$.   Then $\widetilde E_k$ is
completely defined on $T_8^4$.  It is continuous
on the closure of $C_k(\epsilon_k)$ and on the
interior of $C_k(1-\epsilon_k)$, but not globally
continuous.  With this definition, we say that 
$x \in T_8^4$ has {\it extended itinerary\/}
$\epsilon_1,...,\epsilon_8$ if the maps
$\widetilde F_k$ (when extended) are all
defined on $x$ and $\widetilde F_{k-1}(x)$ lies
in the closure of $C_k(\epsilon_k)$ for
$k=1,...,8$.  Put another way, $x$ has
extended itinerary $\epsilon$
provided there exist points arbitrarily close
to $x$ that have itinerery $\epsilon$.

We define the {\it stretch\/} of a convex polygon
$P$ to be the
maximum distance in $T_8^4$ between consecutive
vertices of $P$.  We denote this by $\sigma(P)$.

\begin{lemma}[Definedness Criterion]
Suppose that
\begin{itemize}
\item $\widetilde \Upsilon_R$ is
defined on some point of $P$, and this
point has itinerary $\epsilon$.
\item All the vertices of $P$ have
extended itinerary $\epsilon$.
\item The stretch of $\widetilde F_k(P)$
is less than $2 \sqrt 2$ for each
$k=0,...,7$.
\end{itemize}
Then $\widetilde \Upsilon_R$ is defined on
all points of $P$, and all these points have
itinerary $\epsilon$.
\end{lemma}

\startproof
Looking at explicit formulas for our maps given in
\S \ref{cubemaps},
we see that the Euclidean distance between separare
components of $Y_k$ is at least $2 \sqrt 2$.  Hence,
if $v_1,v_2$ are two points in the closure of
$C_k(\epsilon_k)$ which are less than $2 \sqrt 2$
apart, then the line segment
$\overline{v_1v_2}$ lies in the closure 
 $\overline C_k(\epsilon_k)$.

Suppose we have shown by induction that
$\widetilde F_{k-1}$ is defined on $P$.
Let $P'=\widetilde F_{k-1}(P)$.
The hypotheses of this lemma say that
every edge of $P'$ has length
less than $2 \sqrt 2$, and the endpoints
of such an edge are in 
$\overline C_k(\epsilon_k)$.   Hence, all the
edges of $P'$ lie in
$\overline C_k(\epsilon_k)$.  Since
$\partial C_k(\epsilon_k)$ is three
dimensional, and (by induction) $P'$
is a planar polygon, we must have
$P' \subset \overline C_k(\epsilon_k)$.

If some point of the open
 $P'$ actually lies in
$\partial C_k(\epsilon_k)$ then
all of $P'$ must lie in
$\partial C_k(\epsilon_k)$. The point
here is that the tangent plane at this
bad point must be contained in the
tangent space to $Y_k$, because
there is no crossing allowed.  
We also know that some point of $P'$ lies
in $C_k(\epsilon_k)$, so the above
bad situation cannot occur.  Hence
all points of $P'$ lie in
$C_k(\epsilon_k)$. This completes
the induction step.
\endproof

We now mention a trick that makes the Definedness Criterion more
useful.    Given a line segment $s \in T_8^4$, with
endpoints $p_1$ and $p_2$,  we let
$s'$ denote the partition of $s$ into the two segments
$[p_1,q]$ and $[q,p_2]$ where
\begin{equation}
\label{subdivide}
q=p_1 \phi^{-2} + p_2 \phi^{-1}.
\end{equation}
Note that $q \in s$ because $\phi^{-1} + \phi^{-2}=1$.
We might have chosen $q$ to be the midpoint of $s$, but our
choice interacts better with $\Z[\phi]$.   Given a polygon
$P$, we let $P'$ denote the polygon, with twice as many
vertices, obtained by subdividing each edge of $P$.
In general, let $P^{(n)}=(P^{(n-1)})'$.   Then
$P^{(n)}$ has $2^n$ times as many vertices as
$P$ and the stretch of $P^{(n)}$ is 
$\phi^{-n}$ times the stretch of $P$. Thus, by
taking a sufficiently large integer $n$, we can
guarantee the last condition of the Discreteness
Criterion without even computing the stretch.

We take each polygon $P$ of our partition
$\widetilde {\cal P\/}$ and perform the
following calculation.  We take the
$10$th subdivision $P^{(10)}$ and check
the Definedness Criterion. (This is overkill.)
This verifies the hypotheses of
Lemma \ref{defined} and thereby completes the proof
of the Arithmetic Graph Lemma.  

We give details of our
calculation in \S \ref{partitioncalc}.

\newpage

%% file: 7routines.tex
\section{The Code in Detail}

\subsection{Arithmetic Operations}

We do all our computations with complex numbers of the form
\begin{equation}
\frac{x_0+x_1 \phi}{2}+I\ \frac{x_2+ x_3 \phi}{2}.
\end{equation}
Here $I=\sqrt{-1}$.
We call such a number an {\tt IntegerComplex\/} and represent it as a sequence
\begin{equation}
\{x_0,x_1,x_2,x_3\}.
\end{equation}
Given $a$ and $b$, both {\tt IntegerComplex\/} objects, we perform the
following operations:
\begin{itemize}
\item $a+b$ and $a-b$ are computed by adding or subtracting components.
\item $2ab$, another {\tt IntegerComplex\/}, is computed by expanding out all the terms and grouping them.
\item The conjugate $\overline a$ is obtained by negating the third and fourth
coordinates of $a$.
\item We have a routine {\tt interpolate\/}($a$,$b$), which computes
$a\phi^{-1}+b \phi^{-2}$.
\end{itemize}

Now we explain how we test the sign of the number $a=a_0+a_1\phi$.
Since $\phi$ is irrational, we need a trick.
The idea is to consider
\begin{equation}
f_{50}=12586269025; \hskip 20 pt
f_{51}=20365011074 \hskip 20 pt
f_{52}=32951280099.
\end{equation}
Here $f_n$ is the $n$th Fibonacci number.
The quantity $a_0+a_1\phi$ is positive (respectively negative) provided
that both sums
\begin{equation}
s_1(a)=a_0f_{50}+a_1 f_{51}; \hskip 50 pt
s_2(a)=a_1f_{51}+a_2 f_{52}
\end{equation}
are positive (respectively negative).  This works because
the sign of the difference $\phi-f_{n+1}/f_n$ alternates with $n$.
We implement our routine, called {\tt sign\/}, using the {\tt BigInteger\/} class in
Java, which does integer arithmetic correctly on huge integers.
The routine {\tt sign\/} certainly could fail for some inputs, but it
never fails for the inputs we give it.  The point is that we always give
it inputs involving pretty small integers.

Our routine {\tt dec\/} takes the decimal part of a {\tt IntegerComplex\/} and
returns a {\tt IntegerComplex\/} in $[-1/2,1/2]$.  In defining {\tt dec\/}
we take the decimal part of a {\tt IntegerComplex\/} using floating point
arithmetic to compute the nearest integer. 
The computation does not use exact integer arithmetic, but then we
use the {\tt sign\/} routine to check rigorously that our guess always lies in
$[-1/2,1/2]$, thereby guaranteeing that get the same answer {\it as if\/} we
had used purely integer arithmetic calculations.

With these preliminaries in place, we describe the
{\tt IntegerComplex\/} version of our map $\Psi$
given in Equation \ref{map2}.
\newline
\newline
\noindent
{\tt psi\/}($a$,$b$): \newline
let $x=4-12a+4b$ and  $y=8a-2$.  (the {\tt IntegerComplex\/} version of Eq \ref{map}) \newline
Let $c_1=\{(-x+y)/2,x/2\}$.  (division by $2\phi$; works because $x,y$ are even) \newline
Let $c_2=\{x/2,y/2\}$. (division by $2$; works because $x,y$ are even) \newline
Let $d_j=${\tt dec\/}($c_j$) for $j=1,2$. \newline
return($d_1+Id_2)$.
\newline

\subsection{Classification into Types}
\label{types}

An {\tt IntegerPolygon\/} is a finite list of {\tt IntegerComplex\/}es,
namely the vertices.
Given an {\tt IntegerComplex\/} $z$ and an {\tt IntegerPolygon\/} $P$,
let $T_i(z,P)$ denote the triangle determined by the ordered triple of
{\tt IntegerComplex\/}es $z$, $P(i)$ and $P(i+1)$.  The indices
are taken cyclically.  Using our {\tt sign\/} routine, we have a
straightforward routine {\tt signArea\/} that computes the orientation
of a triangle of {\tt IntegerComplex\/}es and vanishes of the
points are collinear.
The following routine returns a $1$ if $z$ is
contained in the interior of $P$ and a $0$ otherwise.  A straightforward
variant, {\tt IsContainedClosed\/} checks if
$z \in \overline P$.
\newline
\newline
{\tt isContainedOpen\/}($z$,$P$): \newline
loop for the number of vertices of $P$: \newline
\indent check that $T_i(z,P)$ and $T_{i+1}(z,P)$ have nonzero {\tt signArea\/} \newline
\indent if(false) return(0) \newline
\indent check that $T_i(z,P)$ and $T_{i+1}(z,P)$ have the same {\tt signArea\/} \newline
\indent if(false) return(0) \newline
endloop \newline
return(1)
\newline

Recall that the arithmetic graph $\Gamma$ is a certain 
subset of $\Z^2 \cap H$, where $H$ is a certain half-plane.
Each point $(x,y) \in \Z^2$ has a type, which we compute
by determining which open polygon of $\cal P$ contains
$\Psi(x,y) \in \T^2$.   Here we explain how we compute this
in practice.  Lifting $\cal P$ to $\R^2$ we get a
tiling $\widetilde{\cal P}$ of $\R^2$.  We treat
{\tt psi\/}($x$,$y$) as a point in $\R^2$ 
and check which tile of $\widetilde{\cal P}$ it
lands in.  We have arranged that the image of {\tt psi\/}
is fairly close to the origin, and so we only have to
check a smallish portion of the tiling.   

We accomplish
our goal using two routines.
The first of our routines checks whether or not an {\tt IntegerComplex\/} $z$
near the origin in $\R^2$ is contained in a union of integral translates
of a given {\tt IntegerPolygon\/} $P$.   The program returns a
$1$ if the answer is yes. 
\newline
\newline
{\tt isLatticeContainedOpen\/}($z$,$P$): \newline
loop over $i$ from $-3$ to $3$ and
over $j$ from $-3$ to $3$ \newline
\indent if({\tt isContainedOpen\/}($z+(i,j)$,$P$)$=1$ then return(1) \newline
endloop
\newline

We always take $P$ as one of the polygons from our
partition $\cal P$.  These polygons are listed in \S \ref{polylist}.
(See also Figure 2.5.)
Our routine {\tt classify\/} combines our routine {\tt psi\/} with the routine
{\tt isLatticeContainedOpen\/} to classify each point of $\Z^2$ into the
local types.

We use the routine {\tt isContainedClosed\/} in place of
 {\tt isContainedOpen\/} in case we want to
check that a given point is contained in a closed polygon.

\subsection{Types and Dynamical Translations}

There are $23$ nontrivial local types of vertex in
the arithmetic graph, as shown in Figure 2.4.  We list
these types here (rather than in the appendix.)  The array

\begin{equation}
\label{type}
\matrix{T_j \cr a_1&b_1 \cr a_2&b_2}
\end{equation}
indicates that one of the edges emanating from a vertex of type $T_j$.
is $(a_1,a_2)$ and the other one is $(b_1,b_2)$.  In other words,
the two columns of the matrix encode the type.  The ordering
of the two columns is arbitrary.  We list the types
$T_1,T_3,...,T_{23}$.  Similar to what happens above,
$T_{2j}$ is obtained by negating all the entries of $T_{2j-1}$,
for $j=1,...,11$.  Here is the list:

$$\matrix{T_1 \cr 1&0 \cr 1&1} \hskip 40 pt
\matrix{T_3 \cr -1 &1 \cr 1&1} \hskip 40 pt
\matrix{T_5 \cr -1&0 \cr 1 & -1}\hskip 40 pt
\matrix{T_7 \cr -1&0 \cr 1 & -1}$$

$$\matrix{T_9 \cr 0 & -1 \cr 1 & -1} \hskip 40 pt
\matrix{T_{11} \cr 1 & -1 \cr 1 & 0} \hskip 40 pt
\matrix{T_{13} \cr 1 & -1 \cr 1 &0 }\hskip 40 pt
\matrix{T_{15} \cr 0 & -1 \cr -1 & 0}$$

$$\matrix{T_{17} \cr 0 & -1 \cr -1 & 0} \hskip 40 pt
\matrix{T_{19} \cr 0 & -1 \cr 1 & 0} \hskip 40 pt
\matrix{T_{21} \cr 0 & 0 \cr 1 & -1 }\hskip 40 pt
\matrix{T_{23} \cr 0 & 0 \cr 1 & -1}$$

For each of these types, there are two associated dynamical translations of $\T^2$.
Referring to the matrix in Equation \ref{type}, the two maps are:
\begin{equation}
   (x,y) \to [(a_1 \phi^{-4}+a_2 \phi^{-1},a_1\phi^{-3})]; \hskip 30 pt
     (x,y) \to [(b_1 \phi^{-4}+b_2 \phi^{-1},b_1\phi^{-3})].
\end{equation}
These maps relate the points $\Psi(v)$ and $\Psi(v')$ where $v$ is the
vertex and $v'$ is one of the vertices connected to $v$ by the
arithmetic graph. 
The following routine starts with a pair $(x,y)$ obtained
from one of the columns of the above matrices and returns
the {\tt IntegerComplex\/} $z$ which effects the corresponding
dynamical translation. That is, on $\R^2$ the dynamical
translation is given by $w \to w+z$.  The next routine gets the
{\tt IntegerComplex\/} by which we translate.  Here
$x,y \in \{-1,0,1\}$.
\newline
\newline
{\tt getMap\/}($x$,$y$): \newline
return the integer complex with coordinates $(10x-2y,-6x+2y,-6x,4x)$

\subsection{The Sequence Generator}

If $v_1$ and $v_2$ are two consecutive vertices of the arithmetic graph,
then one of the two maps associated to $v_1$ coincides with one of the
two maps associated to $v_2$.  For instance, if a vertex of type $1$ is
connected to a vertex of type $21$ we could write
$$\matrix{T_1&T_{21} \cr 1&0&0 \cr 1 & 1 & -1} \hskip 30 pt
{\rm or\/} \hskip 30 pt
\matrix{1&21 \cr 1&0 \cr 1&1}$$
The second notation system is a simplification of the first.
We drop off the last column because it is not of interest to us.
The information in the right hand side 
array is enough to generate both $v_1$ and $v_2$ and
also to determine the type of $v_1$.  

We can encode longer sequences of types using longer arrays, as
we now illustrate by example:
Our first gene is located at the point $p=(3,4)$.
The two arrays are
$$A((3,4),1,3)=\matrix{11&9&23 \cr 1&0&0 \cr 1&1&1}; \hskip 30 pt
A((3,4),2,3)=\matrix{11&14&10 \cr -1&-1&0 \cr 0 &-1 &-1}$$
This tells us that the $5$ vertex types we see along the
gene $A_0$ are (in one of the two orders)
$10,14,11,9,23.$
Looking at the first column
of $A((3,4),1,3)$, we can see that
that the vector $(1,1)$ connects the point of type $11$ to the
point of type $9$.   In this way, we can draw a copy of the gene
given the above arrays.
The next routine generates the arrays for the point $p=(x,y) \in \Z^2$.
\newline
\newline
{\tt getOrbit\/}($x$,$y$,length,epsilon): \newline
1. let $A$ be the empty array \newline
2. let $X=x$ and $Y=y$.
3. let count$=0$ \newline
4. while(count$<$length): \newline
\indent $a=${\tt classify\/}($X$,$Y$)  \newline
\indent Let $M$ be the matrix associated to $T_a$ \newline
\indent If(count$=0$) then: \newline
\indent \indent append to $A$ the (epsilon)th column of $M$ \newline
\indent if(count$>0$) then: \newline
\indent \indent append to $A$ the col. of $M$ which does not match the 
last col. of $A$.\newline
\indent Let $m_x$ and $m_y$ denote the entries of the column of $M$ used above. \newline
\indent replace $X$ by $X+m_x$\newline
\indent replace $Y$ by $Y+m_y$\newline
\indent increment count \newline
5. return($A$)
\newline

{\tt getOrbit\/} runs until the count equals the length, and then breaks.
At that point, the array $A$ is returned.  For our purposes, we need to
get both sequences associated to $p$ and $k$. So, we run the above
routine once for $\epsilon=1$ and once for $\epsilon=2$.   We
call the resulting pair $(L,R)$ of arrays the $k$-{\it itinerary\/} of 
$p \in \Z^2$.   The genes correspond to the case $k=3$.

\subsection{Checking the Dynamical Polygons}
\label{polylist}
Each array $A$ determines a
dynamical polygon $P=P(A)$, as in \S 4.2.  Here we explain
how to check directly that a point $s \in \T^2$ lies
in the closure of $P$ or else in the interior of $P$
(depending on our interest.) 
The point $s$ is always given as the projection
to $\T^2$ of a certain lift $s \in \R^2$ which we give
the same name.  Let $\P_j$ denote the $j$th tile in
the partition $\cal P$ of $\T^2$. As above, we identify
$\P_j$ with a particular lift to $\R^2$.
The following routine returns a $1$ provided that
$s$ is contained in the interior of $P$.
\newline
\newline
{\tt matchOrbitOpen\/}($s$,$A$): \newline
let $S=s$ \newline
loop for $i=1$ to length of $A$ \newline
\indent check that {\tt isLatticeContainedOpen\/}($S$,$\P_{a_i}$)$=1$. \newline
\indent if false return($0$) \newline
\indent if true then: \newline
\indent \indent let $z=${\tt getMap\/}($x_i$,$y_i$) \newline
\indent \indent let $S=S+z$  (act on $S$ by the dynamical translation) \newline
endloop
\newline

As a variant, we check that $s$ is contained in the closure of 
$P$ by using 
{\tt isLatticeContainedClosed\/} in place of
{\tt isLatticeContainedOpen\/}.   

We can now explain exactly how we check our
guess $P_G$ for the $75$ dynamical polygons associated to the
$75$ genes. Here is the $3$-step process.
\begin{itemize}
\item For the $j$th gene $A_j$ we let $(x_j,y_j)$ be the location
of the center $a_j$ of $A_j$.  We use 
{\tt getOrbit\/}($x_j$,$y_j$) to
generate the $3$-itinerary $A=(L_j,R_j)$.  
\item We check that each {\tt IntegerComplex\/} $v$ representing 
a vertex of $P_G$ is contained in the closure of $P=P_A$
using  {\tt matchOrbitClosed\/}($v$,$L_j$) and then
{\tt matchOrbitClosed\/}($v$,$R_j$).
\item For each edge $e$ of $P_G$ we produce a {\tt IntegerComplex\/} $v$
contained in the interior of $e$ by applying {\tt interpolate\/} to
the two endpoints of $e$.  Given $v$ we then check that one of the
two runs  {\tt matchOrbitOpen\/}($v$,$L_j$) and
{\tt matchOrbitOpen\/}($v$,$R_j$) returns a $0$. This is to say that
$v$ is not contained in the interior of $P$.
\end{itemize}

We list the dynamical polygons $P_0,...,P_{74}$ in
\S \ref{polylist2}.  
We list the centers for our genes in \S \ref{centerlist}.

\subsection{The Shadowing Property}

We begin by recalling the basic setup.
Let $A$ be a gene from our inflation generator and let $B$ be the gene core.
Let $a$ be the center vertex of $A$.  Let $A'$ be the strand which
shadows $\Phi(A)$.  Recall that $\Phi$ is dilation by $\phi^3$.
Let $P$ be the dynamical polygon associated to $A$ and let $P'$ be the
dynamical polygon associated to $A'$.  Let $\beta$ be the inflation map
and let $\gamma: P \to \T^2$ be the associated special similarity.  We
want to verify that $\gamma(P) \subset P'$.  

The first main task is to compute the dymamical sequences associated
to $A'$
We store $A'$ as a triple $\{(x_j,y_j)\}$ of points.
Here $(x_2,y_2)$ is the point $a'$ that comes close to $\Phi(a)$, and
$(x_1,y_1)$ and $(x_3,y_3)$ are the two endpoints of $A'$.
In \S \ref{shadowlist} we give the list of triples.

A variant of {\tt getOrbit\/} recovers the sequence of types from the triple.
This variant starts at $a'$ and moves out in either direction
until it encounters the two endpoints.  In other words, we just
replace line $4$ of {\tt getOrbit\/} with a check that the current
point equals neither $(x_1,y_1)$ nor $(x_3,y_3)$.
Technically we could avoid this reconstruction problem just by saving the sequences
attached to $A_1',...,A_{75}'$, but we prefer not to store so much data.
Below we will list out the $75$ triples.

Now we explain how to compute $\gamma$.
We have already defined the routine {\tt psi\/} above, which chooses some
lift of $\Psi(a)$ to $\R^2$.   Here we explain an improved version, which
chooses the lift of $\Psi(a)$ that is contained in $P$.  The idea is to
find a ``correction vector'' $\lambda \in \{-7,7\}^2$ such that
{\tt psi\/}($a$)+$\lambda \in P$.   This always works, and we call the
result {\tt canonicalPsi\/}.  It is a better lift of $\Psi$.
From the way we have
constructed $A'$, we know that $\gamma$ maps {\tt canonicalPsi\/}($a$) into
$P'$.  Thus, we can use our knowledge of {\tt canonicalPsi\/}($a$) to
compute $\gamma$.    We seek a pair of integers $(m,n)$ such that
\begin{equation}
\label{gam}
\gamma=\gamma_0+[(m \phi^{-4}+n \phi^{-1},0)]; \hskip 15 pt
\gamma_0(x,y)=(-\phi^{-3}x,-\phi^{-3}y).
\end{equation}
We set $y=${\tt canonicalPsi\/}($a$) and then test small integer choices of
$m$ and $n$ until we find a choice which leads to $\gamma(y) \in P'$.
To test that $\gamma(y) \in P'$ we don't need to compute $P'$.  We
simply check that $\gamma(y)$ follows the dynamical sequences associated
to $A'$.   In other words, we just go back to the definition of $P'$.
We make the same kind of loop as in the routine {\tt canonicalPsi\/},
with $|m|,|n| \leq 7$.  This always works.

Once we have $\gamma$ we list out the vertices of $P$ as
$P_1,...,P_k$ and check that $\gamma(P_i)$ follows the dynanical
sequences associated to $P'$.  Here $k \leq 5$ for each of the $75$ choices.

\subsection{Checking Coherence}

Let $g_1,...,g_{1024}$ be the first $1024$ points along $\Gamma_0$. 
The point $g_j$ is the center of a gene $G_j$.  
Let $\beta_j$
be the inflation map we associate to $G_j$ using our
inflation structure $\chi$.
The list $\beta_1,...,\beta_{1024}$
consists of multiple references to the same
$75$ inflation maps.
For each index $j$ we
produce a point $g'_j \in \Gamma_0$, always
within $3$ units of $\Phi(g_j)$, such that
\begin{equation}
\label{cano}
\Psi(g'_j)=\gamma \circ \Psi(g_j).
\end{equation}
This is to say that $\beta_j(g_j)=g'_j$.
We check Equation \ref{cano} as follows:
\newline
\newline
{\tt verifyCoherence:\/} \newline
loop from $j=1$ to $1024$. \newline
\indent 1. Let $k=I_j$ (the gene type of $G_j$) \newline
\indent 2. let $y_j=${\tt canonicalPsi\/}($g_j$) \newline
\indent 3. let $\gamma$ be the map computed in Equation \ref{gam} for the index $k$. \newline
\indent 4. Check that $\gamma(y_j)$ and {\tt psi\/}($g_j$) agree up to a vector in $\{-7,7\}^2$.
\newline
endloop \newline \newline
Line $4$ really says that $\gamma(y_j)=\Psi(g_j)$ in $\T^2$.

As a final detail,
we check visually that $g_1',...,g'_{1024}$ lie in
order on $\Gamma_0$.  This verifies the coherence
of our inflation generator.

In the interest of space, we don't list the pairs
$(g_j,g'_j)$ in the appendix.  The interested reader can
push the button {\bf shadow points\/} on the {\bf Arithmetic Graph Control Panel\/}
of Billiard King and see these points plotted.
The first few pairs $(g_j,g_j')$ are
$$ ((1,1),(3,4));
\hskip 15 pt ((1,2),(3,9)) \hskip 15 pt
((1,3),(3,13)) \hskip 15 pt
((2,4),(8,18)).$$
For reference, the following routine generates $g'_j$ given $g_j$.
\newline
\newline
{\tt generateShadow\/}($j$): \newline
1. Let $k=I_j \in \{1,...,75\}$ be the type of the gene $G_j$ \newline
2. Let $a_k$ be the center of the $k$th gene in our inflation generator \newline
3. Let $a'_k$ be the point of $A'_k$ that shadows $\Phi(a_k)$ \newline
4. Search for a point $g'_j \in \Z^2$ within $3$ units of $\Phi(g_j)$ such that
$g'_j$ and $a'_k$ are the centers of equivalent genes.
\newline

\subsection{Checking the Partition}
\label{partitioncalc}

Here we explain the calculations from the end of \S \ref{binary}.

\subsubsection{Explicit Formulas}
\label{cubemaps}

Let $\widetilde E_j$ be the extension of
$E_j$ to $T_8^4$ guaranteed by the
Extension Lemma.  We can
rigorously compute the formulas for these
maps just by making a few calculations.
(We checked the formulas on millions of points.)

One general feature of the map $\widetilde E_j$ is
that the corresponding ``undefined set'' $Y_j$
consists of two parallel $3$-tori.  The
complementary region $T_8^4-Y_j$
consists of two regions.
This feature is not necessarily
clear from our description below, because we describe
our maps in terms of
coordinates on the set $(-4,4)^4$, which is the interior
of a fundamental domain for $T^4_8$.  The action on $T^4$
is obtained by piecing together the definitions across
the boundaries of $[-4,4]^4$.
In all cases, we will list the following data:
\begin{itemize}
\item The linear part $L$ of $\widetilde E_j$.
\item A {\it determiner function\/} $d:(-4,4)^4 \to \R$.
\item A partition of $\R$ either into $3$ or $5$ intervals.
In the $3$-interval case, the dividing points are
$0$ and $4$ and the intervals are
$(-\infty,0)$ and $(0,4)$ and $(4,\infty)$.
In the $5$-interval case, the dividing points
are $-8,-4,0,4$ and the intervals are determined similarly.
\item For each interval $I$ of the partition we give
an array of the form
$$v(I)=\left[\matrix{
a_{11}&a_{21}&a_{31}&a_{41} \cr
a_{12}&a_{22}&a_{32}&a_{42}} \right]
$$
This array stands for the vector
$$\bigg(
\frac{a_{11}+a_{12}\phi}{2},
\frac{a_{21}+a_{22}\phi}{2},
\frac{a_{31}+a_{32}\phi}{2},
\frac{a_{41}+a_{42}\phi}{2}\bigg)
$$
\end{itemize}
In case all the entries of the array are $0$ we will
save space by just writing $[0]$ in place of the array.
Our map is then given by
$\widetilde E(x)=L(x)+v(I)$, where $I$ is such that
$d(x) \in I$.  This definition makes sense except for the
points $x$ where $d(x)$ lies in one of the dividing points
listed above.

$$
\widetilde E_1: \hskip 20 pt
L_1(x)=\left(\matrix{x_1 \cr x_2 \cr x_2/\phi+x_3+x_4 \cr x_2/\phi}\right);
\hskip 20 pt
d_1(x)=x_1+1;$$
\begin{equation}
\left[\matrix{1&1&1&1\cr0&0&-1&1}\right];
\hskip 8 pt [0]
\hskip 8 pt
\left[\matrix{1&1&1&1\cr 0&0&1&-1}\right].
\end{equation}
\newline
\newline
$$ \widetilde E_2: \hskip 20 pt
L_2(x)=\left(\matrix{x_1 \cr x_2 \cr x_3 \cr -x_1+x_2/\phi+x_3 \phi}\right); \hskip 15 pt
d_2(x)=x_2+x_3-\frac{1}{\phi^2};$$
\begin{equation}
\left[\matrix{0&0&0&0\cr 0&0&0&2}\right] \hskip 4 pt
\left[\matrix{0&1&0&1\cr 0&0&0&1}\right] \hskip 4 pt
[0] \hskip 4 pt
\left[\matrix{0&1&0&1\cr 0&0&0&-1}\right] \hskip 4 pt
\left[\matrix{0&0&0&0\cr 0&0&0&-2}\right]
\end{equation}
\newline
\newline
$$\widetilde E_3: \hskip 20 pt
L_3(x)=\left(\matrix{
x_1 \phi - x_2-x_3 +x_4\phi \cr
x_1/\phi - x_3 + x_4 \phi \cr
-x_1/\phi+x_2+2x_3 - x_4\phi \cr
-x_1/\phi + x_2 + x_3 - x_4/\phi}\right); \hskip 30 pt
d_3(x)=x_1+x_4-\frac{1}{\phi^2};$$
{\small
\begin{equation}
\left[\matrix{0&0&0&0\cr 2&2&-2&-2}\right] \hskip 4 pt
\left[\matrix{1&1&0&0\cr 1&1&-1&-1}\right] \hskip 4 pt
[0]    \hskip 4 pt
\left[\matrix{1&1&0&0\cr -1&-1&1&1}\right] \hskip 4 pt
\left[\matrix{0&0&0&0\cr -2&-2&2&2}\right]
\end{equation}
\/}
\newline
\newline
$$\widetilde E_4: \hskip 20 pt
L_4(x)=\left(\matrix{x_1\cr x_2 \cr x_1/\phi \cr x_4}\right); \hskip 30 pt
d_4(x)=x_0+1$$
\begin{equation}
\left[\matrix{1&0&1&0\cr0&0&1&0}\right]
\hskip 8 pt
[0]
\hskip 8 pt
\left[\matrix{1&0&1&0\cr 0&0&-1&0}\right].
\end{equation}

\subsubsection{Reducing to Planar Polygons}

The polygons in the partition $\widetilde {\cal P\/}$ are subsets
of the torus $T_8^4$.  In this section we explain how we arrange
our computation so that we just have to consider polygons which
are subsets of $\R^4$.

Let
$$\T^2_4=\R^2/(4\Z^2).$$
Then $\T^2_4$ is a square torus which is naturally
a $16$-fold cover of $\T^2$.  Let
$\widehat {\cal P\/}$ denote the lift of
$\cal P$ to $T_2^4$.  If $P_k$ is a polygon
of $\cal P$, then
$P_k$ is the projection
to $\T^2$ of the convex hull of certain vertices
$\{v_{k1},...,v_{kn}\}$. 
For each 
$(\epsilon_1,\epsilon_2) \in \{0,1,2,3\}^2$
we can form the convex polygon
$P_k(\epsilon_1,\epsilon_2)$ with vertex list
\begin{equation}
\label{liftlist}
\{v_{k1}+(\epsilon_1,\epsilon_2),...,
v_{kn}+(\epsilon_1,\epsilon_2)\}
\end{equation}
and project it into $\T^2_4$.  Call the resulting
polygon $\widehat P_{k}(\epsilon_1,\epsilon_2)$.
The union of the polygons
$\widehat P_k(\epsilon_1,\epsilon_2)$,
taken over all indices, gives us the
partition $\widehat {\cal P\/}$.
So, just to be clear, $P_k(\epsilon_1,\epsilon_2)$
is some concrete lift to $R^2$ of a polygon
of the partition $\widehat P$.

We are interested in the polygons of the
partition $\widetilde {\cal P\/}$, which is a subset of
the different torus
$T_8^4$.  We now explain how to translate between
$\widehat {\cal P\/}$ and $\widetilde {\cal P\/}$.
The map $\delta_8(x)=8\delta(x/8)$ computes $x$ mod $8\Z$.
This map is the building block for our map $\widetilde \psi$.
We introduce the maps $\mu_-$ and $\mu_+$, defined by
$$\mu_-(x,y)=((\delta_8(2x+1),\delta_8(2x-1),\delta_8(2y+\phi^{-1}),\delta_8(2y-\phi^{-1}))$$
$$\mu_+(x,y)=((\delta_8(2x-1),\delta_8(2x+1),\delta_8(2y-\phi^{-1}),\delta_8(2y+\phi^{-1}))$$
We also introduce the map
\begin{equation}
\widehat \psi(x)=
(\delta_4(x/(2\phi)),\delta_4(x/2))=(4\delta(x/(8\phi)),4\delta(x/8))
\in T_4^2.
\end{equation}
We compute that
$$\mu_+ \circ \widehat \psi(x,1)=\widetilde \psi(x,1); \hskip 30 pt
\mu_- \circ \widehat \psi(x,-1)=\widetilde \psi(x,-1).$$
It follows from this last equation that
$\mu_+$ maps the polygons of $\widehat {\cal P\/}$ to the polygons
of $\widetilde {\cal P\/}$ which partition $\widetilde T^2(+)$.
A similar statement holds for $\mu_-$.

We can think of $\mu_+$ and $\mu_-$ as being defined on
$\R^2$, and thus the polygons of $\widetilde {\cal P\/}$ all
have the form 
$$\mu_{\pm}(P_k(\epsilon_1,\epsilon_2))$$
for $k=1,...,26$ and $\epsilon_j \in \{0,1,2,3\}$.

\subsubsection{The Calculation}

Here's the routine that produces a lift
of the $a_2$nd point of $\mu_{\pm}(Q)$, where $Q=Q_{a_1}(a_3,a_4)$.  The
integer $a_5 \in \{0,1\}$ toggles $\mu_+$ and $\mu_-$.  
\newline
\newline
{\tt getVertexPlanar\/}($a_1$,$a_2$,$a_3$,$a_4$,$a_5$): \newline
Let $v=x_0+iy_0$ be the $a_2$ vertex of our particular lift of $P_{a_1}$. \newline
Let $x_2=x_1+a_3$ and $y_2=y_1+a_4$. (Equation \ref{liftlist}.) \newline
If $a_5=0$ then return $(2x_1+1,2x_2-1,2y_2+1/\phi,2y_2- 1/\phi)$. (apply $\mu_+$.) \newline
If $a_5=1$ then return $(2x_1-1,2x_2+1,2y_2-1/\phi,2y_2+ 1/\phi)$. (apply $\mu_-$.)
\newline

Next, we have a routine {\tt dec8\/}, which computes $x$ mod $8\Z$, where
$x \in \frac{1}{2}\Z$.   Here {\tt dec8\/} works just like our
routine {\tt dec\/} above, and is checked in the same way.  In
brief {\tt dec8\/}($x$)=$8${\tt dec\/}($x/8$).

Our routine {\tt subdivide\/} starts with a polygon $P \subset \R^4$
and adds one new point {\tt interpolate\/}($v_i$,$v_{i+1})$ between
consecutive vertices $v_i$ and $v_{i+1}$ of $P$.  This, the
new polygon has twice as many vertices as the old one.
Our routine {\tt getPolygon\/} subdivides a polygon
in $\R^4$ ten times and then projects it into
$T_8^4$ using {\tt dec8\/} componentwise on all the vertices.

We use a routine {\tt getTracePoint\/} to get an interior
point of $Q$.  This routine returns the point
${\tt interpolate\/}(x_3,{\tt interpolate\/}(x_1,x_2))$
where $x_1,x_2,x_3$ are the first three vertices of $Q$.

For each of the relevant indices $(a_1,a_3,a_4,a_5)$ we let
$P$ be the polygon produced by {\tt getPolygon\/} and we let
$x \in P$ be the point produced by {\tt getTracePoint\/}.
Our routine {\tt getItinerary\/} computes the itinerary of
$x$ in a straightforward way, just by computing the orbit and
checking which component $C_k(\epsilon_k)$ contains the
relevant point $x_k$ of the orbit for $k=1,...,8$.
The itinerary is then $\epsilon_1,...,\epsilon_8$, as
in \S \ref{binary}.
Finally, for each vertex $y$ of $P$ we perform the following routine
\newline
\newline
{\tt verifyItinerary\/}($\epsilon$,$y$): \newline
Let $y_1=y$.
loop from $k=1$ to $8$. \newline
\indent Verify that $y_k \in C_k(\epsilon_k)$. \newline
\indent If false then return({\tt false\/}).  Otherwise continue. \newline
\indent Let $y_{k+1}= \widetilde E_k(y_k)$ using the $\epsilon_k$ extension of
$\widetilde E_k$. \newline
return({\tt true\/})
\newline

\newpage

%% file: 8appendix.tex
\section{Appendix}

\subsection{The Polygons in the Partition}
\label{polylist1}

Here we list out the coordinates of the polygons in our partition
$\cal P$ of $\T^2$.  These polygons are drawn
in Figure 2.5.
Each vertex of $P_j$ has
the form
$$\left(\frac{a_0+a_1 \phi}{2},\frac{a_2 + a_3 \phi}{2}\right).$$
We will simply list such a vertex as a row 
$$\matrix{a_0&a_1&a_2&a_3}$$
in an array whose other rows correspond to
the other vertices of $P_j$. 

{\tiny
$$
\matrix{P{1} \cr
5&-3&0&0 \cr
10&-6&-6&4 \cr
2&-1&0&0
} \hskip 15 pt
\matrix{P{2} \cr
-5&3&0&0 \cr
-10&6&6&-4 \cr
-2&1&0&0
} \hskip 15 pt
\matrix{P{3} \cr
0&0&0&0 \cr
5&-3&0&0 \cr
10&-6&-6&4 \cr
5&-3&-6&4
}$$

$$\matrix{P{4} \cr
0&0&0&0 \cr
-5&3&0&0 \cr
-10&6&6&-4 \cr
-5&3&6&-4
} \hskip 15 pt
\matrix{P{5} \cr
0&0&0&0 \cr
13&-8&-6&4 \cr
5&-3&-6&4
} \hskip 15 pt
\matrix{P{6} \cr
0&0&0&0 \cr
-13&8&6&-4 \cr
-5&3&6&-4
} \hskip 15 pt
\matrix{P{7} \cr
9&-5&-4&4 \cr
4&-2&2&0 \cr
1&0&2&0 \cr
14&-8&-4&4
}$$

$$\matrix{P{8} \cr
-9&5&4&-4 \cr
-4&2&-2&0 \cr
-1&0&-2&0 \cr
-14&8&4&-4
} \hskip 15 pt
\matrix{P{9} \cr
3&-2&-4&2 \cr
-5&3&2&-2 \cr
3&-2&-2&0 \cr
-2&1&-2&0 \cr
6&-4&-8&4
} \hskip 15 pt
\matrix{P{10} \cr
-3&2&4&-2 \cr
5&-3&-2&2 \cr
-3&2&2&0 \cr
2&-1&2&0 \cr
-6&4&8&-4
} \hskip 15 pt
\matrix{P{11} \cr
0&0&0&0 \cr
5&-3&-4&2 \cr
-3&2&0&0
}$$

$$\matrix{P{12} \cr
0&0&0&0 \cr
-5&3&4&-2 \cr
3&-2&0&0
} \hskip 15 pt
\matrix{P{13} \cr
10&-6&-6&4 \cr
5&-3&-6&4 \cr
-3&2&4&-2
} \hskip 15 pt
\matrix{P{14} \cr
-10&6&6&-4 \cr
-5&3&6&-4 \cr
3&-2&-4&2
} \hskip 15 pt
\matrix{P{15} \cr
13&-8&-6&4 \cr
5&-3&-6&4 \cr
-3&2&4&-2 \cr
5&-3&-2&2
}$$

$$\matrix{P{16} \cr
-13&8&6&-4 \cr
-5&3&6&-4 \cr
3&-2&-4&2 \cr
-5&3&2&-2
} \hskip 15 pt
\matrix{P{17} \cr
12&-8&-8&4 \cr
7&-5&-8&4 \cr
4&-3&-4&2
} \hskip 15 pt
\matrix{P{18} \cr
-12&8&8&-4 \cr
-7&5&8&-4 \cr
-4&3&4&-2
} \hskip 15 pt
\matrix{P{19} \cr
-4&2&-2&0 \cr
1&-1&-2&0 \cr
4&-3&-6&2
}$$

$$\matrix{P{20} \cr
4&-2&2&0 \cr
-1&1&2&0 \cr
-4&3&6&-2
} \hskip 15 pt
\matrix{P{21} \cr
1&-1&0&0 \cr
6&-4&-6&4 \cr
-2&1&0&0
} \hskip 15 pt
\matrix{P{22} \cr
-1&1&0&0 \cr
-6&4&6&-4 \cr
2&-1&0&0
} \hskip 15 pt
\matrix{P{23} \cr
6&-3&-2&2 \cr
-2&2&2&0 \cr
-4&3&4&-2 \cr
4&-2&0&0
}$$

$$\matrix{P{24} \cr
1&-1&0&0 \cr
-2&1&0&0 \cr
3&-2&-4&2 \cr
1&-1&-2&0 \cr
-1&0&-2&0 \cr
4&-3&-4&2
} \hskip 15 pt
\matrix{P{25} \cr
-1&1&0&0 \cr
2&-1&0&0 \cr
-3&2&4&-2 \cr
-1&1&2&0 \cr
1&0&2&0 \cr
-4&3&4&-2
} \hskip 15 pt
\matrix{P{26} \cr
0&0&0&0 \cr
5&-3&-4&2 \cr
0&0&-2&0 \cr
-5&3&2&-2}
$$
\/}

\subsection{The Dynamical Polygons}
\label{polylist2}

Here we list the $75$ dynamical polygons associated to the
$75$ genes.  We use the same notation as above. 

{\tiny
$$\matrix{P{0} \cr
13&-8&0&0 \cr
-21&13&16&-10 \cr
13&-8&-10&6 \cr
5&-3&0&0
} \hskip 15 pt
\matrix{P{1} \cr
-15&9&8&-6 \cr
19&-12&-8&4 \cr
11&-7&-8&4 \cr
19&-12&-18&10
} \hskip 15 pt
\matrix{P{2} \cr
9&-5&-6&4 \cr
17&-10&-6&4 \cr
-17&11&10&-6 \cr
17&-10&-16&10
} \hskip 15 pt
\matrix{P{3} \cr
13&-8&-6&4 \cr
5&-3&-6&4 \cr
13&-8&-16&10 \cr
-21&13&10&-6
}$$

$$\matrix{P{4} \cr
16&-10&-8&4 \cr
3&-2&-2&0 \cr
-10&6&-2&0 \cr
24&-15&-18&10
} \hskip 15 pt
\matrix{P{5} \cr
14&-8&-4&4 \cr
1&0&2&0 \cr
-12&8&2&0 \cr
22&-13&-14&10
} \hskip 15 pt
\matrix{P{6} \cr
19&-11&-8&6 \cr
-15&10&18&-10 \cr
-7&5&8&-4 \cr
-15&10&8&-4
} \hskip 15 pt
\matrix{P{7} \cr
10&-6&-6&4 \cr
5&-3&-6&4 \cr
-3&2&4&-2
}$$

$$\matrix{P{8} \cr
-10&6&12&-8 \cr
3&-2&-4&2 \cr
16&-10&-14&8 \cr
24&-15&-14&8
} \hskip 15 pt
\matrix{P{9} \cr
-20&13&14&-8 \cr
1&0&-2&2 \cr
-12&8&14&-8
} \hskip 15 pt
\matrix{P{10} \cr
-18&11&6&-4 \cr
3&-2&-10&6 \cr
-10&6&6&-4
} \hskip 15 pt
\matrix{P{11} \cr
23&-14&-16&10 \cr
10&-6&0&0 \cr
2&-1&0&0
}$$

$$\matrix{P{12} \cr
-16&10&14&-8 \cr
5&-3&-2&2 \cr
26&-16&-12&8 \cr
5&-3&-12&8
} \hskip 15 pt
\matrix{P{13} \cr
-17&10&4&-4 \cr
-4&2&4&-4 \cr
17&-11&-12&6
} \hskip 15 pt
\matrix{P{14} \cr
-12&8&14&-8 \cr
-20&13&14&-8 \cr
14&-8&-12&8 \cr
22&-13&-12&8
} \hskip 15 pt
\matrix{P{15} \cr
-15&10&18&-10 \cr
19&-11&-8&6 \cr
-2&2&2&0
}$$

$$\matrix{P{16} \cr
-23&14&16&-10 \cr
-10&6&0&0 \cr
11&-7&-10&6
} \hskip 15 pt
\matrix{P{17} \cr
16&-10&-14&8 \cr
-5&3&2&-2 \cr
-26&16&12&-8 \cr
-5&3&12&-8
} \hskip 15 pt
\matrix{P{18} \cr
17&-10&-4&4 \cr
4&-2&-4&4 \cr
-4&3&6&-2
} \hskip 15 pt
\matrix{P{19} \cr
13&-8&0&0 \cr
-8&5&10&-6 \cr
0&0&0&0
}$$

$$\matrix{P{20} \cr
24&-15&-18&10 \cr
3&-2&-8&4 \cr
16&-10&-8&4
} \hskip 15 pt
\matrix{P{21} \cr
1&0&-4&4 \cr
22&-13&-14&10 \cr
14&-8&-4&4
} \hskip 15 pt
\matrix{P{22} \cr
-19&12&8&-4 \cr
-6&4&8&-4 \cr
2&-1&-2&2
} \hskip 15 pt
\matrix{P{23} \cr
9&-5&-2&2 \cr
-12&8&14&-8 \cr
-4&3&4&-2
}$$

$$\matrix{P{24} \cr
-5&3&4&-2 \cr
16&-10&-6&4 \cr
-5&3&-6&4 \cr
-26&16&20&-12
} \hskip 15 pt
\matrix{P{25} \cr
19&-12&-16&10 \cr
6&-4&0&0 \cr
-15&9&10&-6
} \hskip 15 pt
\matrix{P{26} \cr
17&-10&-16&10 \cr
-17&11&10&-6 \cr
4&-2&0&0
} \hskip 15 pt
\matrix{P{27} \cr
13&-8&-16&10 \cr
-21&13&10&-6 \cr
0&0&0&0
}$$

$$\matrix{P{28} \cr
13&-8&0&0 \cr
-21&13&16&-10 \cr
0&0&0&0
} \hskip 15 pt
\matrix{P{29} \cr
19&-12&-8&4 \cr
6&-4&-8&4 \cr
-15&9&8&-6
} \hskip 15 pt
\matrix{P{30} \cr
-9&5&2&-2 \cr
12&-8&-14&8 \cr
4&-3&-4&2
} \hskip 15 pt
\matrix{P{31} \cr
5&-3&-4&2 \cr
-16&10&6&-4 \cr
5&-3&6&-4 \cr
26&-16&-20&12
}$$

$$\matrix{P{32} \cr
-6&4&0&0 \cr
-19&12&16&-10 \cr
2&-1&0&0
} \hskip 15 pt
\matrix{P{33} \cr
4&-3&-2&0 \cr
-17&10&14&-10 \cr
-4&2&-2&0
} \hskip 15 pt
\matrix{P{34} \cr
24&-15&-14&8 \cr
16&-10&-14&8 \cr
3&-2&2&-2
} \hskip 15 pt
\matrix{P{35} \cr
14&-8&-12&8 \cr
22&-13&-12&8 \cr
1&0&4&-2
}$$

$$\matrix{P{36} \cr
-17&10&4&-4 \cr
-4&2&4&-4 \cr
4&-3&-6&2
} \hskip 15 pt
\matrix{P{37} \cr
-13&8&0&0 \cr
8&-5&-10&6 \cr
0&0&0&0
} \hskip 15 pt
\matrix{P{38} \cr
-16&10&8&-4 \cr
-3&2&2&0 \cr
10&-6&2&0 \cr
-24&15&18&-10
} \hskip 15 pt
\matrix{P{39} \cr
-14&8&4&-4 \cr
-1&0&-2&0 \cr
12&-8&-2&0 \cr
-22&13&14&-10
}$$

$$\matrix{P{40} \cr
-10&6&6&-4 \cr
-5&3&6&-4 \cr
3&-2&-4&2
} \hskip 15 pt
\matrix{P{41} \cr
-24&15&14&-8 \cr
-11&7&14&-8 \cr
10&-6&-12&8
} \hskip 15 pt
\matrix{P{42} \cr
26&-16&-20&12 \cr
5&-3&6&-4 \cr
-8&5&6&-4
} \hskip 15 pt
\matrix{P{43} \cr
-19&12&16&-10 \cr
-6&4&0&0 \cr
15&-9&-10&6
}$$

$$\matrix{P{44} \cr
-13&8&16&-10 \cr
21&-13&-10&6 \cr
0&0&0&0
} \hskip 15 pt
\matrix{P{45} \cr
-24&15&18&-10 \cr
-3&2&8&-4 \cr
-16&10&8&-4
} \hskip 15 pt
\matrix{P{46} \cr
-1&0&4&-4 \cr
-22&13&14&-10 \cr
-14&8&4&-4
} \hskip 15 pt
\matrix{P{47} \cr
19&-12&-8&4 \cr
6&-4&-8&4 \cr
-2&1&2&-2
}$$

$$\matrix{P{48} \cr
18&-11&-10&6 \cr
5&-3&6&-4 \cr
-16&10&6&-4
} \hskip 15 pt
\matrix{P{49} \cr
-10&6&-2&0 \cr
24&-15&-18&10 \cr
11&-7&-2&0
} \hskip 15 pt
\matrix{P{50} \cr
-12&8&2&0 \cr
22&-13&-14&10 \cr
9&-5&2&0
} \hskip 15 pt
\matrix{P{51} \cr
24&-15&-14&8 \cr
11&-7&-14&8 \cr
-10&6&12&-8
}$$

$$\matrix{P{52} \cr
9&-5&-12&8 \cr
22&-13&-12&8 \cr
-12&8&14&-8
} \hskip 15 pt
\matrix{P{53} \cr
-26&16&20&-12 \cr
-5&3&-6&4 \cr
8&-5&-6&4
} \hskip 15 pt
\matrix{P{54} \cr
6&-4&0&0 \cr
19&-12&-16&10 \cr
-2&1&0&0
} \hskip 15 pt
\matrix{P{55} \cr
-4&3&2&0 \cr
17&-10&-14&10 \cr
4&-2&2&0
}$$

$$\matrix{P{56} \cr
-24&15&14&-8 \cr
-16&10&14&-8 \cr
-3&2&-2&2
} \hskip 15 pt
\matrix{P{57} \cr
-13&8&0&0 \cr
21&-13&-16&10 \cr
0&0&0&0
} \hskip 15 pt
\matrix{P{58} \cr
-19&12&8&-4 \cr
-6&4&8&-4 \cr
15&-9&-8&6
} \hskip 15 pt
\matrix{P{59} \cr
-18&11&10&-6 \cr
-5&3&-6&4 \cr
16&-10&-6&4
}$$

$$\matrix{P{60} \cr
10&-6&2&0 \cr
-24&15&18&-10 \cr
-11&7&2&0
} \hskip 15 pt
\matrix{P{61} \cr
12&-8&-2&0 \cr
-22&13&14&-10 \cr
-9&5&-2&0
} \hskip 15 pt
\matrix{P{62} \cr
10&-6&-12&8 \cr
-3&2&4&-2 \cr
-16&10&14&-8 \cr
-24&15&14&-8
} \hskip 15 pt
\matrix{P{63} \cr
18&-11&-6&4 \cr
-3&2&10&-6 \cr
10&-6&-6&4
}$$

$$\matrix{P{64} \cr
-23&14&16&-10 \cr
-10&6&0&0 \cr
-2&1&0&0
} \hskip 15 pt
\matrix{P{65} \cr
-13&8&0&0 \cr
21&-13&-16&10 \cr
-13&8&10&-6 \cr
-5&3&0&0
} \hskip 15 pt
\matrix{P{66} \cr
15&-9&-8&6 \cr
-19&12&8&-4 \cr
-11&7&8&-4 \cr
-19&12&18&-10
} \hskip 15 pt
\matrix{P{67} \cr
-13&8&6&-4 \cr
-5&3&6&-4 \cr
-13&8&16&-10 \cr
21&-13&-10&6
}$$

$$\matrix{P{68} \cr
-26&16&12&-8 \cr
-5&3&12&-8 \cr
8&-5&-4&2
} \hskip 15 pt
\matrix{P{69} \cr
17&-10&-4&4 \cr
4&-2&-4&4 \cr
-17&11&12&-6
} \hskip 15 pt
\matrix{P{70} \cr
-5&3&12&-8 \cr
-18&11&12&-8 \cr
16&-10&-14&8
} \hskip 15 pt
\matrix{P{71} \cr
14&-8&-12&8 \cr
-20&13&14&-8 \cr
-7&5&14&-8
}$$

$$\matrix{P{72} \cr
26&-16&-12&8 \cr
5&-3&-12&8 \cr
-8&5&4&-2
} \hskip 15 pt
\matrix{P{73} \cr
23&-14&-16&10 \cr
10&-6&0&0 \cr
-11&7&10&-6
} \hskip 15 pt
\matrix{P{74} \cr
5&-3&-12&8 \cr
18&-11&-12&8 \cr
-16&10&14&-8
}$$
\/}

\subsection{The Gene Locations}
\label{centerlist}

Here we list the coordinates for the
points $a_0,...,a_{74}$.  The point
$a_j$ is the central vertex of the
gene $A_j$.

$$\matrix{a0\cr 3 \cr 4 } \hskip 15 pt
\matrix{a1\cr 4 \cr 5 } \hskip 15 pt
\matrix{a2\cr 4 \cr 6 } \hskip 15 pt
\matrix{a3\cr 4 \cr 7 } \hskip 15 pt
\matrix{a4\cr 4 \cr 10 } \hskip 15 pt
\matrix{a5\cr 4 \cr 11 } \hskip 15 pt
\matrix{a6\cr 3 \cr 13 } \hskip 15 pt
\matrix{a7\cr 5 \cr 15 } \hskip 15 pt
\matrix{a8\cr 6 \cr 16 } \hskip 15 pt
\matrix{a9\cr 6 \cr 17 }$$

$$\matrix{a10\cr 7 \cr 19 } \hskip 15 pt
\matrix{a11\cr 8 \cr 16 } \hskip 15 pt
\matrix{a12\cr 10 \cr 17 } \hskip 15 pt
\matrix{a13\cr 11 \cr 16 } \hskip 15 pt
\matrix{a14\cr 14 \cr 20 } \hskip 15 pt
\matrix{a15\cr 16 \cr 31 } \hskip 15 pt
\matrix{a16\cr 16 \cr 33 } \hskip 15 pt
\matrix{a17\cr 14 \cr 32 } \hskip 15 pt
\matrix{a18\cr 13 \cr 33 } \hskip 15 pt
\matrix{a19\cr 12 \cr 31 }$$

$$\matrix{a20\cr 17 \cr 49 } \hskip 15 pt
\matrix{a21\cr 17 \cr 50 } \hskip 15 pt
\matrix{a22\cr 15 \cr 50 } \hskip 15 pt
\matrix{a23\cr 14 \cr 49 } \hskip 15 pt
\matrix{a24\cr 13 \cr 50 } \hskip 15 pt
\matrix{a15\cr 12 \cr 50 } \hskip 15 pt
\matrix{a26\cr 12 \cr 51 } \hskip 15 pt
\matrix{a27\cr 12 \cr 52 } \hskip 15 pt
\matrix{a28\cr 16 \cr 64 } \hskip 15 pt
\matrix{a29\cr 17 \cr 65 }$$

$$\matrix{a30\cr 18 \cr 66 } \hskip 15 pt
\matrix{a31\cr 19 \cr 65 } \hskip 15 pt
\matrix{a32\cr 20 \cr 65 } \hskip 15 pt
\matrix{a33\cr 20 \cr 64 } \hskip 15 pt
\matrix{a34\cr 22 \cr 64 } \hskip 15 pt
\matrix{a35\cr 22 \cr 65 } \hskip 15 pt
\matrix{a36\cr 32 \cr 79 } \hskip 15 pt
\matrix{a37\cr 33 \cr 81 } \hskip 15 pt
\matrix{a38\cr 33 \cr 78 } \hskip 15 pt
\matrix{a39\cr 33 \cr 77 }$$

$$\matrix{a40\cr 32 \cr 73 } \hskip 15 pt
\matrix{a41\cr 31 \cr 72 } \hskip 15 pt
\matrix{a42\cr 32 \cr 70 } \hskip 15 pt
\matrix{a43\cr 41 \cr 73 } \hskip 15 pt
\matrix{a44\cr 41 \cr 71 } \hskip 15 pt
\matrix{a45\cr 41 \cr 68 } \hskip 15 pt
\matrix{a46\cr 41 \cr 67 } \hskip 15 pt
\matrix{a47\cr 43 \cr 67 } \hskip 15 pt
\matrix{a48\cr 45 \cr 67 } \hskip 15 pt
\matrix{a49\cr 46 \cr 68 }$$

$$\matrix{a50\cr 46 \cr 69 } \hskip 15 pt
\matrix{a51\cr 69 \cr 137 } \hskip 15 pt
\matrix{a52\cr 69 \cr 138 } \hskip 15 pt
\matrix{a53\cr 68 \cr 139 } \hskip 15 pt
\matrix{a54\cr 67 \cr 139 } \hskip 15 pt
\matrix{a55\cr 67 \cr 140 } \hskip 15 pt
\matrix{a56\cr 65 \cr 140 } \hskip 15 pt
\matrix{a57\cr 63 \cr 137 } \hskip 15 pt
\matrix{a58\cr 62 \cr 136 } \hskip 15 pt
\matrix{a59\cr 55 \cr 142 }$$

$$\matrix{a60\cr 54 \cr 141 } \hskip 15 pt
\matrix{a61\cr 54 \cr 140 } \hskip 15 pt
\matrix{a62\cr 52 \cr 135 } \hskip 15 pt
\matrix{a63\cr 51 \cr 132 } \hskip 15 pt
\matrix{a64\cr 50 \cr 135 } \hskip 15 pt
\matrix{a65\cr 63 \cr 213 } \hskip 15 pt
\matrix{a66\cr 62 \cr 212 } \hskip 15 pt
\matrix{a67\cr 62 \cr 210 } \hskip 15 pt
\matrix{a68\cr 61 \cr 207 } \hskip 15 pt
\matrix{a69\cr 55 \cr 214 }$$

$$\matrix{a70\cr 48 \cr 210 } \hskip 15 pt
\matrix{a71\cr 48 \cr 211 } \hskip 15 pt
\matrix{a72\cr 73 \cr 282 } \hskip 15 pt
\matrix{a73\cr 76 \cr 275 } \hskip 15 pt
\matrix{a74\cr 86 \cr 279 } $$

\subsection{The Inflation Data}
\label{shadowlist}

Here we list the data for the strands $A'_0,...,A'_{74}$ that
are associated to the genes $A_0,...,A_{74}$.   The listing
$$\matrix{A_j \cr x_0 \cr y_0 \cr x_1 \cr y_1 \cr x_2 \cr y_2}$$
indicates that the two endpoints of $A_j$ are
$(x_0,y_0)$ and $(x_2,y_2)$, and the point $(x_1,y_1)$ is the
one lying within $3$ units from the center point of
$D(A_j)$.   (There might be several such points, but we
make some choice in each case.)
Here is the listing.

{\tiny
$$\matrix{A_{0} \cr 9 \cr 17 \cr 11 \cr 17 \cr 15 \cr 22 } \hskip 13 pt
\matrix{A_{1} \cr 11 \cr 16 \cr 16 \cr 22 \cr 17 \cr 25 } \hskip 13 pt
\matrix{A_{2} \cr 15 \cr 22 \cr 17 \cr 25 \cr 16 \cr 30 } \hskip 13 pt
\matrix{A_{3} \cr 17 \cr 25 \cr 16 \cr 30 \cr 12 \cr 34 } \hskip 13 pt
\matrix{A_{4} \cr 12 \cr 39 \cr 16 \cr 43 \cr 17 \cr 46 } \hskip 13 pt
\matrix{A_{5} \cr 15 \cr 43 \cr 17 \cr 46 \cr 12 \cr 50 } \hskip 13 pt
\matrix{A_{6} \cr 12 \cr 50 \cr 12 \cr 56 \cr 12 \cr 60 } \hskip 13 pt
\matrix{A_{7} \cr 15 \cr 64 \cr 21 \cr 63 \cr 25 \cr 68 } \hskip 13 pt
\matrix{A_{8} \cr 21 \cr 63 \cr 25 \cr 69 \cr 24 \cr 71 } \hskip 13 pt
\matrix{A_{9} \cr 25 \cr 68 \cr 24 \cr 72 \cr 24 \cr 76 } \hskip 13 pt
\matrix{A_{10} \cr 24 \cr 76 \cr 31 \cr 80 \cr 34 \cr 76 } \hskip 13 pt
\matrix{A_{11} \cr 33 \cr 73 \cr 34 \cr 68 \cr 37 \cr 72 } \hskip 13 pt
\matrix{A_{12} \cr 37 \cr 72 \cr 41 \cr 72 \cr 41 \cr 67 }$$
$$\matrix{A_{13} \cr 41 \cr 67 \cr 46 \cr 69 \cr 45 \cr 71 } \hskip 13 pt
\matrix{A_{14} \cr 59 \cr 81 \cr 58 \cr 85 \cr 58 \cr 89 } \hskip 13 pt
\matrix{A_{15} \cr 67 \cr 128 \cr 67 \cr 132 \cr 67 \cr 136 } \hskip 13 pt
\matrix{A_{16} \cr 67 \cr 136 \cr 67 \cr 140 \cr 62 \cr 135 } \hskip 13 pt
\matrix{A_{17} \cr 62 \cr 136 \cr 59 \cr 137 \cr 58 \cr 140 } \hskip 13 pt
\matrix{A_{18} \cr 58 \cr 140 \cr 54 \cr 140 \cr 54 \cr 136 } \hskip 13 pt
\matrix{A_{19} \cr 54 \cr 136 \cr 51 \cr 132 \cr 46 \cr 136 } \hskip 13 pt
\matrix{A_{20} \cr 67 \cr 204 \cr 71 \cr 208 \cr 71 \cr 212 } \hskip 13 pt
\matrix{A_{21} \cr 71 \cr 208 \cr 71 \cr 213 \cr 66 \cr 217 } \hskip 13 pt
\matrix{A_{22} \cr 66 \cr 217 \cr 62 \cr 211 \cr 63 \cr 208 } \hskip 13 pt
\matrix{A_{23} \cr 63 \cr 208 \cr 58 \cr 208 \cr 58 \cr 212 } \hskip 13 pt
\matrix{A_{24} \cr 58 \cr 212 \cr 54 \cr 212 \cr 51 \cr 208 } \hskip 13 pt
\matrix{A_{25} \cr 51 \cr 208 \cr 50 \cr 211 \cr 51 \cr 216 }$$
$$\matrix{A_{26} \cr 52 \cr 211 \cr 51 \cr 216 \cr 50 \cr 221 } \hskip 13 pt
\matrix{A_{27} \cr 51 \cr 216 \cr 50 \cr 221 \cr 46 \cr 225 } \hskip 13 pt
\matrix{A_{28} \cr 64 \cr 271 \cr 68 \cr 272 \cr 71 \cr 276 } \hskip 13 pt
\matrix{A_{29} \cr 67 \cr 272 \cr 71 \cr 276 \cr 72 \cr 279 } \hskip 13 pt
\matrix{A_{30} \cr 72 \cr 279 \cr 75 \cr 280 \cr 76 \cr 275 } \hskip 13 pt
\matrix{A_{31} \cr 78 \cr 275 \cr 79 \cr 276 \cr 84 \cr 279 } \hskip 13 pt
\matrix{A_{32} \cr 85 \cr 279 \cr 84 \cr 276 \cr 84 \cr 271 } \hskip 13 pt
\matrix{A_{33} \cr 85 \cr 276 \cr 86 \cr 270 \cr 89 \cr 267 } \hskip 13 pt
\matrix{A_{34} \cr 89 \cr 267 \cr 92 \cr 271 \cr 92 \cr 275 } \hskip 13 pt
\matrix{A_{35} \cr 92 \cr 271 \cr 92 \cr 276 \cr 92 \cr 280 } \hskip 13 pt
\matrix{A_{36} \cr 130 \cr 334 \cr 135 \cr 336 \cr 134 \cr 338 } \hskip 13 pt
\matrix{A_{37} \cr 134 \cr 338 \cr 138 \cr 344 \cr 143 \cr 340 } \hskip 13 pt
\matrix{A_{38} \cr 143 \cr 335 \cr 139 \cr 331 \cr 139 \cr 326 }$$
$$\matrix{A_{39} \cr 140 \cr 331 \cr 141 \cr 325 \cr 143 \cr 322 } \hskip 13 pt
\matrix{A_{40} \cr 140 \cr 310 \cr 135 \cr 310 \cr 130 \cr 305 } \hskip 13 pt
\matrix{A_{41} \cr 135 \cr 309 \cr 130 \cr 305 \cr 130 \cr 301 } \hskip 13 pt
\matrix{A_{42} \cr 132 \cr 297 \cr 134 \cr 297 \cr 139 \cr 300 } \hskip 13 pt
\matrix{A_{43} \cr 174 \cr 313 \cr 173 \cr 310 \cr 173 \cr 305 } \hskip 13 pt
\matrix{A_{44} \cr 173 \cr 305 \cr 175 \cr 300 \cr 177 \cr 297 } \hskip 13 pt
\matrix{A_{45} \cr 177 \cr 293 \cr 172 \cr 287 \cr 173 \cr 284 } \hskip 13 pt
\matrix{A_{46} \cr 172 \cr 288 \cr 175 \cr 283 \cr 178 \cr 280 } \hskip 13 pt
\matrix{A_{47} \cr 178 \cr 280 \cr 181 \cr 284 \cr 181 \cr 288 } \hskip 13 pt
\matrix{A_{48} \cr 187 \cr 284 \cr 189 \cr 284 \cr 193 \cr 289 } \hskip 13 pt
\matrix{A_{49} \cr 189 \cr 283 \cr 194 \cr 289 \cr 195 \cr 292 } \hskip 13 pt
\matrix{A_{50} \cr 193 \cr 289 \cr 195 \cr 292 \cr 190 \cr 296 } \hskip 13 pt
\matrix{A_{51} \cr 288 \cr 576 \cr 291 \cr 580 \cr 291 \cr 584 }$$
$$\matrix{A_{52} \cr 291 \cr 580 \cr 291 \cr 585 \cr 291 \cr 589 } \hskip 13 pt
\matrix{A_{53} \cr 291 \cr 589 \cr 287 \cr 589 \cr 284 \cr 585 } \hskip 13 pt
\matrix{A_{54} \cr 284 \cr 585 \cr 283 \cr 588 \cr 284 \cr 593 } \hskip 13 pt
\matrix{A_{55} \cr 285 \cr 588 \cr 284 \cr 593 \cr 279 \cr 597 } \hskip 13 pt
\matrix{A_{56} \cr 279 \cr 597 \cr 274 \cr 593 \cr 274 \cr 589 } \hskip 13 pt
\matrix{A_{57} \cr 271 \cr 580 \cr 265 \cr 581 \cr 261 \cr 575 } \hskip 13 pt
\matrix{A_{58} \cr 266 \cr 581 \cr 261 \cr 575 \cr 262 \cr 572 } \hskip 13 pt
\matrix{A_{59} \cr 237 \cr 600 \cr 231 \cr 602 \cr 228 \cr 598 } \hskip 13 pt
\matrix{A_{60} \cr 232 \cr 602 \cr 228 \cr 598 \cr 228 \cr 593 } \hskip 13 pt
\matrix{A_{61} \cr 229 \cr 598 \cr 230 \cr 592 \cr 232 \cr 589 } \hskip 13 pt
\matrix{A_{62} \cr 224 \cr 576 \cr 219 \cr 572 \cr 219 \cr 568 } \hskip 13 pt
\matrix{A_{63} \cr 219 \cr 563 \cr 216 \cr 559 \cr 211 \cr 563 } \hskip 13 pt
\matrix{A_{64} \cr 211 \cr 568 \cr 211 \cr 572 \cr 206 \cr 567 }$$
$$\matrix{A_{65} \cr 271 \cr 902 \cr 265 \cr 903 \cr 261 \cr 897 } \hskip 13 pt
\matrix{A_{66} \cr 266 \cr 903 \cr 261 \cr 897 \cr 262 \cr 894 } \hskip 13 pt
\matrix{A_{67} \cr 262 \cr 894 \cr 264 \cr 889 \cr 267 \cr 885 } \hskip 13 pt
\matrix{A_{68} \cr 261 \cr 877 \cr 258 \cr 878 \cr 257 \cr 881 } \hskip 13 pt
\matrix{A_{69} \cr 236 \cr 907 \cr 233 \cr 906 \cr 231 \cr 903 } \hskip 13 pt
\matrix{A_{70} \cr 206 \cr 890 \cr 203 \cr 891 \cr 202 \cr 894 } \hskip 13 pt
\matrix{A_{71} \cr 203 \cr 890 \cr 202 \cr 894 \cr 202 \cr 898 } \hskip 13 pt
\matrix{A_{72} \cr 303 \cr 1195 \cr 308 \cr 1195 \cr 308 \cr 1190 } \hskip 13 pt
\matrix{A_{73} \cr 321 \cr 1170 \cr 322 \cr 1165 \cr 325 \cr 1169 } \hskip 13 pt
\matrix{A_{74} \cr 359 \cr 1182 \cr 363 \cr 1182 \cr 363 \cr 1177 }$$
\/}

\newpage

%% file: refs.tex
\section{References}

[{\bf B\/}] P. Boyland, {\it Dual Billiards, twist maps, and impact oscillators\/},
Nonlinearity {\bf 9\/} (1996) 1411-1438
\newline
\newline
[{\bf D\/}], R. Douady, {\it These de 3-eme cycle\/}, Universite de Paris 7, 1982
\newline
\newline
[{\bf DT\/}] F. Dogru and S. Tabachnikov, {\it Dual Billiards\/},
Math Intelligencer vol. 27 No. 4 (2005) 18--25
\newline
\newline
[{\bf G\/}] D. Genin, {\it Regular and Chaotic Dynamics of
Outer Billiards\/}, Penn State Ph.D. thesis (2005)
 \newline \newline
[{\bf GS\/}] E. Gutkin and N. Simanyi, {\it Dual polygonal
billiard and necklace dynamics\/}, Comm. Math. Phys.
{\bf 143\/} (1991) 431--450
\newline
\newline
[{\bf Ke\/}] R. Kenyon, {\it Inflationary tilings with a similarity
structure\/}, Comment. Math. Helv. {\bf 69\/} (1994) 169--198
\newline
\newline
[{\bf Ko\/}] Kolodziej, {\it The antibilliard outside a polygon\/},
Bull. Polish Acad Sci. Math.
{\bf 37\/} (1989) 163--168
\newline
\newline
[{\bf M\/}] J. Moser, {\it Stable and Random Motions in Dynamical Systems, with
Special Emphasis on Celestial Mechanics\/},
Annals of Math Studies {\bf 77\/}, Princeton University Press (1973)
\newline
\newline
[{\bf N\/}] B.H. Neumann, {\it Sharing Ham and Eggs\/},
\newline
summary of a Manchester Mathematics Colloquium, 25 Jan 1959 
\newline
published in Iota, the Manchester University Mathematics students' journal
\newline
\newline
[{\bf T\/}] S. Tabachnikov, {\it Geometry and Billiards\/},
A.M.S. Mathematics Advanced Study Semesters (2005)
\newline
\newline
[{\bf VS\/}] F. Vivaldi, A. Shaidenko, {\it Global stability of a class of discontinuous
dual billiards\/}, Comm. Math. Phys. {\bf 110\/} (1987) 625--640 
\newline
\newline